\documentclass[11pt]{article}
\usepackage{amssymb}
\usepackage{amsmath}
\usepackage{amsthm}
\newcommand{\cA}{{\cal A}}
\newcommand{\cB}{{\cal B}}
\newcommand{\cC}{{\cal C}}

\newcommand{\cH}{{\cal H}}

\newcommand{\cI}{{\cal I}}

\newcommand{\cO}{{\cal O}}
\newcommand{\cL}{{\cal L}}
\newcommand{\cM}{{\cal M}}

\newcommand{\cF}{{\cal F}}
\newcommand{\cK}{{\cal K}}
\newcommand{\cP}{{\cal P}}

\newcommand{\cS}{{\cal S}}
\newcommand{\cT}{{\cal T}}
\newcommand{\cU}{{\cal U}}
\newcommand{\cV}{{\cal V}}
\newcommand{\cW}{{\cal W}}

\newcommand{\cY}{{\cal Y}}

\renewcommand{\AA}{{\mathbb A}}

\newcommand{\NN}{{\mathbb N}}
\newcommand{\ZZ}{{\mathbb Z}}

\newcommand{\OO}{{\mathbb O}}

\newcommand{\gp}{\mathfrak{p}}
\newcommand{\gq}{\mathfrak{q}}

\newcommand{\gt}{\mathfrak{t}}
\newcommand{\gr}{\mathfrak{r}}

\newcommand{\on}{\operatorname}

\newcommand{\Rep}{{\on{Rep}}}
\newcommand{\Sch}{{\on{Sch}}}

\newcommand{\Qlb}{\mathbb{\bar Q}_\ell}
\newcommand{\Gm}{\mathbb{G}_m}

\newcommand{\A}{\mathbb{A}}

\newcommand{\toup}[1]{\stackrel{#1}{\to}}
\newcommand{\hook}[1]{\stackrel{#1}{\hookrightarrow}}

\newcommand{\getsup}[1]{\stackrel{#1}{\gets}}
\newcommand{\Sp}{\on{\mathbb{S}p}}

\newcommand{\IC}{\on{IC}}

\newcommand{\Hom}{\on{Hom}}

\newcommand{\Ext}{\on{Ext}}

\newcommand{\Sym}{\on{Sym}}
\newcommand{\SO}{\on{S\mathbb{O}}}

\newcommand{\Aut}{\on{Aut}}

\newcommand{\RG}{\on{R\Gamma}}

\newcommand{\const}{{const}}

\newcommand{\Pic}{\on{Pic}}

\newcommand{\Bun}{\on{Bun}}

\newcommand{\Bunt}{\on{\widetilde\Bun}}

\newcommand{\pro}{\on{pro}}
\newcommand{\Spec}{\on{Spec}}

\newcommand{\supp}{\on{supp}}

\newcommand{\Gr}{\on{Gr}}

\newcommand{\GL}{\on{GL}}

\newcommand{\pr}{\on{pr}}
\newcommand{\id}{\on{id}}
\newcommand{\tr}{\on{tr}}
\newcommand{\QED}{$\square$} 
\newcommand{\Fq}{\mathbb{F}_q}  
\newcommand{\Fp}{\mathbb{F}_p}  
\newcommand{\iso}{{\widetilde\to}}

\newcommand{\comp}{\circ}
\newcommand{\Four}{\on{Four}}
\renewcommand{\H}{{\on{H}}}   

\newcommand{\DD}{\mathbb{D}}  
\newcommand{\D}{\on{D}}       
\newcommand{\DP}{\on{DP}}
\newcommand{\wt}{\widetilde}
\newcommand{\ov}[1]{\overline{#1}}
\newcommand{\select}[1]{{\it{#1}}}

\newcommand{\und}[1]{\underline{#1}}

\renewcommand{\P}{{\on{P}}}

\newcommand{\<}{\langle}
\renewcommand{\>}{\rangle}

\newcommand{\ev}{\mathit{ev}}

\newcommand{\Loc}{\on{Loc}}

\newcommand{\Sph}{\on{Sph}}
\newcommand{\Res}{\on{Res}}
\newcommand{\gRes}{\on{gRes}}
\newcommand{\ugRes}{\on{g\und{Res}}}

\newcommand{\ttimes}{\tilde\times}

\newcommand{\act}{\on{act}}
\newcommand{\dimrel}{\on{dim.rel}}

\newcommand{\Funct}{\on{Funct}}

\newcommand{\SL}{\on{SL}}

\newcommand{\tboxtimes}{\,\tilde\boxtimes\,}

\newcommand{\glob}{\on{glob}}
\newcommand{\Ind}{\on{Ind}}

\newcommand{\ra}{\rightarrow}
\newcommand{\la}{\leftarrow}

\newcommand{\Heis}{\on{Hs}}

\newcommand{\Weil}{\on{{\mathcal{W}}eil}}
\newcommand{\Char}{\on{Char}}
\newcommand{\Cr}{\on{Cr}}

\newcommand{\temp}{\on{temp}}

\newcommand{\LW}{\on{LW}}
\newcommand{\colim}{\on{colim}}
\newcommand{\DGCat}{\on{DGCat}}
\newcommand{\DG}{\on{DG}}
\newcommand{\Spc}{\on{Spc}}
\newcommand{\PreStk}{\on{PreStk}}
\newcommand{\Cat}{\on{Cat}}

\newtheorem{Lm}{Lemma}
\newtheorem{Th}{Theorem}
\newtheorem{Pp}{Proposition}
\newtheorem{Cor}{Corollary}
\newtheorem{Con}{Conjecture}

\theoremstyle{remark}
\newtheorem{Rem}{Remark}

\newtheorem{Ex}{Example}
\theoremstyle{definition}
\newtheorem{Def}{Definition}

\newenvironment{Prf}{\par\noindent {\it Proof }}{\QED}
\newcommand{\Step}[1]{\par\noindent{\bf Step {#1}}.}

\begin{document}

\author{Sergey Lysenko}
\title{Geometric theta-lifting for the dual pair $\SO_{2m}, \Sp_{2n}$}
\date{}
\maketitle
\begin{abstract}
\noindent{\scshape Abstract}\hskip 0.8 em 
Let $X$ be a smooth projective curve over an algebraically closed field of characteristic $>2$. Consider the dual pair $H=\SO_{2m}, G=\Sp_{2n}$ over $X$ with $H$ split. Write $\Bun_G$ and $\Bun_H$ for the stacks of $G$-torsors and $H$-torsors on $X$. The theta-kernel $\Aut_{G,H}$ on $\Bun_G\times\Bun_H$ yields theta-lifting functors $F_G: \D(\Bun_H)\to\D(\Bun_G)$ and $F_H: \D(\Bun_G)\to\D(\Bun_H)$ between the corresponding derived categories. We describe the relation of these functors with Hecke operators. 

 In two particular cases these functors realize the geometric Langlands functoriality for the above pair (in the non ramified case). Namely, we show that for $n=m$ the functor $F_G: \D(\Bun_H)\to\D(\Bun_G)$ commutes with Hecke operators with respect to the inclusion of the Langlands dual groups $\check{H}\,\iso\, \SO_{2n}\hook{} \SO_{2n+1}\,\iso\,\check{G}$. For $m=n+1$ we show that the functor $F_H: \D(\Bun_G)\to\D(\Bun_H)$ commutes with Hecke operators with respect to the inclusion of the Langlands dual groups $\check{G}\,\iso\, \SO_{2n+1}\hook{} \SO_{2n+2}\,\iso\, \check{H}$.

  In other cases the relation is more complicated and involves the $\SL_2$ of Arthur.  As a step of the proof, we establish the geometric theta-lifting for the dual pair $\GL_m, \GL_n$. Our global results are derived from the corresponding local ones, which provide a geometric analog of a theorem of Rallis.
\end{abstract} 

\medskip

{\centerline{\scshape 1. Introduction}}

\bigskip\noindent
0.1 This is a corrected version of \cite{Ly2}. In the published version the compatibility of the functoriality isomorphisms (in main theorems) with the tensor structure on the category of representations of the dual group was not justified. We do this in the present version. We also pass from the formalism of derived categories to that of $\DG$-categories. 

\medskip\noindent
1.1 The Howe correspondence for dual reductive pairs is known to realize the Langlands functoriality in some particular cases (cf. \cite{R}, \cite{Ad}, \cite{Ku}). In this paper, which is a continuation of \cite{L2}, we depelop a similar geometric theory for the dual reductive pairs $(\Sp_{2n}, \SO_{2m})$ and $(\GL_n,\GL_m)$. We consider only the everywhere unramified case. 

 Recall the classical construction of the theta-lifting operators. Let $X$ be a smooth projective geometrically connected curve over $k=\Fq$. Let $F=k(X)$, $\AA$ be the ad\`eles ring of $X$, $\cO$ be the integer ad\`eles. Let $G$, $H$ be split connected reductive groups over $\Fq$ that form a dual pair inside some symplectic group $\Sp_{2r}$. Assume that the metaplectic covering $\wt\Sp_{2r}(\AA)\to\Sp_{2r}(\AA)$ splits over $G(\AA)\times H(\AA)$. Let $S$ be the corresponding Weil representation of $G(\AA)\times H(\AA)$. A choice of a theta-functional $\theta:S\to\Qlb$ yields a morphism of modules over the global non ramified Hecke algebras $\cH_G\otimes\cH_H$
$$
S^{(G\times H)(\cO)}\to \Funct((G\times H)(F)\backslash (G\times H)(\AA)/(G\times H)(\cO))
$$
sending $\phi$ to the function $(g,h)\mapsto \theta((g,h)\phi)$. The space  $S^{(G\times H)(\cO)}$ has a distinguished non ramified vector, its image $\phi_0$ under the above map is the classical theta-function. Viewing $\phi_0$ as a kernel of integral operators, one gets the classical theta-lifting operators 
$$
F_G: \Funct(H(F)\backslash H(\AA)/H(\cO))\to \Funct(G(F)\backslash G(\AA)/G(\cO))
$$
and 
$$
F_H: \Funct(G(F)\backslash G(\AA)/G(\cO))\to
\Funct(H(F)\backslash H(\AA)/H(\cO))
$$
For the dual pairs $(\Sp_{2n}, \SO_{2m})$ and $(\GL_n,\GL_m)$ these operators realize the Langlands functoriality between the corresponding automorphic representations (as we will see below, its precise formulation involves the $\SL_2$ of Arthur). We establish a geometric analog of this phenomenon. 
 
  Recall that $S\,\iso\, \otimes'_{x\in X} S_x$ is the restricted tensor product of local Weil representations of $G(F_x)\times H(F_x)$. Here $F_x$ denotes the completion of $F$ at $x\in X$. The above functoriality in the classical case is a consequence of a local result describing the space of invariants $S_x^{G(\cO_x)\times H(\cO_x)}$ as a module over the tensor product $_x\cH_G\otimes {_x\cH_H}$ of local (non ramified) Hecke algebras. In the geomeric setting the main step is also to prove a local analog of this and then derive the global functoriality.  The proof of this local result due to Rallis (\cite{R}) does not geometrise in an obvious way, as it makes essential use of functions with infinite-dimensional support. Their geometric counterparts should be perverse sheaves, however the notion of a perverse sheaf with  infinite-dimensional support is not known. We get around this difficulty using inductive systems of perverse sheaves rather then perverse sheaves themselves.  

 Let us underline the following phenomenon in the proof that we find striking. Let $G=\Sp_{2n}$, $H=\SO_{2m}$. The Langlands dual groups are $\check{G}\,\iso\, \SO_{2n+1}$ and $\check{H}\,\iso\, \SO_{2m}$ over $\Qlb$. Write $\Rep(\check{G})$ for the category of finite-dimensional representations of $\check{G}$ over $\Qlb$, and similarly for $\check{H}$. There will be ind-schemes $Y_H, Y_G$ over $k$ and fully faithful functors $f_H: \Rep(\check{H})\to \P(Y_H)$ and $f_G:\Rep(\check{G})\to P(Y_G)$ taking values in the categories of perverse sheaves (pure of weight zero) on $Y_H$ (resp., $Y_G$) with the following properties. Extend $f_H$ to a functor
$$
f_H: \Rep(\check{H}\times\Gm)\to \oplus_{i\in\ZZ} \, P(Y_H)[i]\subset \D(Y_H)
$$
as follows. If $V$ is a representation of $\check{H}$ and $I$ is the standard representation of $\Gm$ then 
$f_H(V\boxtimes (I^{\otimes i}))\,\iso\, f_H(V)[i]$ is placed in perverse cohomological degree $-i$. 
For $n\ge m$ there will be a ind-proper map $\pi:Y_G\to Y_H$ such that the following diagram is commutative
$$
\begin{array}{ccc}
\Rep(\check{G}) & \toup{f_G} & P(Y_G)\\
\downarrow\lefteqn{\scriptstyle \Res^{\kappa}} && \downarrow\lefteqn{\scriptstyle \pi_!}\\
\Rep(\check{H}\times\Gm) & \toup{f_H} & \oplus_{i\in\ZZ} \, P(Y_H)[i]
\end{array}
$$
for some homomorphism $\kappa: \check{H}\times\Gm\to \check{G}$. For $n=m$ the restriction of $\kappa$ to $\Gm$ is trivial, so $\pi_!f_G$ takes values in the category of perverse sheaves in this case. Both $f_G$ and $f_H$ send an irreducible representation to an irreducible perverse sheaf. So, for $V\in\Rep(\check{G})$ the decomposition of $\Res^{\kappa}(V)$ into irreducible ones can be seen via the decomposition theorem of Beilinson, Bernstein and Deligne (\cite{BBD}). Actually here $\pi: \Pi(K)\times\Gr_G\to \Pi(K)$ is the projection, where $K=k((t))$, $\Pi$ is a finite-dimensional $k$-vector space, and $\Gr_G$ is the affine grassmanian for $G$. There will also be an analog of the above result for $n<m$ (and also for the dual pair $\GL_n,\GL_m$).

 The above phenomenon is a part of our main local results (Proposition~4 in Section~5.1, Theorem~7 in Section~6.2). 
They provide a geometric analog of the local theta correspondence for these dual pairs. The key technical tools in the proof are \select{the weak geometric analogs of the Jacquet functors} (cf. Section~4.7). 

\medskip\noindent
1.2 In the global setting let $\Omega$ denote the canonical line bundle on $X$. Let $G$ be the group scheme over $X$ of automorphisms of $\cO_X^n\oplus\Omega^n$ preserving the natural symplectic form $\wedge^2(\cO_X^n\oplus\Omega^n)\to\Omega$. Let $H=\SO_{2m}$. Write $\Bun_H$ for the stack of $H$-torsors on $X$, similarly for $G$. Using the construction from \cite{L1}, we introduce a geometric analog $\Aut_{G,H}$ of the above function $\phi_0$, this is an object of the $\DG$-category of $\ell$-adic sheaves on $\Bun_G\times\Bun_H$. It yields the theta-lifting functors
$$
F_G:\D(\Bun_H)\to\D(\Bun_G)
$$
and 
$$
F_H:\D(\Bun_G)\to\D(\Bun_H)
$$
Our mains global results for the pair $(G,H)$ are Theorems~\ref{Th_main_global_symplectic_orthogonal} and \ref{Th_Hecke_property_Aut_GH} describing the relation between the theta-lifting functors and the Hecke functors on $\Bun_G$ and $\Bun_H$. They agree with the conjectures of Adams (\cite{Ad}). 
One of the advantages of the geometric setting compared to the classical one is that the $\SL_2$ of Arthur appears naturally. 

An essential difficulty in the proof was the fact that the complex $\Aut_{G,H}$ is not perverse (it has infinitely many perverse cohomologies), it is not even a direct sum of its perverse cohomologies (cf. Section~8.3). 

 We also establish the global theta-lifting for the dual pair $(\GL_n,\GL_m)$ (cf. Theorem~\ref{Th_main_global_GL_m_GL_n}).  

\medskip\noindent
1.3 Let us briefly discuss how the paper is organized. Our main results are collected in Section~2. In Section~3 we remind some constructions at the level of functions, which we have in mind for geometrization. 

 In Section~4 we develop a geometric theory for the following classical objects. Let $K=\Fq((t))$ and $\cO=\Fq[[t]]$. Given a reductive group $G$ over $\Fq$ and a finite dimensional representation $M$, the space of invariants in the Schwarz space $\cS(M(K))^{G(\cO)}$ is a module over the non ramified Hecke algebra $\cH_G$. We introduce the geometric analogs of the Fourier transform on this space and (some weak analogs) of the Jacquet functors. A way to relate this with the global case is proposed in Section~4.6. 
 
 In Section~5 we develop the local theta correspondence for the dual pair $(\GL_n,\GL_m)$. The key ingredients here are decomposition theorem from \cite{BBD}, the dimension estimates from \cite{MV} and hyperbolic localization results from \cite{B}. 
 
 In Section~6 we develop the local theta correspondence for the dual pair $(\Sp_{2n},\SO_{2m})$. In addition to the above tools, we use the classical result (Proposition~\ref{Pp_local_main_classical}) in the proof of our Proposition~\ref{Pp_functor_cP_injection}.
  
 In Section~7 we derive the global theta-lifting results for the dual pair $(\GL_n,\GL_m)$. 
 
 In Section~8 we prove our main global results (Theorems~\ref{Th_main_global_symplectic_orthogonal} and \ref{Th_Hecke_property_Aut_GH}) about theta-lifting for the dual pair $(\Sp_{2n},\SO_{2m})$. The relation between the local theory and the theta-kernel $\Aut_{G,H}$ comes from the results of \cite{L2}. In that paper we have introduced a scheme $\cL_d(M(F_x))$ of discrete lagrangian lattices in a symplectic Tate space $M(F_x)$ and a certain $\mu_2$-gerbe $\wt\cL_d(M(F_x))$ over it. The complex $\Aut_{G,H}$ on $\Bun_{G,H}$ comes from the stack $\wt\cL_d(M(F_x))$ simply as the inverse image. The key observation is that it is much easier to prove the Hecke property of $\Aut_{G,H}$ on $\wt\cL_d(M(F_x))$, because over the latter stack it is perverse. 
     
\medskip\noindent 
{\bf Acknowledgements.} I am very grateful to Vincent Lafforgue for stimulating discussions which we have regularly for about last two years. They have contributed to this paper. He has also read the first version of the manuscript and indicated several mistakes. I also thank Alexander Braverman for nice discussions. Most of this work was realized in the University Paris 6 and the final part in the Institute for Advanced Study (Princeton), where the author was supported by NSF grant No. DMS-0111298.

 
\bigskip

{\centerline{\scshape 2. Main results}}

\bigskip\noindent
2.1.1 {\scshape Notation} Let $k$ be an algebraically closed field of characteristic $p>2$ (except in Section~3,  where $k=\Fq$). All the schemes (or prestacks) we consider are defined over $k$. 

 Let $X$ be a smooth projective connected curve. Set $F=k(X)$. For a closed point $x\in X$ write $F_x$ for the completion of $F$ at $x$, let $\cO_x\subset F_x$ be the ring of integers. Let $D_x=\Spec\cO_x$ denotes the disc around $x$. Write $\Omega$ for the canonical line bundle on $X$. Let $\Sch^{aff}$ (resp., $\Sch^{aff}_{ft}$) be the category of affine schemes (resp., affine schemes of finite type) over $k$. Let $\Spc$ be the $\infty$-category of spaces as introduced in (\cite{GaRo}, ch. I.1, 1.1.2). Fix a prime $\ell\ne p$. 
 
 We adapt the conventions regarding sheaf theory from (\cite{G5}, Section~1.2) taking as sheaf theory the \'etale $\Qlb$-sheaves on affine schemes of finite type. So, $\DGCat_{cont}$ denotes the $\infty$-category of $\Qlb$-linear cocomplete presentable $\DG$-categories with continuous exact $\Qlb$-linear functors (defined in \cite{GaRo}, ch. I.1, 10.3.3). The functor
\begin{equation}
\label{Shv_sheaf_theory_initial}
Shv: (\Sch^{aff}_{ft})^{op}\to\DGCat_{cont}
\end{equation}
attaches to $S$ the ind-completion of the $\DG$-category of constructible \'etale $\Qlb$-sheaves on $S$. 

 Recall that $\PreStk$ is the category of accessible functors
$$
(\Sch^{aff})^{op}\to\Spc
$$
The full subcategory $\PreStk_{lft}\subset \PreStk$ of prestacks locally of finite type consists of functors preserving filtered colimits. It identifies via restriction with the category of all functors $(\Sch^{aff}_{ft})^{op}\to\Spc$. We define
$$
Shv: (\PreStk_{lft})^{op}\to \DGCat_{cont}
$$
as the right Kan extension of (\ref{Shv_sheaf_theory_initial}). For a map $f: Y\to Y'$ in $\PreStk_{lft}$ it gives a functor $f^!: Shv(Y')\to Shv(Y)$, it admits a left adjoint $f_!: Shv(Y)\to Shv(Y')$. For $f$ ind-schematic of ind-finite type we have the functor $f_*: Shv(Y)\to Shv(Y')$. 

Recall the category $\DGCat^{non-cocmpl}$ defined in (\cite{GaRo}, ch. I.1, 10.3.1), it admits filtered colimits\footnote{The projection $\DGCat^{non-cocmpl}\to 1-\Cat$ preserves filtered colimits, where $1-\Cat$ is the $\infty$-category of $\infty$-categories.}. For $C\in\DGCat_{cont}$ we write $C^c\subset C$ for the full subcategory of compact objects, one has $C^c\in\DGCat^{non-cocmpl}$ naturally. 

 For an algebraic $k$-stack $S$ locally of finite type or an ind-scheme of ind-finite type we also write $\D(S)=Shv(S)$, it is equipped with a perverse t-structure. For an ind-scheme of ind-finite type $S$ we also set $\D(S)^{constr}=\D(S)^c$. 
 
 For an algebraic stack $S$ locally of finite type with an affine diagonal we define $\D(S)^{constr}\subset \D(S)$ as the full subcategory consisting of objects that pull back to an object of $Shv(S')^c$ for any $S'\to S$, where $S'\in\Sch^{aff}_{ft}$. 
 
 Recall that $\DGCat^{non-cocmpl}$ has a canonical involution $\cC\mapsto\cC^{op}$ sending $\cC$ to the opposite category (cf. \cite{GaRo}, ch. I.1, 10.3.2). By (\cite{AGKRRV}, Appendix C), if $S$ is an ind-scheme of ind-finite type or an algebraic stack  locally of finite type with an affine diagonal then the Verdier duality gives an equivalence $\DD: (\D(S)^{constr})^{op}\,\iso\, \D(S)^{constr}$ in $\DGCat^{non-cocmpl}$. 

 Write $\P(S)\subset \D(S)$ for the full subcategory of perverse sheaves. Set $\DP(S)=\oplus_{i\in\ZZ}\,  \P(S)[i]\subset \D(S)$. For $K,K'\in \P(S), i,j\in\ZZ$ we define
$$
\Hom_{\DP(S)}(K[i], K'[j])=\left\{
\begin{array}{ll}
\Hom_{\P(S)}(K,K'), & \mbox{for}\;\; i=j\\
0, & \mbox{for}\;\; i\ne j
\end{array}
\right.
$$

 Write $\P^{ss}(S)\subset\P(S)$ for the full subcategory of semi-simple perverse sheaves. Let $\DP^{ss}(S)\subset\DP(S)$ be the full subcategory of objects of the form $\oplus_{i\in\ZZ} \, K_i[i]$ with $K_i\in\P^{ss}(S)$ for all $i$. 

\medskip\noindent
2.1.2. We freely use the extension of the above sheaf theory to placid ind-schemes as defined in (\cite{G5}, Appendix C). So, for a placid ind-scheme $Y\in\PreStk$ we have $Shv(Y)\in\DGCat_{cont}$ defined in (\cite{G5}, C.2), and for any morphism $f: Y\to Y'$ of placid ind-schemes we have the functor $f_*: Shv(Y)\to Shv(Y')$. 

 Let $Z$ be a placid scheme written as $\lim_{i\in I^{op}} Z_i$, where $I$ is a filtered small $\infty$-category, $Z_i\in\Sch_{ft}$ and for $\alpha: i\to j$ in $I$ the map $f_{\alpha}: Z_j\to Z_i$ is affine, surjective and smooth of some relative dimension $d_{\alpha}$. We define $D(Z)$ as $\colim_{i\in I} D(Z_i)$ with respect to the transition functors $f_{\alpha}^*[\dimrel(f_{\alpha})]$. Let $Y$ be a placid ind-scheme written as $Y\,\iso\, \colim_{j\in J} Y_j$, where $J$ is a filtered category, $Y_j$ is a placid scheme as above and 
for $\alpha: i\to j$ in $J$, $f_{\alpha}: Y_i\to Y_j$ is a placid closed embedding. Then we let $D(Y)=\lim_{j\in J^{op}} D(Y_j)$ with respect to the transition functors $f_{\alpha}^!$. Of course, we have $Shv(Y)\,\iso\, D(Y)$, but our notation $D(Y)$ is supposed to recall that $D(Y)$ is naturally equipped with a perverse t-structure (while for any placid ind-scheme $S$ we need more data to introduce a t-structure on $Shv(S)$, see \cite{G5}, Remark C.2.2). Recall also that $\D(Y)\,\iso\, \colim_{j\in J} \D(Y_j)$ in $\DGCat_{cont}$ with respect to the transition functors $(f_{\alpha})_!$. 

 In addition to $\D(Y)$, define $\D(Y)^{constr}\in \DGCat^{non-cocmpl}$ as follows. First, for a placid scheme $Z$ as above we let $\D(Z)^{constr}=\colim_{i\in I} \D(Z_i)^{constr}$ taken in $\DGCat^{non-cocmpl}$. Recall that for each $i\in I$, $\D(Z_i)$ is compactly generated so that $\Ind(\D(Z_i)^{constr})\,\iso\,\D(Z_i)$ canonically. Moreover, for $\alpha: i\to j$ in $I$ the functor $f_{\alpha}^*:\D(Z_i)\to \D(Z_j)$ preserves compactness, as its right adjoint is continuous, so that we get a diagram $I\to \DGCat^{non-cocmpl}$, $i\mapsto \D(Z_i)^{constr}$. We set $\D(Z)^{constr}=\colim_{i\in I} \D(Z_i)$ taken in $\DGCat^{non-cocmpl}$. Now as in (\cite{GaRo}, ch. I.1, 7.2.7), one gets $\Ind(\D(Z)^{constr})\,\iso\, \D(Z)$ canonically. 
 
 Similarly, for the placid ind-scheme $Y\,\iso\, \colim_{j\in J} Y_j$ as above, and a map $\alpha: i\to j$ in $J$, the functor $(f_{\alpha})_!: \D(Y_i)\to \D(Y_j)$ sends the full subcategory $\D(Y_i)^{constr}$ to $\D(Y_j)^{constr}$. We get a diagram $J\to \DGCat^{non-cocmpl}$, $j\mapsto \D(Y_j)^{constr}$. We let $\D(Y)^{constr}=\colim_{j\in J} \D(Y_j)^{constr}$ taken in $\DGCat^{non-cocmpl}$. As above, $\D(Y)\,\iso\, \Ind(\D(Y)^{constr})$ canonically in $\DGCat_{cont}$. 
 
 Since the transition functors in the definition of $\D(Y)^{constr}$ commute with the Verdier duality, we get an equivalence $\DD: (\D(Y)^{constr})^{op}\,\iso\, \D(Y)^{constr}$ in $\DGCat^{non-cocmpl}$ for a placid ind-scheme $Y$ as above. 

\medskip\noindent
2.1.3 Fix a nontrivial character $\psi: \Fp\to\Qlb^*$ and denote by
$\cL_{\psi}$ the corresponding Artin-Shreier sheaf on $\A^1$. Since we are working over an algebraically closed field, we systematically ignore Tate twists (except in Sections~6.3-6.4). For a morphism of algebraic stacks $f:Y\to Z$ we denote by $\dimrel(f)$ the function of a connected component $C$ of $Y$ given by $\dim C-\dim C'$, where $C'$ is the connected component of $Z$ containing $f(C)$. 

 If $V\to S$ and $V^*\to S$ are dual rank $r$ vector bundles over a stack $S$, we normalize 
the Fourier transform $\Four_{\psi}: Shv(V)\to Shv(V^*)$ by 
$\Four_{\psi}(K)=(p_{V^*})_!(\xi^*\cL_{\psi}\otimes p_V^*K)[r]$,  
where $p_V, p_{V^*}$ are the projections, and $\xi: V\times_S V^*\to \A^1$ is the pairing.
 
 Write $\Bun_k$ for the stack of rank $k$ vector bundles on $X$. For $k=1$ we also write $\Pic X$ for the Picard stack $\Bun_1$ of $X$. We have a line bundle $\cA_k$ on $\Bun_k$ with fibre $\det\RG(X,V)$ at $V\in\Bun_k$. View it as a $\ZZ/2\ZZ$-graded placed in  degree $\chi(V)\!\mod 2$. Our conventions about $\ZZ/2\ZZ$-grading are those of (\cite{L1}, 3.1). 
 
 For a sheaf of groups $G$ on a scheme $S$, $\cF^0_G$ denotes the trivial $G$-torsor on $S$. For a representation $V$ of $G$ and a $G$-torsor $\cF_G$ on $S$ we write $V_{\cF_G}=V\times^{G} \cF_G$ for the induced vector bundle on $S$. 
 
  If $H$ is an algebraic group of finite type and pure dimension, assume given a scheme $Z$ with an action of $H$ and an $H$-torsor $\cF_H$ over a scheme $Y$. Then to $\cS\in \D(Z/H)$, $K\in\D(Y)$ one associates their twisted external product $K\tboxtimes\cS\in \D(\cF_H\times^H Z)$ defined by $K\tboxtimes\cS=p^*K\otimes q^*\cS[\dim H]$, where $Y\getsup{p} \cF_H\times^H Z \toup{q} Z/H$ are the projections. Here $Z/H$ is the stack quotient.
 
\medskip\noindent
2.2.1 {\scshape Hecke operators} For a connected reductive group $G$ over $k$, let $\cH_G$ be the Hecke stack classifying $(x, \cF_G,\cF'_G, \beta)$, where $\cF_G, \cF'_G$ are $G$-torsors on $X$, $x\in X$ and $\beta: \cF_G\mid_{X-x}\,\iso\, \cF'_G\mid_{X-x}$ is an isomorphism. We have a diagram of projections
$$
X\times\Bun_G \getsup{\supp\times h^{\la}_G} \cH_G \toup{h^{\ra}_G}
\Bun_G,
$$
where $h^{\la}_G$ (resp., $h^{\ra}_G$, $\supp$) sends the above collection to $\cF_G$ (resp., $\cF'_G$, $x$). Write $_x\cH_G$ for the fibre of $\cH_G$ over $x\in X$. 

 Let $T\subset B\subset G$ be a maximal torus and Borel subgroup, we write $\Lambda_G$ (resp., $\check{\Lambda}_G$) for the coweights (resp., weights) lattice of $G$. Let $\Lambda^+_G$ (resp., $\check{\Lambda}^+_G$) denote the set of dominant coweights (resp., dominant weights) of $G$. Write $\check{\rho}_G$ (resp., $\rho_G$) for the half sum of the positive roots (resp., coroots) of $G$, $w_0$ for the longest element of the Weyl group of $G$. For $\check{\lambda}\in\check{\Lambda}^+_G$ we write $\cV^{\check{\lambda}}$ for the corresponding Weyl $G$-module.  
 
 For $x\in X$ we write $\Gr_{G,x}$ for the affine grassmanian $G(F_x)/G(\cO_x)$ (cf. \cite{BG}, Section 3.2 for a detailed discussion). It can be seen as an ind-scheme classifying a $G$-torsor $\cF_G$ on $X$ together with a trivialization $\beta:\cF_G\mid_{X-x}\,\iso\, \cF^0_G\mid_{X-x}$ over $X-x$. For $\lambda\in\Lambda^+_G$ let
$\Gr^{\lambda}_{G,x}$ be the $G(\cO_x)$-orbit through $t^{\lambda}\in\Gr_{G,x}$, and $\ov{\Gr}^{\lambda}_{G,x}$ its closure. 

 Let $\cA^{\lambda}_G$ denote the intersection cohomology sheaf of $\ov{\Gr}^{\lambda}_G$. Write $\check{G}$ for the Langlands dual group to $G$ over $\Qlb$. Write $\Sph_G$ for the category of $G(\cO_x)$-equivariant perverse sheaves on $\Gr_{G,x}$. By (\cite{MV}), this is a tensor category, and there is a canonical equivalence of tensor categories $\Loc: \Rep(\check{G})\,\iso\, \Sph_G$, where $\Rep(\check{G})$ is the category of $\check{G}$-representations over $\Qlb$. 
Under this equivalence $\cA^{\lambda}_G$ corresponds to the irreducible representation of $\check{G}$ with h.w. $\lambda$. 

 Write $\Bun_{G,x}$ for the stack classifying $\cF_G\in\Bun_G$ together with a trivialization $\cF_G\,\iso\,\cF^0_G\mid_{D_x}$. 
Following (\cite{BG}, Section~3.2.4), write $\id^l, \id^r$ for the isomorphisms
$$
_x\cH_G\,\iso\, \Bun_{G,x}\times^{G(\cO_x)} \Gr_{G,x}
$$
such that the projection to the first factor corresponds to $h^{\la}_G, h^{\ra}_G$ respectively. Let $_x\ov{\cH}^{\lambda}_G\subset {_x\cH_G}$ be the closed substack that identifies with $\Bun_{G,x}\times^{G(\cO_x)}\ov{\Gr}^{\lambda}_{G,x}$ via $\id^l$. 

 To $\cS\in\Sph_G, K\in\D(\Bun_G)$ one attaches their twisted external products $(K\tboxtimes \cS)^l$ and $(K\tboxtimes \cS)^r$ on $_x\cH_G$, they are normalized to be perverse for $K,\cS$ perverse (cf. \cite{BG}, Section~0.4.4). The Hecke functors
$$
_x\H^{\la}_G, {_x\H^{\ra}_G}: \Sph_G\times\D(\Bun_G)\to\D(\Bun_G)
$$
are given by 
$$
_x\H^{\la}_G(\cS,K)=(h^{\la}_G)_!(\ast\cS\tboxtimes K)^r\;\;\;\mbox{and}\;\;\; {_x\H^{\ra}_G}(\cS,K)=(h^{\ra}_G)_!(\cS\tboxtimes K)^l
$$
We have denoted by $\ast: \Sph_G\,\iso\,\Sph_G$ the covariant equivalence of categories induced by the map $G(F_x)\to G(F_x)$, $g\mapsto g^{-1}$. Write also $\ast: \Rep(\check{G})\,\iso\,\Rep(\check{G})$ for the corresponding functor (in view of $\Loc$), it sends an irreducible $\check{G}$-module with h.w. $\lambda$ to the irreducible $\check{G}$-module with h.w. $-w_0(\lambda)$ (cf. \cite{FGV}, Theorem~5.2.6). 

By (\cite{FG}, Proposition~5.3.9), we have canonically 
$$
_x\H^{\la}_G(\ast\cS, K)\,\iso\, {_x\H^{\ra}_G(\cS, K)}
$$
Besides, the functors $K\mapsto {_x\H^{\la}_G(\cS, K)}$ and $K\mapsto {_x\H^{\ra}_G(\DD(\cS), K)}$ are mutually (both left and right) adjoint. 

 Letting $x$ move around $X$, one similarly defines Hecke functors 
$$
\H^{\la}_G, {\H^{\ra}_G}: \Sph_G\times\D(S\times\Bun_G)\to\D(X\times S\times\Bun_G),
$$
where $S$ is a scheme. The Hecke functors are compatible with the tensor structure on $\Sph_G$ and commute with Verdier duality for locally bounded objects (cf. \select{loc.cit}). Sometimes we write $\Rep(\check{G})$ instead of $\Sph_G$ in the definition of Hecke functors in view of $\Loc$. 
 
\medskip\noindent
2.2.2 We introduce the category 
$$
\D\Sph_G:=\oplus_{r\in\ZZ} \; \Sph_G[r] \, \subset \D(\Gr_G) 
$$
It is equipped with a tensor structure, associativity and commutativity constraints so that the following holds. There is a canonical equivalence of tensor categories $\Loc^{\gr}: \Rep(\check{G}\times\Gm)\,\iso\, \D\Sph_G$ such that $\Gm$ acts on $\Sph_G[r]$ by the character $x\mapsto x^{-r}$. So, the grading by cohomological degrees in $\D\Sph_G$ corresponds to grading by the characters of $\Gm$ in $\Rep(\check{G}\times\Gm)$. In cohomological degree zero the equivalence $\Loc^{\gr}$ specializes to $\Loc$. 

 The action of $\Sph_G$ on $\D(\Bun_G)$ extends to an action of $\D\Sph_G$. Namely, we still denote by 
$$
_x\H^{\la}_G, {_x\H^{\ra}_G}: \D\Sph_G\times\D(\Bun_G)\to\D(\Bun_G)
$$  
the functors given by $_x\H^{\la}_G(\cS[r],K)={_x\H^{\la}_G}(\cS, K)[r]$ and 
${_x\H^{\ra}_G}(\cS[r], K)={_x\H^{\ra}_G}(\cS, K)[r]$ for $\cS\in\Sph_G$ and $K\in\D(\Bun_G)$. 
 
  We still denote by $\ast:\D\Sph_G\to\D\Sph_G$ the functor given by $\ast(\cS[i])=(\ast \cS)[i]$ for $\cS\in\Sph_G$. 
 Write $\Loc_X$ for the tensor category of local systems on $X$. Set 
$$
\D\Loc_X=\oplus_{i\in\ZZ} \; \Loc_X[i]\subset \D(X)
$$ 
We also equip it with a tensor structure so that a choice of $x\in X$ yields an equivalence of tensor categories $\Rep(\pi_1(X,x)\times\Gm)\,\iso\, \D\Loc_X$. The cohomological grading in $\D\Loc_X$ corresponds to grading by the characters of $\Gm$. 
 
 For the standard definition of a Hecke eigen-sheaf we refer the reader to (\cite{G3}, Section~2.7). Since we need to take into account the maximal torus of $\SL_2$ of Arthur, we modify this standard definition as follows. 
 
\begin{Def} 
\label{Def_Hecke_sheaf}
Given a tensor functor $E: \Sph_G\to \D\Loc_X$, a $E$-Hecke eigensheaf is an object $K\in\D(\Bun_G)$ equipped with an isomorphism 
$$
\H^{\la}_G(\cS, K)\,\iso\, E(\cS)\boxtimes K[1]
$$
functorial in $\cS\in \Sph_G$ and satisfying the compatibility conditions (as in \select{loc.cit.}). 
 Note that once $x\in X$ is chosen, a datum of $E$ becomes equivalent to a datum of a homomorphism $\sigma: \pi_1(X,x)\times\Gm\to\check{G}$. In other words, we are given a homomorphism $\Gm\to\check{G}$ of algebraic groups over $\Qlb$, and a continuous homomorphism $\pi_1(X,x)\to Z_{\check{G}}(\Gm)$, where $Z_{\check{G}}(\Gm)$ is the centralizer of $\Gm$ in $\check{G}$. 
 
 Given $\sigma: \pi_1(X,x)\times\Gm\to\check{G}$ as above we also write by abuse of notations $\sigma: \pi_1(X,x)\times\Gm\to\check{G}\times\Gm$ for the homomorphism, whose first component is $\sigma$, and the second component $\pi_1(X,x)\times\Gm\to\Gm$ is the projection. 
\end{Def} 
 
\begin{Ex} The constant perverse sheaf $\Qlb[\dim\Bun_G]$ on $\Bun_G$ is a $\sigma$-Hecke eigensheaf for the homomorphism $\sigma: \pi_1(X,x)\times\Gm\to\check{G}$ given by $2\rho:\Gm\to\check{G}$ and trivial on $\pi_1(X,x)$. 
\end{Ex}
 
\medskip\noindent
2.3 {\scshape Theta-sheaf} Let $G_r$ denote the sheaf of automorphisms of $\cO_X^r\oplus \Omega^r$ preserving the natural symplectic form $\wedge^2(\cO_X^r\oplus \Omega^r)\to \Omega$. The stack $\Bun_{G_r}$ of $G_r$-torsors on $X$ classifies $M\in\Bun_{2r}$ equipped with a symplectic form $\wedge^2 M\to \Omega$. 

Recall the following objects introduced in \cite{L1}. Write $\cA_{G_r}$ for the line bundle on $\Bun_{G_r}$ with fibre $\det\RG(X,M)$ at $M$. We view it as a $\ZZ/2\ZZ$-graded line bundle purely of degree zero. Denote by $\Bunt_{G_r}\to\Bun_{G_r}$ the $\mu_2$-gerbe of square roots of $\cA_{G_r}$. The theta-sheaf $\Aut=\Aut_g\oplus\Aut_s$ is a perverse sheaf on $\Bunt_{G_r}$ (cf. \cite{L1} for details). 

\medskip\noindent
2.4. {\scshape Dual pair $\SO_{2m}, \Sp_{2n}$} 

\medskip\noindent
2.4.1 Let $n,m\ge 1$, $G=G_n$ and $\cA_G=\cA_{G_n}$. Let $H=\SO_{2m}$ be the split orthogonal group of rank $m$ over $k$. The stack $\Bun_H$ of $H$-torsors on $X$ classifies: $V\in\Bun_{2m}$, a nondegenerate symmetric form $\Sym^2 V\to\cO_X$, and a compatible trivialization $\gamma: \det V\,\iso\,\cO_X$. Let $\cA_H$ be the ($\ZZ/2\ZZ$-graded) line bundle on $\Bun_H$ with fibre $\det\RG(X,V)$ at $V$. 

 Write $\Bun_{G,H}=\Bun_G\times\Bun_H$. Let
$$
\tau: \Bun_{G,H}\to \Bun_{G_{2nm}}
$$
be the map sending $(M,V)$ to $M\otimes V$ with the induced symplectic form $\wedge^2(M\otimes V)\to\Omega$. The following is proved in (\cite{L}, Proposition~2).

\begin{Lm} 
\label{Lm_description_line_bundle_on_Bun_GH}
There is a canonical $\ZZ/2\ZZ$-graded isomorphism of line bundles on $\Bun_{G,H}$
\begin{equation}
\label{iso_line_bundle_A}
\tau^*\cA_{G_{2nm}}\,\iso\, \cA_H^{2n}\otimes\cA_G^{2m}\otimes\det\RG(X,\cO)^{-4nm}
\end{equation}
\end{Lm}

 Let $\tilde\tau: \Bun_{G,H}\to\Bunt_{G_{2nm}}$ be the map sending $(\wedge^2 M\to\Omega, \Sym^2 V\to\cO)$ to $(\wedge^2(M\otimes V)\to\Omega, \cB)$, where
$$
\cB=\frac{\det\RG(X,V)^n\otimes\det\RG(X, M)^m}{\det\RG(X,\cO)^{2nm}},
$$
and $\cB^2$ is identified with $\det\RG(X, M\otimes V)$ via (\ref{iso_line_bundle_A}). 

\begin{Def} Set $\Aut_{G,H}=\tilde\tau^*\Aut[\dimrel(\tau)]\in \D(\Bun_{G,H})$. 
As in (\cite{L}, Section~3.2) for the diagram of projections
$$
\Bun_H\getsup{\gq}  \Bun_{G,H}\toup{\gp} \Bun_G
$$
define $F_G: \D(\Bun_H)\to\D(\Bun_G)$ and $F_H: \D(\Bun_G)\to\D(\Bun_H)$
by
$$
F_G(K)=\gp_!(\Aut_{G,H}\otimes \gq^*K)[-\dim\Bun_H]
$$ 
$$
F_H(H)=\gq_!(\Aut_{G,H}\otimes \gp^*K)[-\dim\Bun_G]
$$ 
\end{Def}

The Langlands dual groups are $\check{G}\,\iso\, \SO_{2n+1}$ and $\check{H}\,\iso\, \SO_{2m}$ over $\Qlb$. 
For convenience of the reader, we first formulate our main result in particular cases that yield Langlands functoriality. 

\begin{Th} 
\label{Th_1} 
1) Case $n=m$. There is an inclusion $\check{H}\hook{}\check{G}$ such that there exists an isomorphism 
\begin{equation}
\label{iso_Th1_n_equal_m}
\H^{\la}_G(V, F_G(K))\,\iso\, (\id\boxtimes F_G)(\H^{\la}_H(\Res^{\check{G}}_{\check{H}}(V), K))
\end{equation}
over $X\times\Bun_G$ functorial in $V\in \Rep(\check{G})$, $K\in\D(\Bun_H)$, and compatible with the tensor structures on $\Rep(\check{G})$, $\Rep(\check{H})$. Here we denoted by $\id\boxtimes F_G: \D(X\times\Bun_H)\to\D(X\times\Bun_G)$ the corresponding theta-lifting functor. \\
2) Case $m=n+1$. There is an inclusion $\check{G}\hook{}\check{H}$ such that there exists an isomorphism
\begin{equation}
\label{iso_Th1_m_equal_n+1}
\H^{\ra}_H(V, F_H(K))\,\iso\,(\id\boxtimes F_H)(\H^{\ra}_G(\Res^{\check{H}}_{\check{G}}(V), K))
\end{equation}
over $X\times \Bun_H$ functorial in $V\in\Rep(\check{H})$, $K\in\D(\Bun_G)$ and compatible with the tensor structures on $\Rep(\check{G})$, $\Rep(\check{H})$. Here we denoted by $\id\boxtimes F_H: \D(X\times\Bun_G)\to \D(X\times\Bun_H)$ the corresponding theta-lifting functor. 
\end{Th} 

 We will derive Theorem~\ref{Th_1} from the following Hecke property of $\Aut_{G,H}$. 

\begin{Th} 
1) Case $n=m$. There is an inclusion $\check{H}\to\check{G}$ such that there exists an isomorphism 
\begin{equation}
\label{iso_Th2_n_equal_m}
\H^{\la}_G(V, \Aut_{G,H})\,\iso\, \H^{\ra}_H(\Res^{\check{G}}_{\check{H}}(V), \Aut_{G,H})
\end{equation}
in $\D(X\times\Bun_{G,H})$ functorial in $V\in \Rep(\check{G})$ and compatible with the tensor structures on $\Rep(\check{G})$, $\Rep(\check{H})$.\\
2) Case $m=n+1$. There is an inclusion $\check{G}\hook{}\check{H}$ such that there exists an isomorphism
\begin{equation}
\label{iso_Th2_m_equal_n+1}
\H^{\ra}_H(V, \Aut_{G,H})\,\iso\, \H^{\la}_G(\Res^{\check{H}}_{\check{G}}(V), \Aut_{G,H})
\end{equation}
in $\D(X\times\Bun_{G,H})$ functorial in $V\in \Rep(\check{H})$ and compatible with the tensor structures on $\Rep(\check{G})$, $\Rep(\check{H})$.
\end{Th}

\medskip\noindent
2.4.2 In the case $m\le n$ define the map $\kappa: \check{H}\times\Gm\to\check{G}$ as follows. 

 Set $W_0=\Qlb^n$, write $W_0=W_1\oplus W_2$, where $W_1$ (resp., $W_2$) is the subspace generated by the first $m$ (resp., last $n-m$) base vectors. Equip $W_0\oplus W_0^*\oplus\Qlb$ with the symmetric form given by the matrix 
$$
\left(
\begin{array}{ccc}
0 & E_n & 0\\
E_n & 0 & 0\\
0 & 0 & 1
\end{array}
\right),
$$
where $E_n\in\GL_n(\Qlb)$ is the unity.  Realize $\check{G}$ as $\SO(W_0\oplus W_0^*\oplus \Qlb)$. Equip the subspace $W_1\oplus W_1^*\subset W_0\oplus W_0^*\oplus \Qlb$ with the induced symmetric form, and realize $\check{H}$ as $\SO(W_1\oplus W_1^*)$. This fixes the inclusion $i_{\kappa}: \check{H}\hook{}\check{G}$. The centralizer of $\check{H}$ in $\check{G}$ contains the group $\OO(W_2\oplus W_2^*\oplus\Qlb)$. Let $\check{T}_{\GL(W_2)}$ be the maximal torus of diagonal matrices in $\GL(W_2)$. We have $\Hom(\Gm, \check{T}_{\GL(W_2)})=\ZZ^{n-m}$ canonically, and we let $\alpha_{\kappa}=(2,4,\ldots, 2n-2m)\in \Hom(\Gm, \check{T}_{\GL(W_2)})$. View $\alpha_{\kappa}$ as a map $\Gm\to \check{G}$. Finally, set $\kappa=(i_{\kappa}, \alpha_{\kappa}): \check{H}\times\Gm\to \check{G}$. 

 Another way to think of $\alpha_{\kappa}$ is to say that $W_2\oplus W_2^*\oplus\Qlb$ admits an irreducible representation of the $\SL_2$ of Arthur, and $\alpha_{\kappa}$ is its restriction to the standard maximal torus
$$
\Gm\hook{}\SL_2\toup{\sigma} \SO(W_2\oplus W_2^*\oplus\Qlb)
$$
As predicted by Adams (\cite{Ad}), the representation $\sigma$ corresponds to the principal unipotent orbit in $\SO(W_2\oplus W_2^*\oplus\Qlb)$, so $\alpha_{\kappa}=2\rho_{\SO(W_2\oplus W_2^*\oplus\Qlb)}$ for a suitable choice of positive roots of $\SO(W_2\oplus W_2^*\oplus\Qlb)$. 

 Write $\gRes^{\kappa}: \Sph_G\to \D\Sph_H$ for the geometric restriction functor corresponding to $\kappa$. By this we mean the restriction functor $\Rep(\check{G})\to\Rep(\check{H}\times\Gm)$ composed with the Satake equivalences.

 In the case $m>n$ define $\kappa: \check{G}\times\Gm\to\check{H}$ as follows. Set in this case $W_0=\Qlb^m$, let $W_1$ (resp., $W_2$) be the subspace of $W_0$ generated by the first $n$ (resp., last $m-n$) base vectors. Equip $W_0\oplus W_0^*$ with the symmetric form given by the matrix
$$ 
\left(
\begin{array}{cc}
0 & E_m\\
E_m & 0\\
\end{array}
\right),
$$
where $E_m\in\GL_m(\Qlb)$ is the unity. Realize $\check{H}$ as $\SO(W_0\oplus W_0^*)$. 

 Write $\{e_i\}$ for the standard base of $W_0$, and $\{e_i^*\}$ for the dual base in $W_0^*$. Write $W_2=W_3\oplus W_4$, where $W_3$ (resp., $W_4$) is spanned by $e_{n+1}$ (resp., by $e_{n+2},\ldots, e_m$). Let $\bar W\subset W_3\oplus W_3^*$ be any nondegenerate one-dimensional subspace. Equip $W_1\oplus W_1^*\oplus \bar W$ with the induced symmetric form and set $\check{G}=\SO(W_1\oplus W_1^*\oplus \bar W)$. This fixes the inclusion $i_{\kappa}:\check{G}\hook{}\check{H}$. 

  Let $\bar W^{\perp}$ denote the orthogonal complement of $\bar W$ in $W_2\oplus W_2^*$. The centralizer of $\check{G}$ in $\check{H}$ contains $\OO(\bar W^{\perp})$. Realize $\GL(W_4)$ as the Levi subgroup of $\SO(\bar W^{\perp})$ using the standard inclusion $W_4\oplus W_4^*\subset \bar W^{\perp}$. Let $\check{T}_{\GL(W_4)}$ be the maximal torus of diagonal matrices in $\GL(W_4)$. Set 
$$
\alpha_{\kappa}=(-2, -4, \ldots, 2-2m+2n)\in \ZZ^{m-n-1}=\Hom(\Gm, \check{T}_{\GL(W_4)})
$$  
View $\alpha_{\kappa}$ as a map $\Gm\to\check{H}$, set $\kappa=(i_{\kappa}, \alpha_{\kappa}): \check{G}\times\Gm\to\check{H}$.  

 Another way to think of $\alpha_{\kappa}$ is to say that $\bar W^{\perp}$ can be thought of as an irreducible representation of \select{the $\SL_2$ of Arthur}, and $\alpha_{\kappa}$ is the restriction to the standard maximal torus
$$
\Gm\hook{}\SL_2\toup{\sigma} \SO(\bar W^{\perp})
$$   
As predicted by Adams (\cite{Ad}), the representation $\sigma$ corresponds to the principal unipotent orbit in $\SO(\bar W^{\perp})$, so $\alpha_{\kappa}=2\rho_{\SO(\bar W^{\perp})}$ for a suitable choice of positive roots of $\SO(\bar W^{\perp})$. 
As above, the geometric restriction functor corresponding to $\kappa$ is denoted $\gRes^{\kappa}: \Sph_H\to\D\Sph_G$. 

 Here is our main global result. 
 
\begin{Th} 
\label{Th_main_global_symplectic_orthogonal}
1) Case $m\le n$. There exists an isomorphism
\begin{equation}
\label{iso_Th_main_global_m_less_n}
\H^{\la}_G(\cS, F_G(K))\,\iso\, (\id\boxtimes F_G)(\H^{\la}_H(\gRes^{\kappa}(\cS), K))
\end{equation}
in $\D(X\times\Bun_G)$ functorial in $\cS\in \Sph_G$, $K\in\D(\Bun_H)$ and compatible with the tensor structures on $\Sph_G$, $\Sph_H$. Here we denoted by $\id\boxtimes F_G: \D(X\times\Bun_H)\to\D(X\times\Bun_G)$ the corresponding theta-lifting functor. 

\smallskip\noindent
2) Case $m>n$. There exists an isomorphism
\begin{equation}
\label{iso_Th_main_global_m_grater_n}
\H^{\la}_H(\cS, F_H(K))\,\iso\,(\id\boxtimes F_H)(\H^{\ra}_G(\gRes^{\kappa}(\ast\cS), K))
\end{equation}
in $\D(X\times \Bun_H)$ functorial in $\cS\in\Sph_H$, $K\in\D(\Bun_G)$ and compatible with the tensor structures on $\Sph_H$, $\Sph_G$. Here we denoted by $\id\boxtimes F_H: \D(X\times\Bun_G)\to \D(X\times\Bun_H)$ is the corresponding theta-lifting functor. 
\end{Th} 
 
 We will derive Theorem~\ref{Th_main_global_symplectic_orthogonal} from the following Hecke property of $\Aut_{G,H}$. 
 
\begin{Th}
\label{Th_Hecke_property_Aut_GH}
1) Case $m\le n$. There exists an isomorphism
\begin{equation}
\label{iso_Th_Hecke_Aut_GH_m_less_n}
\H^{\la}_G(\cS, \Aut_{G,H})\,\iso\, \H^{\la}_H(\ast\gRes^{\kappa}(\cS), \Aut_{G,H})
\end{equation}
in $\D(X\times\Bun_{G,H})$ functorial in $\cS\in \Sph_G$ and compatible with the tensor structures on $\Sph_G$,$\Sph_H$.\\
2) Case $m>n$. There exists an isomorphism
\begin{equation}
\label{iso_Th_Hecke_Aut_GH_m_grater_n}
\H^{\la}_H(\cS, \Aut_{G,H})\,\iso\, \H^{\la}_G(\gRes^{\kappa}(\ast\cS), \Aut_{G,H})
\end{equation}
in $\D(X\times\Bun_{G,H})$ functorial in $\cS\in \Sph_H$ and compatible with the tensor structures on $\Sph_H$, $\Sph_G$.
\end{Th}

\begin{Rem}
The compatibility of the isomorphism (\ref{iso_Th_Hecke_Aut_GH_m_less_n}) with the tensor structures means the following. If we denote (\ref{iso_Th_Hecke_Aut_GH_m_less_n})  by $\alpha_{\cS}$ then for $\cS,\cS'\in\Sph_G$ the diagram commutes
$$
\begin{array}{ccc}
\H^{\la}_G(\cS, \H^{\la}_G(\cS', \Aut_{G,H})) & \toup{\alpha_{\cS'}} & \H^{\la}_G(\cS, \H^{\la}_H(\ast\gRes^{\kappa}(\cS'), \Aut_{G,H}))\\
\downarrow\lefteqn{\scriptstyle \alpha_{\cS\star\cS'}} && \downarrow\lefteqn{\scriptstyle \alpha_{\cS}}\\
\H^{\la}_H(\ast\gRes^{\kappa}(\cS\star \cS'), \Aut_{G,H})) & \iso & \H^{\la}_H(\ast\gRes^{\kappa}(\cS), \H^{\la}_H(\ast\gRes^{\kappa}(\cS'), \Aut_{G,H}))
\end{array}
$$
Here $\cS\star\cS'$ denotes the convolution in $\Sph_G$. We used here the canonical commutation of the Hecke functors for $G$ and for $H$. For (\ref{iso_Th_Hecke_Aut_GH_m_grater_n}) the compatibility with the tensor structure is similar. 
\end{Rem}

\smallskip\noindent
2.4.3 There is an automorphism $\sigma_H: \check{H}\,\iso\,\check{H}$ inducing the functor $\ast: \Rep(\check{H})\,\iso\,\Rep(\check{H})$ defined in Section~2.2.2. For $m>n$ write $\tilde\kappa=\sigma_H\comp\kappa$. Note also that the functor $\ast: \Rep(\check{G})\,\iso\,\Rep(\check{G})$ is isomorphic to the identity functor. From Theorem~\ref{Th_main_global_symplectic_orthogonal} one derives the following. 

\begin{Cor} 
1) For $m\le n$ let $K\in \D(\Bun_H)$ be a $\sigma$-Hecke eigensheaf for some $\sigma: \pi_1(X,x)\times\Gm\to \check{H}$. Let $\tau$ be the composition
$$
\pi_1(X,x)\times\Gm \,\toup{\sigma} \,\check{H}\times\Gm\,\toup{\kappa}\,
\check{G},
$$
where we used our convention from Definition~\ref{Def_Hecke_sheaf}. Then $F_G(K)$ is equipped with a structure of a $\tau$-Hecke eigensheaf.  

\smallskip
\noindent
2) For $m>n$ let $K\in\D(\Bun_G)$ be a $\sigma$-Hecke eigensheaf for some $\sigma:\pi_1(X,x)\times\Gm\to\check{G}$. Let $\tau$ be the composition
$$
\pi_1(X,x)\times\Gm\,\toup{\sigma}\,\check{G}\times
\Gm\,\toup{\tilde\kappa}\, \check{H},
$$
Then $F_H(K)$ is equipped with a structure of a $\tau$-Hecke engeinsheaf.
\end{Cor}

\medskip\noindent
2.5 {\scshape Dual pair $\GL_m, \GL_n$}

\medskip\noindent
Let $n,m\ge 0$. Recall that $\Bun_n$ denotes the stack of rank $n$ vector bundles on $X$. Our convention is that $\GL_0=\{1\}$ and $\Bun_0=\Spec k$. 

 Let $\cW_{n,m}$ denote the stack classifying $L\in\Bun_n, U\in\Bun_m$ and a section $s:\cO_X\to L\otimes U$. We have a diagram
$$
\Bun_n\,\getsup{h_n}\,\cW_{n,m}\,\toup{h_m}\Bun_m,
$$
where $h_m$ (resp., $h_n$) sends $(L,U,s)$ to $U$ (resp., to $L$). Let $\cW'_{n,m}$ be the stack classifying $L\in\Bun_n, U\in\Bun_m$ and a section $s': L\otimes U\to\Omega$. We have a diagram
$$
\Bun_n\,\getsup{h'_n}\,\cW'_{n,m}\,\toup{h'_m}\,\Bun_m,
$$
where $h'_m$ (resp., $h'_n$) sends $(L,U, s')$ to $U$ (resp., to $L$). 
\begin{Def} The theta-lifting functors $F_{n,m}, F'_{n,m}:\D(\Bun_n)\to\D(\Bun_m)$ are given by
$$
F_{n,m}(K)= (h_m)_!h_n^*K[\dim\Bun_m+a_{n,m}]\;\;\;\;\;\mbox{and}\;\;\;\;\;
F'_{n,m}(K)=(h'_m)_!(h'_n)^*K[\dim\Bun_m-a_{n,m}]
$$
Here $a_{n,m}$ is a function of a connected component of $\Bun_n\times\Bun_m$ given by $a_{n,m}=\chi(L\otimes U)$ for $L\in\Bun_n, U\in\Bun_m$. By restriction under $h_n\times h_m$ (resp., under $h'_n\times h'_m$), we view $a_{n,m}$ in the above formulas as a function on $\cW_{n,m}$ (resp., on $\cW'_{n,m}$). 
\end{Def}

 Since $h_m$ and $h'_m$ are not representable, a priori  $F_{n,m}$ and $F'_{n,m}$ may send a bounded complex to an unbounded one. The following result can be thought of as a functional equation for the theta-lifting functors. 
 
\begin{Lm} There is a canonical isomorphism of functors $F'_{n,m}\,\iso\, F_{n,m}$. 
\end{Lm} 
\begin{Prf}
Write $\phi, \phi'$ for the projections from $\cW_{n,m}$ and from $\cW'_{n,m}$ to $\Bun_n\times\Bun_m$. As in (\cite{BG}, Lemma~7.3.6) one shows that $\phi_!\Qlb[a_{n,m}]\,\iso\,\phi'_!\Qlb[-a_{n,m}]$ canonically. The assertion follows. 
\end{Prf}

\medskip

 For the rest of Section~2.5 assume $m\ge n$ and set $G=\GL(L_0)$ and $H=\GL(U_0)$ for $U_0=k^m$, $L_0=k^n$. 
Write $U_0=U_1\oplus U_2$, where $U_1$ (resp., $U_2$) is the subspace generated by the first $n$ (resp., last $m-n$) base vectors. 
Let $M=\GL(U_1)\times\GL(U_2)\subset H$ be the corresponding Levi factor. 

 Define $\kappa:\check{G}\times\Gm\to\check{H}$ as the composition
$$
\check{G}\times\Gm\toup{\id\times 2\check{\rho}_{\GL(U_2)}}
\check{G}\times \check{\GL}(U_2)=\check{M}\hook{}\check{H}
$$
Write $\gRes^{\kappa}:\Sph_H\to\D\Sph_G$ for the corresponding geometric restriction functor. 

 The analog of Theorem~\ref{Th_main_global_symplectic_orthogonal} for the dual pair $(G,H)$ is as follows. 
 
\begin{Th}
\label{Th_main_global_GL_m_GL_n} 
We assume $m\ge n$.  There exists an isomorphism
\begin{equation}
\label{iso_Th_main_global_GL_m_GL_n}
\H^{\la}_H(\cS, F_{n,m}(K))\,\iso\, (\id\boxtimes F_{n,m})(\H^{\ra}_G(\gRes^{\kappa}(\cS), K))
\end{equation}
in $\D(X\times\Bun_m)$ functorial in $\cS\in \Sph_H$, $K\in\D(\Bun_n)$ and compatible with the tensor structures on $\Sph_H$, $\Sph_G$. Here we denoted by $\id\boxtimes F_{n,m}: \D(X\times\Bun_n)\to\D(X\times\Bun_m)$ the corresponding theta-lifting functor. 
\end{Th} 

 If $n=m$ or $m=n+1$ then the restriction of $\kappa$ to $\Gm$ is trivial, so Theorem~\ref{Th_main_global_GL_m_GL_n} in this case says that $F_{n,m}$ realizes the (non ramified) geometric Langlands functoriality with respect to an inclusion $\check{G}\hook{}\check{H}$. For example, for $n=m$ one may show the following. For an irreducible rank $n$ local system $E$ on $X$ write $\Aut_E$ for the automorphic sheaf on $\Bun_n$ corresponding to $E$ (cf. \cite{FGV}). Then $F_{n,n}(\Aut_E)$ is isomorphic to $\Aut_{E^*}$ tensored by some constant complex. 
 
 Write $_{\infty}\cW_{n,m}$ for the stack classifying $x\in X$, $L\in\Bun_n, U\in\Bun_m$ and a section $s:\cO_X\to L\otimes U(\infty x)$, which is allowed to have an arbitrary pole at $x$. This is an ind-algebraic stack. For a closed point $x\in X$ let $_{x,\infty}\cW_{n,m}\subset {_{\infty}\cW_{n,m}}$ be the closed stack given by fixing $x$.
  
  In Section~7 we will define Hecke functors 
\begin{equation}
\label{Hecke_functor_H_forW_mn}
_x\H^{\la}_H, {_x\H^{\ra}_H}:\; \Sph_H\times\D(_{x,\infty}\cW_{n,m})\to \D(_{x,\infty}\cW_{n,m})
\end{equation}
\begin{equation}
\label{Hecke_functor_G_forW_mn}
_x\H^{\la}_G, {_x\H^{\ra}_G}: \; \Sph_G\times\D(_{x,\infty}\cW_{n,m})\to \D(_{x,\infty}\cW_{n,m})
\end{equation}
and their family versions acting on $\D(_{\infty}\cW_{n,m})$.
 Set 
\begin{equation}
\label{object_cI}
\cI=(\Qlb)_{\cW_{n,m}}[\dim\Bun_m+\dim\Bun_n+a_{n,m}],
\end{equation}
where $a_{n,m}$ is a function of a connected component of $\cW_{n,m}$ defined above. View $\cI$ as a complex on $\cW_{n,m}$ extended by zero to $_{x,\infty}\cW_{n,m}$. Write also $_{\infty}\cI$ for $\cI\boxtimes\Qlb[1]$ over $\cW_{n,m}\times X$ extended by zero to $_{\infty}\cW_{n,m}$. We will derive Theorem~\ref{Th_main_global_GL_m_GL_n} from the following `Hecke property' of $\cI$. 

\begin{Th}
\label{Th_global_Hecke_property_I_GL_n_GL_m}
The two functors $\Sph_H\to \D(_{\infty}\cW_{m,n})$ given by
$$
\cT\mapsto {_x\H^{\la}_H}(\cT, {_{\infty}\cI})\;\;\;\;\mbox{and}\;\;\;\;
\cT\mapsto {_x\H^{\la}_G}(\gRes^{\kappa}(\cT), {_{\infty}\cI})
$$
are isomorphic in a way compatible with the tensor structures on $\Sph_H, \Sph_G$.
\end{Th}

\bigskip\noindent

{\centerline{\scshape 3. Classical setting and motivations}}

\bigskip\noindent
In Section~3 we assume $k=\Fq$. 

\medskip\noindent
3.1 {\scshape Weil representation of $\GL_m\times\GL_n$.}
 
\medskip\noindent 
Let $U_0$ (resp., $L_0$) be a $k$-vector space of dimension $m$ (resp., $n$). For Section~3.1 set $G=\GL(L_0)$ and  $H=\GL(U_0)$. Let $\Pi_0=U_0\otimes L_0$ and $\Pi=\Pi_0(\cO)$. 

 Let $x\in X$ be a closed point. Recall that the Weil representation of $G(F_x)\times H(F_x)$ can be realized in the Schwarz space $\cS(\Pi(F_x))$ of locally constant compactly supported $\Qlb$-valued functions on $\Pi(F_x)$. The action of $G(F_x)\times H(F_x)$ on this space comes from its natural action on $\Pi(F_x)$. 

 Write $\cH_x(G)$ for the Hecke algebra of the pair $(G(F_x), G(\cO_x))$, and similarly for $\cH_x(H)$. Recall that $\cH_x(G)$ identifies canonically with the Grothendieck group $K(\Rep(\check{G}))$ of the category $\Rep(\check{G})$ of $\check{G}$-representations over $\Qlb$. 
 
 The space of invariants $\cS(\Pi(F_x))^{G(\cO_x)\times H(\cO_x)}$ is naturally a module over $\cH_x(G)\otimes\cH_x(H)$. Let $\phi_0\in\cS(\Pi(F))$ be the characteristic function of $\Pi(\cO)$. The following result is well-known (cf. \cite{MVW}, \cite{R}), in Section~5 we prove its geometric version.   
 
\begin{Lm} Assume $m\ge n$. The map $\cH_x(G)\to \cS(\Pi(F))^{(G\times H)(\cO_x)}$ sending $h$ to $h\phi_0$ is an isomorphism of $\cH_x(G)$-modules. There is a homomorphism $\kappa: \cH_x(H)\to\cH_x(G)$ such that the $\cH_x(H)$-action on 
$\cS(\Pi(F_x))^{G(\cO_x)\times H(\cO_x)}$ factors through $\kappa$.
\end{Lm} 

 For $n=m$ the homomorphism $\kappa$ comes from the functor $\Rep(\check{H})\to\Rep(\check{G})$ of restriction with respect to an isomorphism $\check{G}\,\iso\, \check{H}$. For $m\ge n$ we will see that $\kappa$ comes from the functor $\Rep(\check{H})\to \Rep(\check{G}\times\Gm)\,\iso\,\D\Sph_G$ of restriction with respect to a homomorphism $\check{G}\times\Gm\to\check{H}$. For $m>n+1$ the restriction of this homomorphism to $\Gm$ is nontrivial.

\bigskip\noindent
3.2 {\scshape Weil representation of $\SO_{2m}\times\Sp_{2n}$}

\medskip\noindent
3.2.1 In this sebsection we introduce some objects on the level of functions whose geometric analogs are used in the proof of Theorem~\ref{Th_main_global_symplectic_orthogonal}. Keep the notation of Section~2.4. Let $U_0=\cO^m_X$ and $V_0=U_0\oplus U_0^*$, we equip $V_0$ with the symmetric form $\Sym^2 V_0\to\cO_X$ given by the pairing between $U_0$ and $U_0^*$, so $U_0$ and $U_0^*$ are isotropic subbundles in $V_0$. Think of $V_0$ as the standard representation of $H$. 

 Let $P(H)\subset H$ be the parabolic subgroup preserving $U_0$, let $U(H)\subset P(H)$ be its unipotent radical, so $U(H)\,\iso\, \wedge^2 U_0$ canonically.

 Let $L_0=\cO_X^n$ and $M_0=L_0\oplus L_0^*\otimes\Omega$. We equip $M_0$ the symplectic form $\wedge^2 M_0\to\Omega$ given by the pairing $L_0$ with $L_0^*\otimes\Omega$. So, $L_0$ and $L_0^*\otimes\Omega$ are lagrangian subbundles in $M_0$. Recall that $G$ is the group scheme over $X$ of automorphisms of $M_0$ preserving the symplectic form. 

 Let $P(G)\subset G$ be the parabolic subgroup preserving $L_0$, write $U(G)\subset P(G)$ for its unipotent radical. We have $U(G)\,\iso\, \Omega^{-1}\otimes\Sym^2 L_0$ canonically.  
 
  Set $\cM_0=V_0\otimes M_0$, it is equipped with a symplectic form, which is the tensor product of forms on $V_0$ and $M_0$. 
 
 Set $F=k(X)$. Let $\AA$ be the ad\`eles ring of $F$, $\cO\subset\AA$ be the entire adeles. Let $\chi:\Omega(\AA)/\Omega(F)\to \Qlb^*$ denote the character
$$
\chi(\omega)=\psi(\sum_{x\in X} \tr_{k(x)/k} \Res \omega_x)
$$ 
Let $\Heis=\cM_0\oplus \Omega$ be the Heisenberg group over $X$ constructed out of the symplectic bundle $\cM_0$. The product in $\Heis$ is given by
$$
(m_1,\omega_1)(m_2,\omega_2)=(m_1+m_2, \omega_1+\omega_2+\frac{1}{2}\<m_1,m_2\>)
$$
For the generalities on the metaplectic extension $\wt\Sp(\cM_0)$ of $\Sp(\cM_0)$ and its Weil representation we refer the reader to \cite{L1}. 
The natural map $G(\AA)\times H(\AA)\to \Sp(\cM_0)(\AA)$ lifts naturally to a homomorphism $G(\AA)\times H(\AA)\to\wt\Sp(\cM_0)(\AA)$. We use two Schr\"odinger models of the corresponding Weil representation of $G(\AA)\times H(\AA)$. 

 Set $\cL_0=V_0\otimes L_0\subset V_0\otimes M_0$, this is a Lagrangian subbundle in $\cM_0$. Let 
$$
\chi_{\cL}: \cL_0(\AA)\oplus \Omega(\AA)\to\Qlb^*
$$
denote the character $\chi_{\cL}(u, \omega)=\chi(\omega)$. Let $\cS_{\cL,\psi}$ denote the induced representation of $(\cL_0(\AA)\oplus \Omega(\AA), \chi_{\cL})$ to $\Heis(\AA)$. By definition, $\cS_{\cL,\psi}$ is the space of functions $f: \Heis(\AA)\to\Qlb$ satisfying:
\begin{itemize}
\item $f(ah)=\chi_{\cL}(a)f(h)$ for $a\in \cL_0(\AA)\oplus\Omega(\AA)$, $h\in\Heis(\AA)$;
\item there is an open subgroup $U\subset \cM_0(\AA)$ such that $f(h(u,0))=f(h)$ for $u\in U$, $h\in\Heis(\AA)$;
\item $f$ is of compact support modulo $\cL_0(\AA)\oplus\Omega(\AA)$.
\end{itemize}

 For a free $\AA$-module (or free $F_x$-module) $R$ of finite type denote by $\cS(R)$ the Schwarz space of locally constant compactly supported $\Qlb$-valued functions on $R$. We have an isomorphism $\cS_{\cL,\psi}\,\iso\,\cS(V_0\otimes L_0^*\otimes\Omega(\AA))$ sending $f$ to $\phi$ given by $\phi(v)=f(v,0)$, $v\in V_0\otimes L_0^*\otimes\Omega(\AA)$.  

 The theta-functional 
$$
\Theta_{\cL}: \cS(V_0\otimes L_0^*\otimes\Omega(\AA))\to\Qlb
$$ 
is given by 
$$
\Theta_{\cL}(\phi)=\sum_{v\in V_0\otimes L_0^*\otimes\Omega(F)} \phi(v)
\;\;\;\;\; \mbox{for}\;\; \phi\in \cS(V_0\otimes L_0^*\otimes\Omega(\AA))
$$

 Set $\cU_0=U_0\otimes M_0$, this is a Lagrangian subbundle in $\cM_0$. Let 
$$
\chi_{\cU}: \cU_0(\AA)\oplus\Omega(\AA)\to\Qlb^*
$$
be the character $\chi_{\cU}(u,\omega)=\chi(\omega)$. Let $\cS_{\cU,\psi}$ denote the induced representation of $(\cU_0(\AA)\oplus \Omega(\AA), \chi_{\cU})$ to $\Heis(\AA)$. As above, we identify it with the Schwarz space $\cS(U_0^*\otimes M_0(\AA))$. 
 
 The theta-functional $\Theta_{\cU}: \cS(U_0^*\otimes M_0(\AA))\to\Qlb$ is given by  
$$
\Theta_{\cU}(\phi)=\sum_{t\in U_0^*\otimes M_0(F)} \phi(t)
\;\;\;\;\; \mbox{for}\;\;\phi\in \cS(U_0^*\otimes M_0(\AA))
$$

 For a locally free $\cO_X$-module $\cY$ of finite type write $\chi(\cY)$ for the Euler characteristic of $\cY$. Set $\epsilon=q^{\chi(U_0^*\otimes L_0)}$. Let us construct a diagram of $H(\AA)\times G(\AA)$-representations
\begin{equation}
\label{diag_one}
\begin{array}{ccc}
\cS(U_0^*\otimes M_0(\AA)) & \toup{\theta_{\cU}} & \Funct((H\times G)(F)\backslash (H\times G)(\AA))\\
\uparrow\lefteqn{\scriptstyle \zeta} & \nearrow\lefteqn{\scriptstyle \epsilon\theta_{\cL}}\\
\cS(V_0\otimes L_0^*\otimes\Omega(\AA)),
\end{array}
\end{equation}
where $H(\AA)\times G(\AA)$ acts on the space of functions $\Funct((H\times G)(F)\backslash (H\times G)(\AA))$ by right translations. 
The map $\theta_{\cU}$ sends $\phi$ to $\theta_{\cU,\phi}$ given by $\theta_{\cU,\phi}(h,g)=\Theta_{\cU}((h,g)\phi)$. The map $\theta_{\cL}$ sends $\phi$ to $\theta_{\cL,\phi}$ given by 
$$
\theta_{\cL,\phi}(h,g)=\Theta_{\cL}((h,g)\phi)
$$ 

 For $\phi\in \cS(V_0\otimes L_0^*\otimes\Omega(\AA))$ let $\zeta \phi\in \cS(U_0^*\otimes M_0(\AA))$ be given by
\begin{equation}
\label{map_zeta}
(\zeta \phi)(b)=\int_{U_0\otimes L_0^*\otimes\Omega(\AA)} \chi(\<a, b_1\>)\phi(a+b_2)da
\end{equation}
Here for $b\in U_0^*\otimes M_0(\AA)$ we write $b=b_1+b_2$ with $b_1\in U_0^*\otimes L_0(\AA)$ and $b_2\in U_0^*\otimes L_0^*\otimes\Omega(\AA)$, and 
$da$ is the Haar measure on $U_0\otimes L_0^*\otimes\Omega(\AA)$ normalized by requiring that the volume of $U_0\otimes L_0^*\otimes\Omega(\cO)$ is one. It is known that $\zeta$ is an isomorphism of $G(\AA)\times H(\AA)$-modules (cf. \cite{MVW}). 

 Let $\phi_{0,\cU}$ (resp., $\phi_{0,\cL}$) be the characteristic function of $U_0^*\otimes M_0(\cO)$ (resp., of $V_0\otimes L_0^*\otimes\Omega(\cO)$). An easy calculation shows that $\zeta \phi_{0,\cL}=\phi_{0,\cU}$. 

\begin{Lm} The diagram (\ref{diag_one}) commutes. 
\end{Lm}
\begin{Prf} 
We have $\Theta_{\cL}(\phi_{0,\cL})=q^{\dim \H^0(X, V_0\otimes L_0^*\otimes\Omega)}$ and $\Theta_{\cU}(\phi_{0,\cU})=q^{\dim \H^0(X, U_0^*\otimes M_0)}$. Since 
$$
\dim \H^0(X, U_0^*\otimes M_0)=\chi(U_0^*\otimes L_0)+
\dim \H^0(X, V_0\otimes L_0^*\otimes\Omega),
$$ 
we get $\Theta_{\cU}\comp \zeta=\epsilon\Theta_{\cL}$. Since $\zeta$ is an isomorphism of $G(\AA)\times H(\AA)$-modules,
our assertion follows. 
\end{Prf}

 \medskip
 
 Write $\cH(H)$ for the Hecke algebra of the pair $(H(\cO), H(\AA))$, and similarly for $G$. Passing to the $(G\times H)(\cO)$-invariants, one gets from (\ref{diag_one}) the commutative diagram 
\begin{equation}
\label{diag_invariants}
\begin{array}{ccc}
\cS(U_0^*\otimes M_0(\AA))^{(H\times G)(\cO)} & \toup{\theta_{\cU}} & \Funct(\Bun_{G,H}(k))\\
\uparrow\lefteqn{\scriptstyle \zeta} & \nearrow\lefteqn{\scriptstyle \epsilon\theta_{\cL}}\\
\cS(V_0\otimes L_0^*\otimes\Omega(\AA))^{(H\times G)(\cO)}
\end{array}
\end{equation}
of $\cH(H)\otimes\cH(G)$-modules. The notation $\Bun_{G,H}$ is that of Section~2.4.1. 

 Let $\phi_0\in \Funct(\Bun_{G,H}(k))$ be the function trace of Frobenius of $\Aut_{G,H}$. Then $\theta_{\cU}\phi_{0,\cU}$ equals $\phi_0$ up to a multiple. 
 
 For $x\in X$ let $\phi_{0,\cU,x}\in\cS(U_0^*\otimes M_0(F_x))$ be the characteristic function of $U_0^*\otimes M_0(\cO_x)$, let $\phi_{0,\cL,x}\in \cS(V_0\otimes L_0^*\otimes\Omega(F_x))$ be the characteristic function of $V_0\otimes L_0^*\otimes\Omega(\cO_x)$. 

 Denote by $\cH_x(G)$ the Hecke algebra of the pair $(G(\cO_x), G(F_x))$, and similarly for $H$. Recall the decomposition as a restricted tensor product 
$$
\cH(G)\,\iso\,
\mathop{\otimes'}\limits_{x\in X} \cH_x(G)
$$ 
Similarly, we have
$$
\cS(U_0^*\otimes M_0(\AA))\,\iso\, \mathop{\otimes'}\limits_{x\in X} \cS(U_0^*\otimes M_0(F_x))
$$
In view of this isomorphism $\cS(U_0^*\otimes M_0(\AA))$ is generated as a $\Qlb$-vector space by functions of the form $\otimes_x \phi_x$ with $\phi_x\in \cS(U_0^*\otimes M_0(F_x))$, where $\phi_x=\phi_{0,\cU,x}$ for all but finite number of $x\in X$.  

 In particular, we have a canonical diagram
$$
\begin{array}{ccc}
 \cS(U_0^*\otimes M_0(F_x))^{(H\times G)(\cO_x)} & \hook{} &  \cS(U_0^*\otimes M_0(\AA))^{(H\times G)(\cO)}\\
\uparrow\lefteqn{\scriptstyle \zeta_x} & & \uparrow\lefteqn{\scriptstyle \zeta} \\
\cS(V_0\otimes L_0^*\otimes\Omega(F_x))^{(H\times G)(\cO_x)} & \hook{} & \cS(V_0\otimes L_0^*\otimes\Omega(\AA))^{(H\times G)(\cO)},
\end{array}
$$
where $\zeta_x$ is given by (\ref{map_zeta}) with $U_0\otimes L_0^*\otimes \Omega(\AA)$ replaced by $U_0\otimes L_0^*\otimes \Omega(F_x)$. 

 Set 
\begin{multline*}
\Weil_{G,H}(k)=\{(f_1,f_2)\mid 
f_1\in \cS(V_0\otimes L_0^*\otimes\Omega(F_x))^{(H\times G)(\cO_x)}, \\
f_2\in  \cS(U_0^*\otimes M_0(F_x))^{(H\times G)(\cO_x)}\;\;
\mbox{such that}\;\; \zeta_x(f_1)=f_2\}
\end{multline*}

  The Hecke property of $\phi_0$ (a classical analogue of Theorem~\ref{Th_Hecke_property_Aut_GH}) is as follows. 
\begin{Pp}
\label{Pp_classical_analogue_Th_Main}
 1) Assume $m\le n$. There is a homomorphism $\kappa: \cH_x(G)\to \cH_x(H)$ such that for $h\in \cH_x(G)$ we have 
$$
_x\H^{\la}_G(h, \phi_0)={_x\H^{\ra}_H}(\kappa(h), \phi_0)
$$
2) Assume $m>n$. There is a homomorphism $\kappa: \cH_x(H)\to \cH_x(G)$ such that for $h\in \cH_x(H)$ we have 
$$
_x\H^{\la}_H(h, \phi_0)={_x\H^{\ra}_G}(\kappa(h), \phi_0)
$$
\end{Pp}

 The above discussion reduces the proof of Proposition~\ref{Pp_classical_analogue_Th_Main} to the following local result.
 
\begin{Pp}  
\label{Pp_local_main_classical}
1) Assume $m\le n$. There is a homomorphism $\kappa: \cH_x(G)\to \cH_x(H)$ such that for $h\in \cH_x(G)$ we have 
$$ 
_x\H^{\la}_G(h, \phi_{0,\cU,x})=\zeta_x( {_x\H^{\ra}_H}(\kappa(h), \phi_{0,\cL,x}))
$$
Moreover, $\Weil_{G,H}(k)$ is a free module of rank one over $\cH_x(H)$ generated by $\phi_{0,\cL,x}$. \\
2) Assume $m>n$. There is a homomorphism $\kappa: \cH_x(H)\to \cH_x(G)$ such that for $h\in \cH_x(H)$ we have 
$$
\zeta_x({_x\H^{\la}_H}(h, \phi_{0,\cL, x}))={_x\H^{\ra}_G}(\kappa(h), \phi_{0,\cU,x})
$$
Moreover, $\Weil_{G,H}(k)$ is a free module of rank one over $\cH_x(G)$ generated by $\phi_{0,\cU,x}$. 
\end{Pp}

 To the author's best knowledge, there are three different prooofs of Proposition~\ref{Pp_local_main_classical} available in the literature. First part of both statements 1) and 2) is proved by Rallis  (\cite{R}) by some explicit calculation based on the following description of the Jacquet module. By (\cite{Ku2}, Lemma~5.1) we have an isomorphism of $\GL(L_0)(F_x)\times H(F_x)$-representations
\begin{equation}
\label{Jacquet_module_iso} 
 \cS(V_0\otimes L_0^*\otimes\Omega(F_x))_{U(G)(F_x)}
\,\iso\, \cS(\Cr(V_0\otimes L_0^*\otimes\Omega(F_x)))
\end{equation}
Here $\Cr(V_0\otimes L_0^*\otimes\Omega(F_x))\subset V_0\otimes L_0^*\otimes\Omega(F_x)$ is the subset of maps $v: L_0(F_x)\to V_0\otimes\Omega(F_x)$ such that 
$s_{\cL}(v)=0$, where $s_{\cL}(v)$ denotes the composition
$$
s_{\cL}(v):\Sym^2 L_0(F_x)\toup{\Sym^2 v} \Sym^2 (V_0\otimes\Omega(F_x))\to \Omega^2(F_x)
$$
A different proof due to Howe is found in \cite{MVW}, where the space $\Weil_{G,H}$ is described completely (a revisited version is given in \cite{H}). One more proof is given by Kudla in \cite{Ku}. Namely, in \cite{Ku2} it was shown that the Howe correspondence is compatible with the parabolic induction, this allows one to describe explicitely the image of a principal series representation under the Howe correspondence (cf. \cite{Ku}, Proposition~3.2, p.96), hence, to derive the functoriality (\cite{Ku}, Theorem on p. 105).

\medskip\noindent
3.2.2  In Section~6 we prove Theorem~\ref{Th_main_local_Sp_SO}, which is a geometric analogue of Proposition~\ref{Pp_local_main_classical}. The main difficulty is that the existing proofs of proof Proposition~\ref{Pp_local_main_classical} do not geometrize in an obvious way. Our approach, though inspired by \cite{R}, is somewhat different.   
  
 One more feature is that classical proofs of Proposition~\ref{Pp_local_main_classical} do not reveal a relation with the $\SL_2$ of Arthur, though it is believed to be relevant here (cf. the conjectures of Adams in \cite{Ad}). In our approach at least the maximal torus of $\SL_2$ of Athur appears naturally. In Section~8 we derive Theorem~\ref{Th_Hecke_property_Aut_GH} from Theorem~\ref{Th_main_local_Sp_SO}. 
 
\bigskip

{\centerline{\scshape 4. Geometric model of  the Schwarz space and Hecke functors}}

\bigskip\noindent
4.1 Set $\cO=k[[t]]\subset F=k((t))$, write $D^*=\Spec F\subset D=\Spec\cO$. Let $\Omega$ be the completed module of relative differentials of $\cO$ over $k$. 

 For a free $\cO$-module $M$ of finite rank we introduce the categories $P(M(F))\subset \D(M(F))$ as follows. For $N,r\ge 0$ set $_{N,r}M=t^{-N}M/t^rM$. Given positive integers $N_1\ge N_2$, $r_1\ge r_2$ 
we have a cartesian diagram
\begin{equation}
\label{diag_subspaces_limit}
\begin{array}{ccc}
_{N_2,r_1}M & \hook{i} & _{N_1, r_1}M\\
\downarrow\lefteqn{\scriptstyle  p} && \downarrow\lefteqn{\scriptstyle  p}\\
_{N_2,r_2}M & \hook{i} & _{N_1,r_2}M,
\end{array}
\end{equation}
where $i$ is the natural closed immersion, and $p$ is the projection.  

 By (\cite{G}, Lemma~4.8), the functor $\D(_{N,r_2}M)\to 
\D(_{N,r_1}M)$ given by $K\mapsto p^*K[\dimrel (p)]$ 
is fully faithful and exact for the perverse t-structures, and similarly for the functor $i_*$. These functors yield a diagram of full triangulated subcategories 
\begin{equation}
\label{diag_categories_limit}
\begin{array}{ccc}
\D(_{N_2,r_1}M) & \hook{} & \D(_{N_1,r_1}M)\\
\uparrow && \uparrow\\
\D(_{N_2,r_2}M) & \hook{} & \D(_{N_1,r_2}M)
\end{array}
\end{equation}
We let $\D(M(F))\,\iso\, \colim_{r,N}  \D(_{N,r}M)$ in $\DGCat_{cont}$. Then $\P(M(F))\subset \D(M(F))$ is the heart of the perverse t-structure. 

The category $\P(M(F))$ is a geometric analog of the Schwarz space of locally constant functions with compact support on $M(F)$.

 Set $_NM=t^{-N}M$ viewed as a $k$-scheme (not of finite type). 

\medskip\noindent
4.2.1 Let $G$ be a connected reductive group over $k$, assume that $M=M_0\otimes_k \cO$, where $M_0$ is a given finite-dimensional representation of $G$. 

 For $N+r>0$ the group $G(\cO)$ acts on $_{N,r}M$ via its finite-dimensional quotient $G(\cO/t^{N+r}\cO)$. For $r_1\ge N+r>0$ the kernel of $G(\cO/t^{r_1}\cO)\to G(\cO/t^{N+r}\cO)$ is a contractible unipotent group. So, the projection between the stack quotients
$$
q: G(\cO/t^{r_1}\cO)\backslash {_{N,r}M}\to
G(\cO/t^{N+r}\cO)\backslash {_{N,r}M}
$$
yields an (exact for the perverse t-structures) equivalence  
\begin{equation}
\label{trans_functor_for_4.2.1}
\D_{G(\cO/t^{N+r}\cO)}(_{N,r}M)\to \D_{G(\cO/t^{r_1}\cO)}(_{N,r}M)
\end{equation}
in $\DGCat_{cont}$. Denote by $\D_{G(\cO)}({_{N,r}M})$ the $\DG$-category $\D_{G(\cO/t^{r_1}\cO)}({_{N,r}M})$ for any $r_1\ge N+r$. 

 The stack quotient of (\ref{diag_subspaces_limit}) by $G(\cO/t^{N_1+r_1}\cO)$ yields a diagram
\begin{equation}
\label{diag_transition_functors_M(F)}
\begin{array}{ccc}
\D_{G(\cO)}(_{N_2,r_1}M) & \hook{} & \D_{G(\cO)}(_{N_1,r_1}M)\\
\uparrow && \uparrow\\
\D_{G(\cO)}(_{N_2,r_2}M) & \hook{} & \D_{G(\cO)}(_{N_1,r_2}M),
\end{array}
\end{equation}
where each arrow is a fully faithfull (and exact for the perverse t-structures) functor. Let $\D_{G(\cO)}(M(F))=\colim_{N,r} \D_{G(\cO)}(_{N,r}M)$ in $\DGCat_{cont}$. 

 Since $G(\cO/t^{N+r}\cO)$ is connected, the category $\P_{G(\cO)}(_{N,r}M)$ of $G(\cO/t^{N+r}\cO)$-equivariant perverse sheaves on $_{N,r}M$ is a full subcategory of $\P(_{N,r}M)$. The category $\P_{G(\cO)}(M(F))$ is defined as the heart of the perverse t-structure on $\D_{G(\cO)}(M(F))$. 
 
 For $N,r\ge 0, r_1\ge N+r$ we have the full subcategory 
$$
\D_{G(\cO/t^{r_1}\cO)}({_{N,r}M})^{constr}\subset \D_{G(\cO/t^{r_1}\cO)}({_{N,r}M})
$$ 
defined as in Section~2.1.1, as the stack $G(\cO/t^{r_1}\cO)\backslash {_{N,r}M}$ has an affine diagonal. The functors (\ref{trans_functor_for_4.2.1}) identify these categories for all $r_1\ge N+r$, so we get a full subcategory $\D_{G(\cO)}({_{N,r}M})^{constr}\subset \D_{G(\cO)}({_{N,r}M})$. The transition functors (\ref{diag_transition_functors_M(F)}) preserve the constructible subcategories. So, as in Section~2.1.2 we define $\D_{G(\cO)}(M(F))^{constr}=\colim_{N,r} \D_{G(\cO)}({_{N,r}M})^{constr}$ in $\DGCat^{non-cocmpl}$. Then $\D_{G(\cO)}(M(F))^{constr}$ is a full subcategory of $\D_{G(\cO)}(M(F))$. 
 
 Since the Verdier duality is compatible with the transition functors in (\ref{diag_transition_functors_M(F)}), it yields an equivalence 
$$
\DD: (\D_{G(\cO)}(M(F))^{constr})^{op}\,\iso\, \D_{G(\cO)}(M(F))^{constr}
$$ 
By Section~2.1.2, we also have the equivalence $\DD: (\D(M(F))^{constr})^{op}\,\iso\, \D(M(F))^{constr}$. 

 Recall also that, by (\cite{AGKRRV}, Th. C.2.6), each stack $G(\cO/t^{r_1}\cO)\backslash {_{N,r}M}$ is duality adapted in the sense of (\cite{AGKRRV}, C.2.4). So, the Verdier duality also gives an equivalence 
$$
\DD: (\D_{G(\cO/t^{r_1}\cO)}(_{N,r}M)^c)^{op}\,\iso\, \D_{G(\cO/t^{r_1}\cO)}(_{N,r}M)^c,
$$ 
for $r_1\ge N+r$. Besides, $\D_{G(\cO)}(_{N,r}M)$ is compactly generated. 
 
 The transition functors (\ref{diag_transition_functors_M(F)}) also preserve compact objects (their right adjoints are continuous), so we may consider $\colim_{N,r} \D_{G(\cO)}(_{N,r}M)^c$ taken in $\DGCat^{non-cocmpl}$, and 
$$
\Ind(\mathop{\colim}_{N,r} \; \D_{G(\cO)}(_{N,r}M)^c)\,\iso\, \D_{G(\cO)}(M(F))
$$ 
canonically in $\DGCat_{cont}$. In partcular, $\D_{G(\cO)}(M(F))$ is compactly generated.
 
\medskip\noindent
4.2.2  Write $\Gr_G$ for the affine grassmanian $G(F)/G(\cO)$ of $G$. Let us define the $\DG$-category $\D_{G(\cO)}(M(F)\times\Gr_G)$. 

 For $s_1, s_2\ge 0$ let 
$$
_{s_1,s_2}G(F)=\{g\in G(F)\mid t^{s_1}M\subset gM\subset t^{-s_2}M\},
$$ 
it is stable by left and right multiplication by $G(\cO)$, and $_{s_1,s_2}\Gr_G:=(_{s_1,s_2}G(F))/G(\cO)$ is closed in $\Gr_G$. For $s'_1\ge s_1$, $s'_2\ge s_2$ we have a closed embedding $_{s_1,s_2}\Gr_G\hook{} {_{s'_1,s'_2}\Gr_G}$, and the union of all $_{s_1,s_2}\Gr_G$ is $\Gr_G$. The map $g\mapsto g^{-1}$ yields an isomorphism $_{s_1,s_2}G(F)\,\iso\, {_{s_2, s_1}G(F)}$.

Assume for simplicity that $M_0$ is a faithful $G$-module, then the action of $G(\cO)$ on $_{s_1,s_2}\Gr_G$ factors through an action of $G(\cO/t^{s_1+s_2})$. 

 For $N,r, s_1,s_2\ge 0$ and $s\ge\max\{N+r, s_1+s_2\}$ we have the $\DG$-category
$$
\D_{G(\cO/t^s)}(_{N,r}M\times {_{s_1,s_2}\Gr_G}),
$$
where the action of $G(\cO/t^s)$ on $_{N,r}M\times {_{s_1,s_2}\Gr_G}$ is diagonal. For $s'\ge s\ge \max\{N+r, s_1+s_2\}$ we have a canonical equivalence (exact for the perverse t-structures)
$$
\D_{G(\cO/t^s)}(_{N,r}M\times {_{s_1,s_2}\Gr_G})\,\iso\, 
\D_{G(\cO/t^{s'})}(_{N,r}M\times {_{s_1,s_2}\Gr_G})
$$
Define $\D_{G(\cO)}(_{N,r}M\times {_{s_1,s_2}\Gr_G})$ as the category $\D_{G(\cO/t^s)}(_{N,r}M\times {_{s_1,s_2}\Gr_G})$ for any $s$ as above. As in Section~4.2.1, set 
$$
\D_{G(\cO)}(M(F)\times\Gr_G)=\colim_{N,r, s_1, s_2} \D_{G(\cO)}(_{N,r}M\times {_{s_1,s_2}\Gr_G})
$$ 
in $\DGCat_{cont}$. The subcategory of perverse sheaves 
$$
\P_{G(\cO)}(M(F)\times\Gr_G)\subset \D_{G(\cO)}(M(F)\times\Gr_G)
$$ 
is defined as the heart of the perverse t-structure.  

 One defines $\D_{G(\cO)}(\Gr_G)^{constr}\in\DGCat^{non-cocmpl}$ as in Section~4.2.1. Namely, let us write $\Gr_G\,\iso\, \colim_{i\in I} Z_i$, where $Z_i\subset \Gr_G$ is a $G(\cO)$-invariant closed subscheme of finite type, $I$ is a small filtered category, and for $i\to j$ in $I$ the map $Z_i\to Z_j$ is a closed immersion. By definition, $\D_{G(\cO)}(\Gr_G)^{constr}=\colim_{i\in I} \D_{G(\cO)}(Z_i)^{constr}$ taken in $\DGCat^{non-cocmpl}$. Here $\D_{G(\cO)}(Z_i)^{constr}$ is defined as $\D_{G(\cO/t^r\cO)}(Z_i)^{constr}$ for $r$ large enough. 
Then $\D_{G(\cO)}(\Gr_G)^{constr}\subset \D_{G(\cO)}(\Gr_G)$ is a full subcategory.  
 
\medskip\noindent
4.3.1 Let $\Sph_G$ denote the category of spherical perverse sheaves on $\Gr_G$. Recall the canonical equivalence of tensor categories $\Loc: \Rep(\check{G})\,\iso\,\Sph_G$ (cf. \cite{MV}). 

 Let us define an action of the tensor category $\Sph_G$ on $\D_{G(\cO)}(M(F))$ by Hecke functors. Consider the map 
\begin{equation}
\label{map_alpha_for_def_q_act}
\alpha: M(F)\times G(F)\,\iso\, M(F)\times G(F)
\end{equation}
sending $(m,g)$ to $(g^{-1}m, g)$. Let $(a,b)\in G(\cO)\times G(\cO)$ act on the source sending $(m,g)$ to $(am, agb)$. Let it also act on the target sending $(m',g')$ to $(b^{-1}m', ag'b)$. The above map is equivariant for these actions, so yields a morphism of stacks
$$
_q\act: G(\cO)\backslash (M(F)\times \Gr_G)\to (M(F)/G(\cO))\times (G(\cO)\backslash \Gr_G),
$$ 
where the action of $G(\cO)$ on $M(F)\times \Gr_G$ is the diagonal one. 

 The connected components of $\Gr_G$ are indexed by $\pi_1(G)$. For $\theta\in\pi_1(G)$ the component 
$\Gr_G^{\theta}$ is the one containing $\Gr_G^{\lambda}$ for any $\lambda\in\Lambda^+_G$ whose image in $\pi_1(G)$ equals $\theta$. For $\theta\in\pi_1(G)$ set $_{s_1,s_2}\Gr_G^{\theta}=\Gr_G^{\theta}\cap {_{s_1,s_2}\Gr_G}$. 

 In the rest of Section~4.3.1 we construct an inverse image functor
\begin{equation}
\label{functor_q_act}
_q\act^*(\cdot,\cdot): \D_{G(\cO)}(M(F))\times \D_{G(\cO)}(\Gr_G)\to \D_{G(\cO)}(M(F)\times\Gr_G)
\end{equation}
satisfying the following properties.
\begin{itemize}
\item[A1)] For $K\in\D_{G(\cO)}(M(F))^{constr}$, $\cT\in \D_{G(\cO)}(\Gr_G)^{constr}$ one has canonically
\begin{equation}
\label{inverse_image_for_Hecke_comm_Verdier}
\DD(_q\act^*(K,\cT))\,\iso\, {_q\act^*}(\DD(K),\DD(\cT))
\end{equation}
\item[A2)] If both $K$ and $\cT$ are perverse then $_q\act^*(K,\cT)$ is perverse.
\end{itemize}

 For non negative integers $N, r, s_1,s_2$, with $r\ge s_1$ and $s\ge \max\{s_1+s_2, N+r\}$ we have a diagram
$$
\begin{array}{ccccc}
&& _{N,r}M \times {_{s_1,s_2}G(F)} & \toup{\act} & _{N+s_1, r-s_1}M\\
&& \downarrow\lefteqn{\scriptstyle q_G} && \downarrow\lefteqn{\scriptstyle q_M}\\
_{N,r}M & \getsup{\pr} & _{N,r}M\times  {_{s_1,s_2}\Gr_G} & \toup{\act_q} & G(\cO/t^s)\backslash {_{N+s_1, r-s_1}M}\\
\downarrow && \downarrow & \nearrow\lefteqn{\scriptstyle \act_{q,s}}\\
G(\cO/t^s)\backslash {_{N,r}M} & \getsup{\pr} & G(\cO/t^s)\backslash (_{N,r}M \times {_{s_1,s_2}\Gr_G}) & \toup{\pr_2} &G(\cO/t^s)\backslash (_{s_1,s_2}\Gr_G)
\end{array}
$$
Here $\act$ sends $(m, g)$ to $g^{-1}m$, the map $q_G$ sends $(m,g)$ to $(m, gG(\cO))$, and $\pr, \pr_2$ denote the projections. All the vertical arrows are the stack quotient maps for the action of a corresponding group. One checks that $\act$ descends to a map $\act_q$ between the corresponding quotients.

  For $s\ge \max\{s_1+s_2, N+r\}$ the group $G(\cO/t^s\cO)$ acts diagonally on $_{N,r}M\times  {_{s_1,s_2}\Gr_G}$, and $\act_q$ is equivariant with respect to this action. Consider the functors 
$$
\begin{array}{ccc}
\D_{G(\cO/t^s)}(_{N+s_1, r-s_1}M)\times \D_{G(\cO/t^s)}(_{s_1,s_2}\Gr_G) & \toup{\temp} & \D_{G(\cO/t^s)}(_{N,r}M \times {_{s_1,s_2}\Gr_G})\\
\parallel && \parallel\\
\D_{G(\cO)}(_{N+s_1, r-s_1}M)\times \D_{G(\cO)}(_{s_1,s_2}\Gr_G) && \D_{G(\cO)}(_{N,r}M \times {_{s_1,s_2}\Gr_G})
\end{array}
$$
sending $(K,\cT)$ to 
\begin{equation}
\label{def_act_q_well_shifted}
\act_{q,s}^*(K)\otimes \pr_2^*\cT[s\dim G+s_1\dim M_0-c],
\end{equation}
where $c$ equals to $\<\theta,\check{\mu}\>$ over $_{s_1,s_2}\Gr_G^{\theta}$. Here $\check{\mu}\in\check{\Lambda}^+_G$ denotes the character $\det M_0$. 

 For $r_1\ge r_2$ and $s\ge\max\{s_1+s_2, N+r_1\}$ the functors $\temp$ are compatible with the transition functors for the diagram
$$
\begin{array}{ccc}
G(\cO/t^s)\backslash (_{N,r_1}M \times {_{s_1,s_2}\Gr_G}) & \toup{\act_{q,s}} & G(\cO/t^s)\backslash {_{N+s_1, r_1-s_1}M}\\
\downarrow && \downarrow\\
G(\cO/t^s)\backslash (_{N,r_2}M \times {_{s_1,s_2}\Gr_G}) & \toup{\act_{q,s}} & G(\cO/t^s)\backslash {_{N+s_1, r_2-s_1}M}
\end{array}
$$
So, they yield a functor
$$
_{N,s_1,s_2}\temp: \D_{G(\cO)}(_{N+s_1}M)\times \D_{G(\cO)}(_{s_1,s_2}\Gr_G)\to \D_{G(\cO)}(_NM\times_{s_1,s_2}\Gr_G)
$$

 For $N_1\ge N_2$ and $s\ge \max\{s_1+s_2, N_1+r\}$ 
we have a diagram, where the vertical maps are closed immersions
$$
\begin{array}{ccc}
G(\cO/t^s)\backslash (_{N_1,r}M \times {_{s_1,s_2}\Gr_G}) & \toup{\act_{q,s}} & 
G(\cO/t^s)\backslash {_{N_1+s_1, r-s_1}M}\\
\uparrow && \uparrow\\
G(\cO/t^s)\backslash (_{N_2,r}M \times {_{s_1,s_2}\Gr_G}) & \toup{\act_{q,s}} & 
G(\cO/t^s)\backslash {_{N_2+s_1, r-s_1}M}
\end{array}
$$
This diagram is not cartesian in general, we come around this as follows. If $K\in\D_{G(\cO)}(_NM)$, $\cT\in\D_{G(\cO)}(_{s_1,s_2}\Gr_G)$ then for any $N_1\ge N+s_2$ the image of $(K,\cT)$ under the composition
\begin{multline*}
\D_{G(\cO)}(_NM)\times \D_{G(\cO)}(_{s_1,s_2}\Gr_G)\subset \D_{G(\cO)}(_{N_1+s_1}M)\times \D_{G(\cO)}(_{s_1,s_2}\Gr_G) \;\toup{_{N_1,s_1,s_2}\temp} \\ \D_{G(\cO)}(_{N_1}M\times{_{s_1,s_2}\Gr_G}) \subset
\D_{G(\cO)}(M(F)\times{_{s_1,s_2}\Gr_G})
\end{multline*}
does not depend on $N_1$, so we get a functor
$$
_{s_1,s_2}\temp: \D_{G(\cO)}(M(F))\times \D_{G(\cO)}(_{s_1,s_2}\Gr_G)\to \D_{G(\cO)}(M(F)\times_{s_1,s_2}\Gr_G)
$$

 For $s'_1\ge s_1, s'_2\ge s_2$ we have the functors 
of extension by zero 
$$
\D_{G(\cO)}(_{s_1,s_2}\Gr_G)\to \D_{G(\cO)}(_{s'_1,s'_2}\Gr_G)
$$ 
They are compatible with $_{s_1,s_2}\temp$. This yields the desired functor (\ref{functor_q_act}). The properties A1, A2 follow easily from the construction.

\medskip\noindent
4.3.2 For nonnegative integers $s_1,s_2,N,r$ and $s\ge \max\{N+r, s_1+s_2\}$ for the projection 
$$
\pr: G(\cO/t^s)\backslash (_{N,r}M\times_{s_1,s_2}\Gr_G)\to  G(\cO/t^s)\backslash {_{N,r}M}
$$ 
the corresponding functors $\pr_!:\D_{G(\cO)}(_{N,r}M\times_{s_1,s_2}\Gr_G)\to \D_{G(\cO)}(_{N,r}M)$ are compatible with the transition functors, so yield a functor
$\pr_!: \D_{G(\cO)}(M(F)\times\Gr_G)\to \D_{G(\cO)}(M(F))$. 
 
  Finally, we define the Hecke functor
\begin{equation}
\label{Hecke_functor_local_final}
\H^{\la}_G(\cT,\cdot):  \D_{G(\cO)}(M(F))\to \D_{G(\cO)}(M(F))
\end{equation}
by $\H^{\la}_G(\cT,K)=\pr_!(_q\act^*(K,\cT))$ for $\cT\in\Sph_G$ and $K\in \D_{G(\cO)}(M(F))$. 
 From (\ref{inverse_image_for_Hecke_comm_Verdier}) one gets that the functors (\ref{Hecke_functor_local_final}) commute with the Verdier duality, namely 
$$
\DD(\H^{\la}_G(\cT,K))\,\iso\, \H^{\la}_G(\DD\cT, \DD K)
$$ 
They are also compatible with the tensor structure on $\Sph_G$ (as in \cite{BG}, Section~3.2.4).  For $\cT\in\Sph_G$ and $K\in \D_{G(\cO)}(M(F))$ set $\H^{\ra}_G(\cT,K)=\H^{\la}_G(\ast \cT, K)$. Then the functors 
$$
K\mapsto \H^{\la}_G(\cT,K)\;\;\;\mbox{and}\;\;\;
K\mapsto \H^{\ra}_G(\DD(\cT),K)
$$
are mutually (both left and right) adjoint. 

 For a $G$-dominant coweight $\lambda$ we set $\H^{\lambda}_G(\cdot)=\H^{\la}_G(\cA^{\lambda}_G,\cdot)$. 
 
\begin{Rem} 
\label{Rem_true_measures}
i) Call $K\in \P_{G(\cO)}(_RM)$ \select{smooth} if it comes from a $G(\cO)$-equivariant local system on $_{R,r}M$ for some $r$. Let us make the above definition explicit in this case. 

 Let us above $\check{\mu}\in\check{\Lambda}^+_G$ denote the character $\det M_0$, so the virtual dimension $\dim(M/gM)=\<\theta, \check{\mu}\>$ for $gG(\cO)\in\Gr_G^{\theta}$. Let $\cT\in \Sph_G$ be the extension by zero from $_{s_1,s_2}\Gr_G$. For $r$ large enough, let $_{R,r}M\ttimes {_{s_1,s_2}\Gr_G}$ be the scheme of pairs $(m, gG(\cO))$ with $gG(\cO)\in {_{s_1,s_2}\Gr_G}$ and $m\in t^{-R}gM/t^rM$. Set 
$$ _{s_1,s_2}\Gr_G^{\theta}={_{s_1,s_2}\Gr_G}\cap\Gr_G^{\theta}$$ 
Then $_{R,r}M\ttimes {_{s_1,s_2}\Gr^{\theta}_G}$ is a locally trivial fibration over $_{s_1,s_2}\Gr^{\theta}_G$ with fibre an affine space of dimension $(R+r)\dim M_0-\<\theta,\check{\mu}\>$. We get a diagram
$$
{_{R+s_2, r}M}\getsup{\pr} {_{R,r}M\ttimes {_{s_1,s_2}\Gr_G}}\toup{\act_q} G(\cO/t^{R+r-s_1})\backslash(_{R,r-s_1}M),
$$ 
where $\pr$ sends $(m, gG(\cO))$ to $m$. Let $K\tboxtimes\cT$ denote the perverse sheaf $\act_q^*K\otimes\pr_2^*\cT[\dim]$ on ${_{R,r}M\ttimes {_{s_1,s_2}\Gr_G}}$, here $\dim$ is the unique integer for which this complex is perverse. Then $\H^{\la}_G(\cT, K)=\pr_!(K\tboxtimes\cT)$. We see once again that indeed the shift in (\ref{def_act_q_well_shifted}) must depend on $\check{\mu}$. 

\medskip\noindent
ii) At the level of functions the formula (\ref{Hecke_functor_local_final}) 
reads as follows. For $f: M(F)\to \Qlb$, which is $G(\cO)$-invariant, $\eta: \Gr_G\to \Qlb$, which is of compact support and $G(\cO)$-invariant one gets $\H^{\la}_G(\eta, f): M(F)\to \Qlb$ given by
$$
m\mapsto \int_{gG(\cO)\in\Gr_G} f(g^{-1}m)\eta(gG(\cO))dg,
$$
where $dg$ is a Haar measure on $G(F)$ such that the volume of $G(\cO)$ is one.
\end{Rem}  
 
\medskip\noindent
4.4 Let $I_0$ denote the constant sheaf $\Qlb$ on $_{0,0}M$, this is an object of $P_{G(\cO)}(M(F))$. For $K=I_0$ 
the above definition of $\H^{\lambda}_G(K)$  
simplifies as follows.
 
Assume that all the weights of $M_0$ are less or equal to a $G$-dominant weight $\check{\lambda}$. Then for a dominant coweight $\lambda$ of $G$ we have $\ov{\Gr}_G^{\lambda}\subset {_{r,N}\Gr_G}$ and 
$\H^{\lambda}_G(I_0)\in \D_{G(\cO)}(_{N,r}M)$ with $r=\<\lambda,\check{\lambda}\>$ and $N=\<-w_0(\lambda), \check{\lambda}\>$. 

Let $M\ttimes \ov{\Gr}_G^{\lambda}$ be the scheme classifying pairs $gG(\cO)\in \ov{\Gr}_G^{\lambda}$, $m\in gM$. Let $\pi_{\pro}: M\ttimes \ov{\Gr}_G^{\lambda}\to {_NM}$ be the map sending $(m, gG(\cO))$ to $m\in {_NM}$. 
Write $_{0,r}M\ttimes \ov{\Gr}_G^{\lambda}$ for the scheme classifying $gG(\cO)\in \ov{\Gr}_G^{\lambda}$, $m\in gM/t^rM$. The map $\pi_{pro}$ gives rise to a proper map
\begin{equation}
\label{map_pi_rem2}
\pi: {_{0,r}M}\ttimes \ov{\Gr}_G^{\lambda}\to {_{N,r}M}
\end{equation}
sending $(m, gG(\cO))$ to $m$. By Remark~\ref{Rem_true_measures}, $\H^{\lambda}_G(I_0)$ identifies with $\pi_!(\Qlb\tboxtimes \cA^{\lambda}_G)[\dim {_{0,r}M}]$. 
 
\medskip\noindent 
4.5.1 The group of automorphisms of the $k$-algebra $\cO$ is naturally the group of $k$-points of a (reduced) affine group scheme $\Aut^0\cO$ over $k$. The group scheme $\Aut^0\cO$ acts naturally on $M(F)$, $G(F)$, $G(\cO)$ and $\Gr_G$. We write $\delta: \Aut^0\cO\times M(F)\to M(F)$ and $\delta: \Aut^0\cO\times G(F)\to G(F)$ for the corresponding action maps, and 
$G(F)\rtimes\Aut^0\cO$ for the corresponding semi-direct product with operation 
$$
(g_1,c_1)(g_2,c_2)=(g_1\delta(c_1, g_2), c_1c_2),\;\;\;  c_i\in \Aut^0\cO, g_i\in G(F).
$$ 
Then $G(F)\rtimes\Aut^0\cO$ acts on $M(F)$ via the map $(G(F)\rtimes\Aut^0\cO)\times M(F)\to M(F)$ sending $((g,c), m)$ to $g\delta(c,m)$. For $r>0$ we similarly have a semi-direct product $G(\cO/t^r)\rtimes \Aut(\cO/t^r)$ and a surjective homomorphism $G(\cO)\rtimes\Aut(\cO)\to G(\cO/t^r)\rtimes \Aut(\cO/t^r)$, whose kernel is pro-unipotent. 

 Let us define the $\DG$-category $\D_{G(\cO)\rtimes\Aut^0\cO}(M(F))$. As in 4.2, for $r_1\ge N+r>0$ the projection between the stack quotients
$$
q: (G(\cO/t^{r_1})\rtimes\Aut(\cO/t^{r_1}))\backslash {_{N,r}M}
\to (G(\cO/t^{N+r})\rtimes\Aut(\cO/t^{N+r}))\backslash {_{N,r}M}
$$
yields an (exact for the perverse t-structures) equivalence 
$$
\D_{G(\cO/t^{N+r})\rtimes\Aut(\cO/t^{N+r})}(_{N,r}M)\to
\D_{G(\cO/t^{r_1})\rtimes\Aut(\cO/t^{r_1})}(_{N,r}M)
$$
in $\DGCat_{cont}$. Denote by $\D_{G(\cO)\rtimes\Aut^0\cO}(_{N,r}M)$ the $\DG$-category $\D_{G(\cO/t^{r_1})\rtimes\Aut(\cO/t^{r_1})}(_{N,r}M)$ for any $r_1\ge N+r$. 

 The stack quotient of (\ref{diag_subspaces_limit}) by $G(\cO/t^{N_1+r_1})\rtimes\Aut(\cO/t^{N_1+r_1})$ yields a diagram
\begin{equation}
\label{diag_transition_functors_M(F)_Gm}
\begin{array}{ccc}
\D_{G(\cO)\rtimes\Aut^0\cO}(_{N_2,r_1}M) & \hook{} & \D_{G(\cO)\rtimes\Aut^0\cO}(_{N_1,r_1}M)\\
\uparrow && \uparrow\\
\D_{G(\cO)\rtimes\Aut^0\cO}(_{N_2,r_2}M) & \hook{} & \D_{G(\cO)\rtimes\Aut^0\cO}(_{N_1,r_2}M),
\end{array}
\end{equation}
where each arrow is a fully faithfull (and exact for the perverse t-structures) functor. Let 
$$
\D_{G(\cO)\rtimes\Aut^0\cO}(M(F))=\colim_{N,r} \D_{G(\cO)\rtimes\Aut^0\cO}(_{N,r}M)
$$ 
in $\DGCat_{cont}$. Similarly, the category of perverse sheaves 
$$
\P_{G(\cO)\rtimes\Aut^0\cO}(M(F))\subset \D_{G(\cO)\rtimes\Aut^0\cO}(M(F))
$$
is the heart of the corresponding t-structure. 

 As in 4.3, one defines a natural action of $\Sph_G$ on $\D_{G(\cO)\rtimes\Aut^0\cO}(M(F))$. For our purposes note that the map (\ref{map_pi_rem2}) is $\Aut^0\cO$-equivariant, so that all the perverse cohomologies of $\H^{\lambda}_G(I_0)$ are objects of $\P_{G(\cO)\rtimes\Aut^0\cO}(M(F))$. 
 
\medskip\noindent
4.5.2 For the projection $q_M: M(F)\to\Spec k$ define $(q_M)_!, (q_M)_*: \D(M(F))\to \D(\Spec k)$ as follows. Recall that in the definition of $\D(M(F))$ we use the transition functors $\phi^*[\dimrel \phi]$ for the natural map $\phi: {_{N,r'}M}\to {_{N,r}M}$. For $N\le N'$ and the closed embeddings $i: {_{N,r}M}\hook{} {_{N', r}M}$ we also use the transition functor $i_*$. 

 For $N, r\ge 0$ consider the functor $\D( {_{N,r}M})\to\D(\Spec k)$ 
$$
K\mapsto \RG_c(_{N,r}M, K)[r\dim(M_0)]
$$ 
These functors are compatible with the above transition functors, so define the desired functor $(q_M)_!$ by passing to the colimit. 

 For $N, r\ge 0$ the functors $\D( {_{N,r}M})\to\D(\Spec k)$ 
$$
K\mapsto \RG(_{N,r}M, K)[-r\dim(M_0)]
$$ 
are compatible with the above transition functors, so define the desired functor $(q_M)_*$ by passing to the colimit.  

\medskip\noindent
4.5.3 Define the functors still denoted $(q_M)_!, (q_M)_*: \D_{G(\cO)}(M(F))\to D_{G(\cO)}(\Spec k)$ as follows.

 Recall that $\D_{G(\cO)}(M(F))$ is defined as 
$$
\mathop{\colim}\limits_{N\in\NN,r\in\NN} \D_{G(\cO)}(_{N,r}M)
$$ 
for the following transition functors. For $N\le N'$ we consider the closed immersions 
$i: G(\cO/t^s)\backslash{_{N,r}M}\to G(\cO/t^s)\backslash{_{N',r}M}$
with $s$ large enough and use $i_*$ as the transition functor. For $r\le r'$ we consider the natural map $\phi: G(\cO/t^s)\backslash{_{N,r'}M}\to G(\cO/t^s)\backslash{_{N,r}M}$ for $s$ large enough, and use $\phi^*[\dimrel(\phi)]$ as transition functors. 

 For $r,N\ge 0$ we have the maps  
$$
q_{N,r}:  G(\cO/t^s)\backslash {_{N,r}M}\to G(\cO/t^s)\backslash\Spec k
$$ 
with $s$ large enough. The functors
$\D_{G(\cO/t^s)}({_{N,r}M})\to \D_{G(\cO/t^s)}(\Spec k)$, 
$$
K\mapsto (q_{N,r})_![r\dim M_0]
$$
are compatible with the above transition functors, so yield the desired functor $(q_M)_!$. The functors
$\D_{G(\cO/t^s)}({_{N,r}M})\to \D_{G(\cO/t^s)}(\Spec k)$, 
$$
K\mapsto (q_{N,r})_*[-r\dim M_0]
$$
are compatible with the above transition functors, so yield the desired functor $(q_M)_*$. 

\medskip\noindent
4.5.4 Define the functor $\bar\epsilon: \D_{G(\cO)}(\Spec k)\to \D(\Spec k)$ exact for the perverse t-structures as follows. Recall that for $s'\ge s$ we have the surjection $G(\cO/t^{s'})\to G(\cO/t^s)$, hence a diagram $\Spec k\toup{v_{s'}} B(G(\cO/t^{s'})\toup{u} B(G(\cO/t^s))$ with $v_s=u v_{s'}$. The category  
$\D_{G(\cO)}(\Spec k)$ is defined as $\D_{G(\cO/t^s)}(\Spec k)$ for any $s>0$ with the transition functors
$$
u^*[\dimrel(u)]: \D_{G(\cO/t^s)}(\Spec k)\iso \D_{G(\cO/t^{s'})}(\Spec k)
$$ 
The functors $v_s^*[\dimrel(v_s)]: \D_{G(\cO/t^s)}(\Spec k)\to \D(\Spec k)$ are compatible with these transition functors, so define the functor $\bar\epsilon$ exact for the perverse t-structures. 

 Define $\bar q_{M, *}, \bar q_{M,!}: \D_{G(\cO)}(M(F))\to \D(\Spec k)$  by $\bar q_{M,!}=\bar\epsilon\comp (q_M)_!$, $\bar q_{M, *}=\bar\epsilon\comp(q_M)_*$. If $I_0$ is the constant perverse sheaf on $M$ then $\bar q_{M, *}I_0\,\iso\,\Qlb$ canonically. 
 
\medskip\noindent
4.5.5 The map (\ref{map_alpha_for_def_q_act}) gives rise to the commutative diagram of stack quotients
$$
\begin{array}{ccccc}
G(\cO)\backslash M(F) & \getsup{\pr_1} & G(\cO)\backslash (M(F)\times G(F)/G(\cO)) & \iso & (M(F)\times G(\cO)\backslash G(F))/G(\cO)\\
\downarrow\lefteqn{\scriptstyle q_M} && \downarrow\lefteqn{\scriptstyle\pr_2} & \swarrow\\
G(\cO)\backslash \Spec k & \gets & G(\cO)\backslash G(F)/G(\cO)
\end{array}
$$
The above diagram gives the following.

\begin{Lm} 1) For $K\in \D_{G(\cO)}(M(F)), \cS\in\Sph_G$ there is a canonical isomorphism in $\D(\Spec k)$
$$
\bar q_{M,*}(\H^{\la}_G(\cS, K))\,\iso\, \bar q_{M,*}(K) \otimes \RG(\Gr_G, \Res^{\kappa_{M*}}(\cS))
$$
where $\kappa_{M*}: \check{G}\times\Gm\to \check{G}$ is the map $(\id, \check{\mu})$. Here $\check{\mu}: G\to\Gm$ is the character $\det M_0$. \\
2) There is a canonical isomorphism in $\D(\Spec k)$
$$
\bar q_{M,!}(\H^{\la}_G(\cS, K))\,\iso\, \bar q_{M,!}(K)\otimes\RG(\Gr_G, \Res^{\kappa_{M !}}(\cS)),
$$
where $\kappa_{M !}: \check{G}\times\Gm\to \check{G}$ is the map $(\id, -\check{\mu})$. 
\QED
\end{Lm}
 
\medskip\noindent
4.6 If $X$ is a smooth projective connected curve and $x\in X$ then one can consider the following global version of the category $\D_{G(\cO)}(M(F))$. 

 Let $_{x,\infty}\cW_G$ be stack classifying a $G$-torsor $\cF_G$ on $X$ together with a section $\cO_X\toup{s} M_{\cF_G}(\infty x)$. The stack $_{x,\infty}\cW_G$ is an ind-algebraic. We have a diagram
$$
_{x,\infty}\cW_G  \getsup{h^{\la}_{\cW}} {_{x,\infty}\cW_G}\times_{\Bun_G} {_x\cH_G} \toup{h^{\ra}_{\cW}} {_{x,\infty}\cW_G},
$$
where we used $h^{\la}_G$ to define the fibre product, the map $h^{\la}_{\cW}$ (resp., $h^{\ra}_{\cW}$) sends $(\cF_G,\cF'_G,\beta, \cO_X\toup{s} M_{\cF_G}(\infty x))$ to $(\cF_G, s)$ (resp., to $(\cF'_G, s')$ with $s'=s\comp\beta$). 
As in Sections~2.2.1 and 4.3, one defines the Hecke functors 
$$
\H^{\la}_G(\cdot,\cdot), \H^{\ra}_G(\cdot,\cdot): \Sph_G\times\D(_{x,\infty}\cW_G)\to \D(_{x,\infty}\cW_G)
$$ 
 
 Let $_{x,\le N}\cW_G\subset {_{x,\infty}\cW_G}$ be the closed substack given by requiring that $\cO_X\toup{s}M_{\cF_G}(Nx)$ is regular.
 
  For $r\ge 1$ let $D_{r,x}=\Spec \cO_x/t_x^r$, where $\cO_x$ is the completed local ring at $x\in X$, and $t_x\in\cO_x$ is a local parameter. Pick a trivialization $\cO_x\,\iso\,\cO$. For $N,r\ge 0$ it yields a map 
$$
_{N,r}p_{\cW}: {_{x,\le N}\cW_G}\to G(\cO/t^{N+r})\backslash {_{N,r}M}
$$
sending $(\cF_G, \cO_X\toup{t} M_{\cF_G}(Nx))$ to $\cF_G\mid_{D_{N+r,x}}$ equipped with the induced $G(\cO/t^{N+r})$-equivariant map $\cF_G\mid_{D_{N+r,x}}\to {_{N,r}M}$. We get a functor $\D_{G(\cO/t^{N+r})}(_{N,r}M)\to \D(_{x,\le N}\cW_G)$ given by 
$$
K\mapsto {_{N,r}p_{\cW}^*K}[a+\dim\Bun_G+N\dim M_0-\dim G(\cO/t^{N+r})\backslash {_{N,r}M}],
$$ 
here $a$ is a function of a connected component of $\Bun_G$ sending $\cF_G$ to $\chi(M_{\cF_G})$. The shift in the above formula should be thought of as `the corrected relative dimension' of $_{N,r}p_{\cW}$, over a suitable open substack of $_{x,\le N}\cW_G$ it is indeed the relative dimension. These functors are compatible with the transition functors in (\ref{diag_transition_functors_M(F)}), thus we get a well-defined functor 
$$
\glob_x: \D_{G(\cO)}(M(F))\to \D(_{x,\infty}\cW_G),
$$
here $\glob$ stands for `globalization'. 
One checks that it commutes with the functors $\H^{\la}_G, \H^{\ra}_G$. Along the same lines, one defines a functor $\D_{G(\cO)\rtimes\Aut^0\cO}(M(F))\to \D(_{x,\infty}\cW_G)$ that does not depend on a choice of a trivialization $\cO_x\,\iso\,\cO$. 

 Write $_{\infty}\cW_G$ for the stack classifying $x\in X$, a $G$-torsor $\cF_G$ on $X$, and a section $s: \cO_X\to M_{\cF_G}(\infty x)$. As above, one defines the Hecke functors
$$
\H^{\la}_G(\cdot,\cdot), \H^{\ra}_G(\cdot,\cdot): \Sph_G\times\D(_{\infty}\cW_G)\to \D(_{\infty}\cW_G)
$$ 
Let $_{\le N}\cW_G\subset {_{\infty}\cW_G}$ be the closed substack given by the condition that $s: \cO_X\to M_{\cF_G}(N x)$ is regular. Along the same lines one gets a map $_{\le N}\cW_G\to (G(\cO/t^{N+r})\rtimes\Aut(\cO/t^{N+r}))\backslash {_{N,r}M}$, the corresponding functors 
$$
\D_{G(\cO/t^{N+r})\rtimes\Aut(\cO/t^{N+r})}(_{N,r}M)\to \D(_{\le N}\cW_G)
$$
are compatible with the transitions functors in (\ref{diag_transition_functors_M(F)_Gm}). The resulting functor 
$$
\glob_{\infty}: \D_{G(\cO)\rtimes \Aut^0(\cO)}(M(F))\to \D(_{\infty}\cW_G)
$$
commutes with the action of $\H^{\la}_G$, $\H^{\ra}_G$.


\bigskip\noindent
4.7 {\scshape Weak analogues of Jacquet functors} 

\medskip\noindent
4.7.1 Let $P\subset G$ be a parabolic subgroup, $U\subset P$ its unipotent radical and $L=P/U$ the Levi quotient. In classical setting, an important tool is the Jacquet module $\cS(M(F))_{U(F)}$ of coinvariants with respect to $U(F)$. We don't know how to geometrize the whole Jacquet module. However, let $V_0\subset M_0$ be a $P$-invariant subspace, on which $U$ acts trivially. Set $V=V_0(\cO)$. We have a surjective map of $L(F)$-representations $\cS(M(F))_{U(F)}\to \cS(V(F))$ given by restriction under $V(F)\hook{} M(F)$. We rather geometrize the composition map $\cS(M(F))\to \cS(M(F))_{U(F)}\to \cS(V(F))$. 
 
  As in 4.2, we have the $\DG$-categories $\D(V(F))$, $\D_{L(\cO)}(V(F))$. We are going to define natural functors 
\begin{equation}
\label{functors_J_P}
J_P^*, J_P^! : \D_{G(\cO)}(M(F))\to \D_{L(\cO)}(V(F))
\end{equation}
To do so, for $N+r\ge 0$ consider the natural closed embedding $i_{N,r}: {_{N,r}V}\hook{} {_{N,r}M}$. Recall that $_{N,r}V=t^{-N}V/t^r V$. 
Consider the diagram of stack quotients
\begin{equation}
\label{diag_Jacquet_functors_def}
\begin{array}{cc}
P(\cO/t^{N+r})\backslash (_{N,r}V) & \toup{i_{N,r}}\; P(\cO/t^{N+r})\backslash (_{N,r}M)
\toup{p} G(\cO/t^{N+r})\backslash (_{N,r}M)\\
\downarrow\lefteqn{\scriptstyle q}\\ 
L(\cO/t^{N+r})\backslash (_{N,r}V),
\end{array}
\end{equation}
where $p$ comes from the inclusion $P\subset G$ and $q$ is the natural quotient map. Using (\ref{diag_Jacquet_functors_def}), define functors 
$$
J_P^*, J_P^!: \D_{G(\cO/t^{N+r})}(_{N,r}M)\to \D_{L(\cO/t^{N+r})}(_{N,r} V)
$$ 
by 
$$
\begin{array}{c}
q^*\comp J_P^*[\dimrel(q)]=(i_{N,r})^* p^*[\dimrel(p)-ra]\\ \\
q^*\comp J_P^![\dimrel(q)]=(i_{N,r})^! p^*[\dimrel(p)+ra]
\end{array}
$$
Since $q^*[\dimrel(q)]: \D_{L(\cO/t^{N+r})}(_{N,r}V)\to \D_{P(\cO/t^{N+r})}(_{N,r}V)$ is an equivalence (exact for the perverse t-structures),
the functors $J_P^*, J_P^!$ are well-defined.  
Here we have set $a=\dim M_0-\dim V_0$. 

 Further, $J_P^*, J_P^!$ are compatible with the transition functors in (\ref{diag_transition_functors_M(F)}), so give rise to the desired functors (\ref{functors_J_P}). We underline that $J^*_P, J^!_P$ do not depend on a choice of a section of $P\to P/U$. Note also that $\DD\comp J^*_P\,\iso\, J^!_P\comp\DD$ naturally.

\smallskip\noindent
4.7.2 Due to its importance, recall the definition of the geometric restriction functor $\gRes^G_L:\Sph_G\to \Sph_L$ from (\cite{BG}, Proposition~4.3.3). The diagram $L\gets P\to G$ yields by functoriality the diagram 
$$
\Gr_L\getsup{\gt_P}\Gr_P\toup{\gt_G}\Gr_G
$$
The connected components of $\Gr_G$ are indexed by $\pi_1(G)$. For $\theta\in\pi_1(G)$ the component 
$\Gr_G^{\theta}$ is the one containing $\Gr_G^{\lambda}$ for any $\lambda\in\Lambda^+_G$ whose image in $\pi_1(G)$ equals $\theta$. 

For $\theta\in\pi_1(L)$ let $\Gr_P^{\theta}$ be the preimage of $\Gr_L^{\theta}$ under $\gt_P$. The following strengthened version of (\cite{BG}, Proposition~4.3.3) is derived from (\cite{MV}, Theorem~3.5) (cf. also \cite{BD}, Sections~5.3.27-5.3.30).

\begin{Pp} 
\label{Pp_tensor_functor_in_Satake}
For any $\cS\in\Sph_G$ and $\theta\in\pi_1(L)$ the complex 
$$
(\gt_P)_!(\cS\mid_{\Gr^{\theta}_P})[\<\theta, 2(\check{\rho}-\check{\rho}_L\>]
$$
lies in $\Sph_L$. The functor $\gRes^G_L: \Sph_G\to\Sph_L$ given by
$$
\cS\mapsto \mathop{\oplus}_{\theta\in\pi_1(L)} (\gt_P)_!(\cS\mid_{\Gr^{\theta}_P})[\<\theta, 2(\check{\rho}-\check{\rho}_L\>]
$$
has a natural structure of a tensor functor. The following diagram is commutative
$$
\begin{array}{ccc}
\Sph_G & \toup{\gRes^G_L} & \Sph_L\\
\uparrow\lefteqn{\scriptstyle \Loc} && \uparrow\lefteqn{\scriptstyle \Loc}\\
\Rep(\check{G}) & \toup{\Res^G_L} & \Rep(\check{L})
\end{array}
$$
\QED
\end{Pp}

 For the purposes of Lemma~\ref{Lm_Jacquet_functors} below we renormalize $\gRes^G_L$ as follows. 
We let $\ugRes^G_L:\Sph_G\to \D\Sph_L$ be given by $\ugRes^G_L(\cT)=(\gt_P)_!\gt_G^*\cT$. 

\begin{Cor} 
\label{Cor_gRes_renormalized}
The diagram is commutative
$$
\begin{array}{ccc}
\Sph_G & \toup{\ugRes^G_L} & \D\Sph_L\\
\uparrow\lefteqn{\scriptstyle \Loc} && \uparrow\lefteqn{\scriptstyle \Loc^{\gr}}\\
\Rep(\check{G}) & \toup{\Res^{\kappa_0}} & \Rep(\check{L}\times\Gm),
\end{array}
$$
where $\kappa_0: \check{L}\times\Gm\to \check{G}$ is the map whose first component is a Levi factor $\check{L}\toup{i_L}\check{G}$, and the second is 
$$
\Gm\toup{2(\check{\rho}-\check{\rho}_L)} Z(\check{L})\hook{} \check{L}\hook{i_L}\check{G}
$$
Here $Z(\check{L})$ is the center of $\check{L}$. \QED
\end{Cor}

 Write $\check{\mu}=\det M_0$ and $\check{\nu}=\det V_0$, view them as cocharacters of the center $Z(\check{L})$ of $\check{L}$. Let $\kappa:\check{L}\times\Gm\to\check{G}$ be the homomorphism, whose first component is $i_L:\check{L}\to\check{G}$, and the second component is 
$2(\check{\rho}-\check{\rho}_L)+\check{\mu}-\check{\nu}$. Let $\gRes^{\kappa}: \Sph_G\to \D\Sph_L$ denote the corresponding geometric restriction functor. 

\begin{Lm} 
\label{Lm_Jacquet_functors}
For $\cT\in \Sph_G, K\in \D_{G(\cO)}(M(F))$ there is a filtration on $J^*_P\H^{\la}_G(\cT, K)$ in the $\DG$-category $\D_{L(\cO)}(V(F))$ such that the corresponding graded complex identifies with
$$
\H^{\la}_L(\gRes^{\kappa}(\cT), J^*_P(K))
$$
For $P=G$ and a $G$-subrepresentation $V_0\subset M_0$ we have canonically
$$
J^*_P\H^{\la}_G(\cT, K)\,\iso\, \H^{\la}_L(\gRes^{\kappa}(\cT), J^*_P(K))
$$
\end{Lm}
\begin{Prf}
For $s_1,s_2\ge 0$ let 
$$
_{s_1,s_2}P(F)=\{p\in P(F)\mid t^{s_1}M\subset pM\subset t^{-s_2}M\},
$$
it is stable by left and right multiplication by $P(\cO)$, and $_{s_1,s_2}\Gr_P:=(_{s_1,s_2}P(F))/P(\cO)$ is closed in $\Gr_P$. We have a natural map $_{s_1,s_2}\Gr_P\to {_{s_1,s_2}\Gr_G}$, and at the level of reduced schemes the connected components of $_{s_1,s_2}\Gr_P$ form a stratification of $_{s_1,s_2}\Gr_G$. Set
$$
_{s_1,s_2}\Gr_L=\{x\in L(\cO)\backslash L(F)\mid t^{s_1}V\subset x V\subset t^{-s_2} V\}
$$
The map $\gt_P: \Gr_P\to\Gr_L$ yields a map still denoted $\gt_P: {_{s_1,s_2}\Gr_P}\to {_{s_1,s_2}\Gr_L}$.  
 
 Let $N,r\ge 0$, assume that $\cT$ is the extension by zero from $_{s_1,s_2}\Gr_G$, and $K\in \D_{G(\cO)}(_{N+s_1, r-s_1}M)$. For the diagram
$$
_{N,r}M\;\getsup{\pr} \; {_{N,r}M}\times{_{s_1,s_2}\Gr_G}\toup{\act_q} G(\cO/t^{N+r})\backslash {_{N+s_1, r-s_1}M}
$$
we calculate the direct image
$$
\pr_!(\act_q^*K\otimes \pr_2^*\cT)[\dim]
$$
with respect to the stratification of $_{s_1,s_2}\Gr_G$ by the connected components of $_{s_1,s_2}\Gr_P$. We have the diagram
$$
\begin{array}{ccccc}
&& _{N,r}V\times {_{s_1,s_2}P(F)} & \toup{\act} & _{N+s_1,r-s_1}V\\
&& \downarrow\lefteqn{\scriptstyle q_P} && \downarrow\lefteqn{\scriptstyle q_U}\\
_{N,r}V & \getsup{\pr}& _{N,r}V\times {_{s_1,s_2}\Gr_P} & \toup{\act_{q,P}} & P(\cO/t^{N+r})\backslash{_{N+s_1, r-s_1} V}\\
\downarrow\lefteqn{\scriptstyle i_{N,r}} && \downarrow\lefteqn{\scriptstyle i_{N,r}\times\id} && \downarrow\lefteqn{\scriptstyle i_{N+s_1,r-s_1}}\\
_{N,r}M & \getsup{\pr} & _{N,r}M\times {_{s_1,s_2}\Gr_P} & \toup{\act_{q,P}} & P(\cO/t^{N+r})\backslash {_{N+s_1, r-s_1} M}\\
&& \downarrow && \downarrow\lefteqn{\scriptstyle p}\\
&& _{N,r}M\times {_{s_1,s_2}\Gr_G}&\toup{\act_q} & G(\cO/t^{N+r})\backslash {_{N+s_1,r-s_1}M},
\end{array}
$$ 
where $\act$ sends $(m,p)$ to $p^{-1}m$, the map
$q_P$ sends $(m,p)$ to $(m, pP(\cO))$, and $q_U$ is the stack quotient under the action of $P(\cO/t^{N+r})$. Moreover, $\act_{q,P}$ fits into the diagram
$$
\begin{array}{ccc}
_{N,r}V\times {_{s_1,s_2}\Gr_P} & \toup{\act_{q,P}} & P(\cO/t^{N+r})\backslash{_{N+s_1, r-s_1} V}\\
\downarrow\lefteqn{\scriptstyle \id\times\gt_P} && \downarrow\lefteqn{\scriptstyle q}\\
_{N,r}V\times {_{s_1,s_2}\Gr_L} & \toup{\act_{q,L}} &  
L(\cO/t^{N+r})\backslash{_{N+s_1, r-s_1} V},
\end{array}
$$
Our assertion follows (the shifts can be checked using Remark~\ref{Rem_true_measures}). 
\end{Prf} 
 
\medskip

Let $\delta_U: \Gm\times M_0\to M_0$ be an action, whose fixed points set is $V_0$. Assume that $\delta_U$ contracts $M_0$ onto $V_0$, let $\pi_U: M_0\to V_0$ be the contraction map. We will apply Lemma~\ref{Lm_Jacquet_functors} under the following form.
  
\begin{Cor}  
\label{Cor_actually_used_Jacquet}
Let $K\in \P_{G(\cO)}(M(F))$ be $\Gm$-equivariant for $\delta_U$-action on $M(F)$. Assume that $\H^{\la}_G(\cT, K)$ is also $\Gm$-equivariant for $\delta_U$-action on $M(F)$. Assume that $K$ admits a $k'$-structure for some finite subfield $k'\subset k$ and, as such, is pure of weight zero. Then $J_P^*(K)$ is also pure of weight zero, and there is a (non-canonical) isomorphism
\begin{equation}
\label{iso_for_Cor3_important}
J^*_P\H^{\la}_G(\cT, K)\,\iso\, \H^{\la}_L(\gRes^{\kappa}(\cT), J^*_P(K))
\end{equation}
in $\D_{L(\cO)}(V(F))$.
\end{Cor}
\begin{Prf}
Under our assumptions, $J_P^*$ is the hyperbolic localization functor with respect to the $\delta_U$-action on $M(F)$, the assertion follows from (\cite{B}, Theorem~2) and Lemma~\ref{Lm_Jacquet_functors}.
\end{Prf}  

\medskip\noindent
\begin{Rem}
Suppose we are in the situation of Corollary~\ref{Cor_actually_used_Jacquet}. Assume in addition given an opposite parabolic $P^-$ so that $L=P\cap P^-$ is their common Levi subgroup. Assume $\pi_U$ is $k$-linear and  $P^-$-equivariant, where $P^-$ acts on $V_0$ via its quotient $P^-\to L$. 
As in (\cite{MV}, Theorem~3.6), we may also get a filtration on $J^*_P\H^{\la}_G(\cT, K)$ coming from the stratification of $\Gr_G$ by the connected components of $\Gr_{P^-}$. It will coincide with the filtration constructed in Lemma~\ref{Lm_Jacquet_functors}. 
\end{Rem}

\medskip\noindent
\begin{Rem} In our applications $\delta_U$ will be of the form $\delta_U(x)=\nu(x) x^{-r}$, $x\in\Gm$, where $\nu: \Gm\to L$ is a cocharacter of the center of $L$ acting on $V_0$ by $x\mapsto x^r$ for some $r\in\ZZ$. If $K\in \P_{G(\cO)}(M(F))$ is $\Gm$-equivariant under homotheties on $M(F)$ then both $K$ and $\H^{\la}_G(\cT, K)$ are $\Gm$-equivariant for $\delta_U$-action on $M(F)$.
\end{Rem}

\medskip\noindent
4.7.3 We similarly have for $q_V: V(F)\to\Spec k$ the functor $\bar q_{V,*}: \D_{L(\cO)}(V(F))\to \D_{L(\cO)}(\Spec k)$. 
\begin{Lm} Under the assumptions of Corollary~\ref{Cor_actually_used_Jacquet} one has a canonical isomorphism 
$$
\bar q_{V,*}(J_P^*(K))\,\iso\, \bar q_{M,*}(K)
$$
in $\D(\Spec k)$ functorial in $K\in \D_{M(\cO)}(M(F))$.
\end{Lm}
\begin{proof}
Recall the map $i_{N,r}: {_{N,r}V}\hook{} {_{N,r}M}$ from Section~4.7.1. By our assumptions, $i_{N,r}^*K\,\iso\, (\pi_{N,r})_*K$ in $\D(_{N,r}V)$, where $\pi_{N,r}: {_{N,r}M}\to {_{N,r}V}$ is the map contracting $_{N,r}M$ on $_{N,r}V$. Assume $K$ comes from an object of $\D_{G(\cO)}(_{N,r}M)$. Then 
$$
\bar q_{M, *}K\,\iso\,\RG(_{N,r}M, K)[-r\dim M_0]\,\iso\, \RG(_{N,r}V, (\pi_{N,r})_*K)[-ra-r\dim V_0], 
$$
where $a=\dim M_0-\dim V_0$. 
\end{proof}

\begin{Rem} 
\label{Rem_expected_splitting_of_filtration}
As in Definition~\ref{Def_Hecke_sheaf}, we convent to denote for $\kappa: \check{L}\times\Gm\to \check{G}$ also by $\kappa: \check{L}\times\Gm\to \check{G}\times\Gm$ the map $(\kappa, \pr_2)$. We expect that under the assumptions of Corollary~\ref{Cor_actually_used_Jacquet} the filtration on $J_P^*\H^{\la}_G(\cT, K)$ obtained in Lemma~\ref{Lm_Jacquet_functors} admits a splitting functorial in $\cT\in\Sph_G, K\in \D_{M(\cO)}(M(F))$) with the following properties: 
\begin{itemize}
\item[i)] the isomorphism of Corollary~\ref{Cor_actually_used_Jacquet} becomes uniquely defined; 
\item[ii)] it fits into the commutative diagram
$$
\begin{array}{ccc}
\bar q_{V,*} J_P^*\H^{\la}_G(\cT, K) & \iso\;\; \bar q_{M,*}\H^{\la}_G(\cT, K)\;\;\iso & \bar q_{M,*}(K)\otimes  \RG(\Gr_G, \Res^{\kappa_{M*}}\cT)\\
\downarrow\lefteqn{\scriptstyle (\ref{iso_for_Cor3_important})} && \downarrow\lefteqn{\scriptstyle \id\otimes u}\\
\bar q_{V,*}\H^{\la}_L(\Res^{\kappa}(\cT), J_P^*(K))
 &\iso & \bar q_{M,*}(K)\otimes  \RG(\Gr_L, \Res^{\kappa_{V*} }\Res^{\kappa}(\cT)), 
\end{array}
$$
where 
$$
u: \RG(\Gr_G, \Res^{\kappa_{M*}}\cT)\,\iso\, \RG(\Gr_L, \Res^{\kappa_{V*}}\Res^{\kappa}(\cT))
$$ 
is the natural isomorphism (cf. Proposition~\ref{Pp_tensor_functor_in_Satake}), and the other arrows are canonical isomorphisms given by the above lemmas. We used here the commutative diagram
$$
\begin{array}{ccc}
\Gm & \toup{2\check{\rho}\;\;}\check{G}\times\Gm \;\;\toup{\kappa_{M*}} & \check{G}\\
\downarrow\lefteqn{\scriptstyle 2\check{\rho}_L} && \uparrow\lefteqn{\scriptstyle \kappa}\\
\check{L}\times\Gm & \toup{\kappa_{V*}} & \check{L}\times\Gm,
\end{array}
$$ 
where $\kappa=(i_L, 2(\check{\rho}-\check{\rho}_L)+\check{\mu}-\check{\nu})$ with $\check{\mu}=\det M_0, \check{\nu}=\det V_0$. 
\end{itemize}
\end{Rem}

\begin{Lm} 
\label{Lm_splitting_unique}
Keep the assumptions of Corollary~\ref{Cor_actually_used_Jacquet}. Let $K=I_0$ be the constant perverse sheaf on $M(\cO)$. Assume that for any $\cT\in \Sph_G$ and any irreducible perverse sheaf $\cK$ appearing as a shifted direct summand in $\H^{\la}_L(\gRes^{\kappa}(\cT), I_0)$ one has $\bar q_{V,*}(\cK)\ne 0$. Then there is a canonical splitting of the filtration on $J_P^*\H^{\la}_G(\cT, I_0)$ obtained in Lemma~\ref{Lm_Jacquet_functors}. This splitting is functorial in $\cT\in\Sph_G$ and satisfying the properties of Remark~\ref{Rem_expected_splitting_of_filtration}.
\end{Lm}
\begin{proof}
We have $\bar q_{M,*}I_0\,\iso\,\Qlb$, $J_P^*(I_0)\,\iso\, I_0$, and $\bar q_{V,*}I_0\,\iso\,\Qlb$. By Corollary~\ref{Cor_actually_used_Jacquet}, there is a splitting of our filtration, maybe non unique. It is uniquely normalized by the property that the diagram in Remark~\ref{Rem_expected_splitting_of_filtration} ii) commutes. Indeed, two splitting differ by an endomorphism, which will act nontrivially on $\bar q_{V,*} J_P^*\H^{\la}_G(\cT, K)$.
\end{proof}

\begin{Rem} i) The assumptions of Lemma~\ref{Lm_splitting_unique} will be satisfied in all the examples we consider in this paper (for dual pairs $(\GL_n, \GL_m)$ and $(\Sp_{2n},\SO_{2m})$. \\
ii) If $L$ is a maximal torus of $G$ then the isomorphism $u$ from Remark~\ref{Rem_expected_splitting_of_filtration} is explained in details in (\cite{BR}, Remark~5.10). 
\end{Rem}

\bigskip\noindent
4.8 {\scshape Fourier transform}

\medskip\noindent
Recall the notation $\Omega$ from Section~4.1. Let us define the Fourier transform functors $\Four_{\psi}: \D(M(F))\to \D(M^*\otimes\Omega(F))$ and 
\begin{equation}
\label{Fourier_transform_def_M}
\Four_{\psi}: \D_{G(\cO)}(M(F))\to \D_{G(\cO)}(M^*\otimes\Omega(F))
\end{equation}

 We actually will use the following a bit more general functor. Given a decomposition $M_0\,\iso\, M_1\oplus M_2$ into direct sum of vector spaces, one defines the Fourier transform 
\begin{equation}
\label{Fourier_transform_for_def}
\Four_{\psi}: \D(M(F))\to \D(M_1^*\otimes\Omega(F)\oplus M_2(F))
\end{equation}
as follows. For $N\ge 0$ we have a natural evaluation map $\ev: {_{N,N}M_1}\times {_{N,N}(M_1^*\otimes\Omega)}\to \A^1$ sending $(m, m^*)$ to $\Res \<m, m_1\>$. It gives rise to the usual Fourier transform functor
$$
\Four_{\psi}: \D(_{N,N}M)\,\iso\,\D(_{N,N}M_1^*\otimes\Omega\oplus {_{N,N} M_2})
$$
For $N'\ge N$ these functors are compatible with the transition functors $\D(_{N,N}M)\to \D(_{N',N'}M)$ in (\ref{diag_categories_limit}), so give rise to the desired functor (\ref{Fourier_transform_for_def}). From the usual properties of the Fourier transform we learn that (\ref{Fourier_transform_for_def}) is an equivalence of triangulated categories, which preserves the perversity. 

 Assume in addition that $M_0\,\iso\, M_1\oplus M_2$ is a decomposition of $M_0$ into a direct sum of $G$-modules. Then similarly the usual Fourier transform functors 
$$
\Four_{\psi}: \D_{G(\cO)}(_{N,N}M)\,\iso\,\D_{G(\cO)}(_{N,N}M_1^*\otimes\Omega\oplus {_{N,N} M_2}),
$$
being compatible with the transition functors in (\ref{diag_categories_limit}), give rise to the functor
\begin{equation}
\label{Fourier_transform_for_def_equiv}
\Four_{\psi}: \D_{G(\cO)}(M(F))\,\iso\,\D_{G(\cO)}(M_1^*\otimes\Omega(F)\oplus M_2(F)),
\end{equation}
which satisfies the same formal properties. 

\begin{Rem} 
\label{Rem_comp_partial_Fourier}
The following diagram commutes
$$
\begin{array}{ccc}
\D_{G(\cO)}(M(F)) & \toup{\Four_{\psi}} & \D_{G(\cO)}(M_1^*\otimes\Omega(F)\oplus M_2(F))\\
 & \searrow\lefteqn{\scriptstyle \Four_{\psi}} & \downarrow\lefteqn{\scriptstyle \Four_{\psi}}\\
 && \D_{G(\cO)}(M^*\otimes\Omega(F)),
\end{array}
$$
that is, the composition of two partial Fourier transforms identifies with the complete Fourier transform. 
\end{Rem}
  
\begin{Lm} 
\label{Lm_Fourier_Hecke_commute}
The functor (\ref{Fourier_transform_def_M}) commutes with Hecke operators. Namely, 
there is an isomorphism functorial in $\cT\in\Sph_G$ and $K\in\D_{G(\cO)}(M(F))$
$$
\Four_{\psi}\H^{\la}_G(\cT, K)\,\iso\, \H^{\la}_G(\cT, \Four_{\psi}(K))
$$
\end{Lm}
\begin{Prf}
\Step 1 Pick $s_1,s_2\ge 0$ so that $\cT$ is the extension by zero from $_{s_1,s_2}\Gr_G$. Pick $r,r_1,N, N_1$ large enough compared to $s_i$ and $K$. In particular, we assume
\begin{equation}
\label{equation_inequalities_assumptions_Lm11}
r-N_1\ge s_1+s_2\;\;\;\;\mbox{and}\;\;\;\; r_1-N\ge s_1+s_2
\end{equation}
Let $s\ge\max\{s_1+s_2, N+r, N_1+r_1\}$. Consider the diagram
$$
\begin{array}{ccccc}
&& G(\cO/t^s)\backslash {_{N+s_1, r-s_1}M}\\
&& \uparrow\lefteqn{\scriptstyle {\act_{q,s}}}\\
G(\cO/t^s)\backslash {_{N,r}M} & \getsup{\pr} & G(\cO/t^s)\backslash ({_{N,r}M\times {_{s_1,s_2}\Gr_G}}) \\
\uparrow\lefteqn{\scriptstyle\alpha} && \uparrow\\
G(\cO/t^s)\backslash ({_{N_1,r_1}(M^*\otimes\Omega)}\times_{N,r}M) 
& \getsup{\pr} & G(\cO/t^s)\backslash (_{N_1,r_1}(M^*\otimes\Omega)\times_{N,r}M
\times {_{s_1,s_2}\Gr_G})\\
\downarrow\lefteqn{\scriptstyle\beta}\\
G(\cO/t^s)\backslash (_{N_1,r_1}(M^*\otimes\Omega)),
\end{array}
$$
where the square is cartesian, all the quotients are taken in the stack sense, the action of $G(\cO/t^s)$ on all the involved schemes is diagonal. We have denoted by $\alpha$ and $\beta$ are the projections.

 By our assumptions,
$$
\H^{\la}_G(\cT,K)\,\iso\, \pr_!(\cT\otimes \act_{q,s}^*K)[\dim]
$$
for a suitable shift. Assuming $K\in\D_{G(\cO)}(_{\tilde N,\tilde r}M)$ with $N_i,r_i$ sufficiently large with respect to $\tilde N,\tilde r$, we get 
$$
\Four_{\psi}(\H^{\la}_G(\cT,K))\,\iso\, 
\beta_!(\ev^*\cL_{\psi}\otimes\alpha^* \H^{\la}_G(\cT,K))[\dimrel(\alpha)]
$$
Here $\ev: G(\cO/t^s)\backslash ({_{N_1,r_1}(M^*\otimes\Omega)}\times_{N,r}M) \to\A^1$ is the evaluation map, it is correctly defined because $r-N_1$ and $r_1-N$ are nonnegative. 

 Consider the diagram
$$
\begin{array}{ccc}
&& G(\cO/t^s)\backslash {_{N+s_1, r-s_1}M}\\
&& \uparrow\lefteqn{\scriptstyle\alpha'}\\
G(\cO/t^s)\backslash (_{N_1,r_1}(M^*\otimes\Omega)\times_{N,r}M
\times {_{s_1,s_2}\Gr_G})
 & \toup{\act_{q,s}} & G(\cO/t^s)\backslash (_{N_1+s_2,r_1-s_2}(M^*\otimes\Omega)\times_{N+s_1,r-s_1}M)\\
\downarrow && \downarrow\lefteqn{\scriptstyle\beta'}\\
G(\cO/t^s)\backslash (_{N_1,r_1}(M^*\otimes\Omega)
\times {_{s_1,s_2}\Gr_G}) &\toup{\act'_{q,s}}& G(\cO/t^s)\backslash {_{N_1+s_2,r_1-s_2}(M^*\otimes\Omega)}\\
\downarrow\lefteqn{\scriptstyle\pr'}\\
G(\cO/t^s)\backslash {_{N_1,r_1}(M^*\otimes\Omega)}
\end{array}
$$
where $\alpha',\beta',\pr'$ are the projections. The square in the above diagram is not cartesian, write 
$$
b:\cY\to G(\cO/t^s)\backslash (_{N_1+s_2,r_1-s_2}(M^*\otimes\Omega)\times_{N+s_1,r-s_1}M)
$$
for the map obtained from $\act'_{q,s}$ by the base change $\beta'$. Then  
\begin{equation}
\label{closed_immersion_for_Lm11}
G(\cO/t^s)\backslash (_{N_1,r_1}(M^*\otimes\Omega)\times_{N,r}M
\times {_{s_1,s_2}\Gr_G})\hook{}\cY
\end{equation}
is naturally a closed substack. Let 
$$
\ev: G(\cO/t^s)\backslash (_{N_1+s_2,r_1-s_2}(M^*\otimes\Omega)\times_{N+s_1,r-s_1}M)\to\A^1
$$
be the evaluation map, it is correctly defined due to (\ref{equation_inequalities_assumptions_Lm11}). 
By our assumptions,
$$
\Four_{\psi}(K)\,\iso\, \beta'_!(\ev^*\cL_{\psi}\otimes \alpha'^*K)[\dimrel(\alpha')]
$$ 
and 
$$
\H^{\la}_G(\cT, \Four_{\psi}(K))\,\iso\, \pr'_!(\cT\otimes (\act'_{q,s})^*\Four_{\psi}(K))[\dim]
$$ 
Since $N,r$ are large enough compared to $\tilde N,\tilde r$, it follows that $b^*(\alpha')^*K$ is the extension by zero under (\ref{closed_immersion_for_Lm11}). The desired result follows now from the base change theorem.
\end{Prf}

\medskip

 Note that for the functor (\ref{Fourier_transform_def_M}) we have $\Four_{\psi}(I_0)\,\iso\, I_0$ canonically. 
 
\medskip\noindent
4.9 {\scshape Extensions of actions}

\medskip\noindent
Let $G$ be a connected reductive group, $P$ and 
$P^-$ two opposite parabolic subgroups in $G$ with common Levi subgroup $L=P\cap P^-$. 

\begin{Lm} 
\label{Lm_extending_actions_f_dim}
Let $Y$ be a scheme (of finite type over $k$) with a $G$-action. Then we have a diagram of equivariant categories of perverse sheaves
$$
\begin{array}{ccc}
P_P(Y) & \subset & P_L(Y)\\
\cup && \cup\\
P_G(Y) & \subset & P_{P^-}(Y),
\end{array}
$$ 
where all the functors are fully faithful embeddings. Moreover, $P_P(Y)\cap P_{P^-}(Y)=P_G(Y)$, that is, if an object $K\in P_L(Y)$ lies in both $P_P(Y)$ and $P_{P^-}(Y)$ then $K\in P_G(Y)$.
\end{Lm}
\begin{Prf}
The natural maps between the stack quotients $Y/L\to Y/P\to Y/G$ are smooth of fixed relative dimension, surjective, and have connected fibres. By (\cite{G}, Lemma~4.8), they induce the corresponding fully faithful embeddings of categories. 

 Now assume $K$ is an object of $P_P(Y)\cap P_{P^-}(Y)$. Let $W$ be the image of the product map $m: P\times P^-\to G$. We have a diagram
$$
\begin{array}{ccc}
P\times P^-\times Y & \toup{\act'} & Y\\
\downarrow\lefteqn{\scriptstyle m\times\id} && \downarrow\lefteqn{\scriptstyle \id}\\
W\times Y & \toup{\act_W} & Y,
\end{array}
$$
where $\act'$ sends $(p_1,p_1,y)$ to $p_1p_2y$. The map $m: P\times P^-\to W$ is smooth and surjective with connected fibres. So, by \select{loc.cit.}, the equivariance isomorphism $(\act')^*K\,\iso\,\Qlb\boxtimes K$ descends to an isomorphism $\act_W^*K\,\iso\, \Qlb\boxtimes K$ over $W\times Y$. 

 Further, the product map $m_W: W\times W\to G$ is smooth and surjective with connected fibres. Indeed, any fibre of $m_W$ is isomorphic to $zPP^-\cap P^-P$ for some $z\in G$. The latter intersection is connected, because it is open in $G$. So, for the action map $\act_{W\times W}: W\times W\times Y\to Y$ the equivariance isomorphism $\act_{W\times W}^*K\,\iso\,\Qlb\boxtimes K$ descends to the desired isomorphism $\act^*K\,\iso\, \Qlb\boxtimes K$ over $G\times Y$.
\end{Prf}
 
\medskip

 Now assume $M_0$ is a finite-dimensional representation of $G$, set $M=M_0\otimes_k \cO$. Let $U$ be the unipotent radical of $P$. The following result will be used in Section~6.2.
 
\begin{Lm} 
\label{Lm_extending_actions_on_Schwarz}
i) We have a diagram of fully faithful embeddings of categories
$$
\begin{array}{ccc}
P_{P(\cO)}(M(F)) & \subset & P_{L(\cO)}(M(F))\\
\cup && \cup\\
P_{G(\cO)}(M(F)) & \subset & P_{P^-(\cO)}(M(F)),
\end{array}
$$
The  intersection $P_{P(\cO)}(M(F))\cap P_{P^-(\cO)}(M(F))$ inside $P_{L(\cO)}(M(F)$ equals $P_{G(\cO)}(M(F))$. 

\smallskip\noindent
ii)  We have fully faithful embeddings $P_{L(\cO)}(M(F))\subset P(M(F)) \supset P_{U(\cO)}(M(F))$. The intersection $P_{L(\cO)}(M(F))\cap P_{U(\cO)}(M(F))$ equals $P_{P(\cO)}(M(F))$. 
\end{Lm}
\begin{Prf} 
i) Given $N,r\ge 0$, by (\cite{G}, Lemma~4.8), we get a diagram of fully faithful embeddings
$$
\begin{array}{ccc}
P_{P(\cO/t^{s})}(_{N,r}M) & \subset & P_{L(\cO/t^{s})}(_{N,r}M)\\
\cup && \cup\\
P_{G(\cO/t^{s})}(_{N,r}M) & \subset & P_{P^-(\cO/t^{s})}(_{N,r}M)
\end{array}
$$
with $s=N+r$. 
Let $K$ be an object of $P_{P(\cO/t^s)}(_{N,r}M)\cap P_{P^-(\cO/t^s)}(_{N,r}M)$. 

 Let $m: P\times P^-\to W$ and $m:W\times W\to G$ be as in Lemma~\ref{Lm_extending_actions_f_dim}. The  induced maps $m: P(\cO/t^s)\times P^-(\cO/t^s)\to W(\cO/t^s)$ and $m: W(\cO/t^s)\times W(\cO/t^s)\to G(\cO/t^s)$ are again smooth and surjective. Indeed, if $Y_1\to Y_2$ is a smooth surjective morphism of affine algebraic varieties, $A$ is an Artin $k$-algebra then $Y_1(A)$ is a scheme, and the induced map $Y_1(A)\to Y_2(A)$ is smooth and surjective. As in Lemma~\ref{Lm_extending_actions_f_dim}, one shows now that $K\in P_{G(\cO/t^s)}(_{N,r}M)$. The first assertion follows.

\smallskip\noindent
ii) Given $N,r\ge 0$ as above for $s=N+r$ one gets a diagram of fully faithful embeddings
$$
P_{L(\cO/t^s)}(_{N,r}M)\subset P(_{N,r}M)\supset P_{U(\cO/t^s)}(_{N,r}M)
$$
If $K\in P_{L(\cO/t^s)}(_{N,r}M)\cap P_{U(\cO/t^s)}(_{N,r}M)$ then the equivariance isomorphisms for $L$ and $U$ yield an isomorphism $a^*K\,\iso\, \Qlb\boxtimes K$, where $a: L(\cO/t^s)\times U(\cO/t^s)\times {_{N,r}M}\to {_{N,r}M}$ is the map sending $(g,u,m)$ to $gum$. The product induces an isomorphism $L(\cO/t^s)\times U(\cO/t^s)\,\iso\, P(\cO/t^s)$, so $K$ lies in the full subcategory $\P_{P(\cO/t^s)}(_{N,r}M)\subset \P(_{N,r}M)$. Our assertion follows.
\end{Prf}

\bigskip\medskip

\centerline{\scshape 5. Geometric model of the Weil representation of $\GL_m\times\GL_n$}

\bigskip\noindent
5.1 Let $U_0=k^m, L_0=k^n$ be the standard $k$-vector spaces of dimensions $m$ and $n$. For Section~5 we let $G=\GL(L_0)$ and $H=\GL(U_0)$. Let $\Pi_0=U_0\otimes L_0$.

 Set $U=U_0(\cO)$, $L=L_0(\cO)$ and $\Pi=\Pi_0(\cO)$. 
Let $T_G\subset B_G\subset G$ be the torus of diagonal matrices and the Borel subgroup of upper-triangular matrices. We identify $\Lambda_G\,\iso\,\ZZ^n$ in the usual way. Write $\check{\omega}_i\in\check{\Lambda}^+_G$ be the h.w. of the representation $\wedge^i L_0$ of $G$. The objects $T_H\subset B_H\subset H$ are defined similarly for $H$. By some abuse of notation, $\check{\omega}_i\in\check{\Lambda}^+_H$ will also denote the h.w. of the $H$-representation $\wedge^i U_0$. 
Keep the notations of Section~4, in particular we write $_{N,r}\Pi=t^{-N}\Pi/t^r\Pi$, and $I_0$ is the constant sheaf on $_{0,0}\Pi$. 

 We are going to describe the submodule over $\Sph_G$ (resp., over $\Sph_H$) in $\D_{(G\times H)(\cO)}(\Pi(F))$ generated by $I_0$. 
Assume $m\ge n$. 

 Let $U_1\oplus U_2\,\iso\,U_0$ be the direct sum decomposition, where $U_1$ (resp., $U_2$) is 
generated by the first $n$ (resp., last $m-n$) base vectors.
Let $P\subset H$ be the parabolic subgroup preserving $U_1$, $U_H\subset P$ be its unipotent radical. Let $M=\GL(U_1)\times\GL(U_2)\subset P$ be the standard Levi factor. Let $\kappa: \check{G}\times\Gm\to\check{H}$ be the composition 
$$
\check{G}\times\Gm\toup{\id\times 2\check{\rho}_{\GL(U_2)}} \check{G}\times\check{\GL}(U_2)=\check{M}\hook{} \check{H}
$$ 


Write $\gRes^{\kappa}:\Sph_H\to\D\Sph_G$ for the functor corresponding (in view of $\Loc$ and $\Loc^{\gr}$) to the restriction $\Rep(\check{H})\to\Rep(\check{G}\times\Gm)$ with respect to $\kappa$. Here is the main result of Section~5. 

\begin{Pp}
\label{Pp_action_Hecke_local_GL_m_GL_n} 
The two functors $\Sph_H\to \D_{(G\times H)(\cO)}(\Pi(F))$
given by
\begin{equation}
\label{functors_GL_n_GL_m_Section5}
\cT\mapsto \H^{\la}_H(\cT, I_0)\;\;\;\;\mbox{and}\;\;\;\;
\cT\mapsto \H^{\la}_G(\gRes^{\kappa}(\cT), I_0)
\end{equation}
are canonically isomorphic (once a decomposition $U_1\oplus U_2\,\iso\,U_0$ is fixed). This isomorphism is compatible with the tensor structures on $\Sph_H, \Sph_G$. 
\end{Pp}
 
  Let $N,r\in\ZZ$ with $N+r\ge 0$. Think of $v\in\Pi(F)$ as a map $v: U^*(F)\to L(F)$. For $v\in {_{N,r}\Pi}$ 
let $U_{v,r}=v(U^*)+t^rL$, this is a $\cO$-lattice in $L(F)$. 
For $\lambda\in\Lambda^+_G$ satisfying 
\begin{equation} 
\label{eq_lambda_stratification_NrPi} 
 \<-w_0^G(\lambda), \check{\omega}_1\>\le N\;\;\;\mbox{and}\;\;\; \<\lambda, \check{\omega}_1\>\le r
\end{equation} 
let $_{\lambda,r}\Pi^0\subset {_{N,r}\Pi}$ be the locally closed subscheme  of those $v\in {_{N,r}\Pi}$ for which $t^{a_n}L/(v(U^*)+t^{a_1}L)$ is isomorphic to $\cO/t^{a_1-a_n}\oplus\ldots\oplus \cO/t^{a_n-a_n}$ as $\cO$-module. Here $\lambda=(a_1\ge\ldots\ge a_n)$. In other words, for $v\in {_{N,r}\Pi}$ we have $v\in {_{\lambda,r}\Pi^0}$ iff $U_{v,r}\in \Gr^{\lambda}_G$. 

 One checks that the $G(\cO)\times H(\cO)$-orbits on $_{N,r}\Pi$ are exactly $_{\lambda,r}\Pi^0$ for $\lambda\in\Lambda^+_G$ satisfying (\ref{eq_lambda_stratification_NrPi}). The fact that the set of $(G\times H)(\cO)$-orbits on $_{N,r}\Pi$ is finite also follows from (\cite{GN}, Theorem~3.2.1), because $\Pi_0$ is a spherical $G\times H$-variety.

 Given $\lambda\in\Lambda^+_G$ let now $N=\<-w_0^G(\lambda), \check{\omega}_1\>$ and $r=\<\lambda, \check{\omega}_1\>$. 
By 4.4, $\H^{\lambda}_G(I_0)\in \D_{(G\times H)(\cO)}(_{N,r}\Pi)$. 
Define the closed subscheme $_{\lambda}\Pi\subset {_N\Pi}$ as follows. A point $v\in {_N\Pi}$ lies in $_{\lambda}\Pi$ iff
for $i=1,\ldots, n$ the map
$$
\wedge^i U^*\toup{\wedge^i v} (\wedge^i L)(-\<w_0(\lambda), \check{\omega}_i\>)
$$
is regular. The scheme $_{\lambda}\Pi$ is stable under translations by $t^r\Pi(\cO)$, so there is a unique closed subscheme $_{\lambda,r}\Pi\subset {_{N,r}\Pi}$ such that $_{\lambda}\Pi$ is the preimage of $_{\lambda,r}\Pi$ under the projection $_N\Pi\to {_{N,r}\Pi}$. Under our assumptions the map (\ref{map_pi_rem2}) factors as
$$
_{0,r}\Pi\ttimes \ov{\Gr}^{\lambda}_G\toup{\pi} {_{\lambda,r}\Pi}\hook{} {_{N,r}\Pi}
$$

\begin{Pp} 
\label{Pp_action_Hecke_GL_n_smallest}
For $\lambda\in\Lambda^+_G$ we have a canonical isomorphism $\H^{\lambda}_G(I_0)\,\iso\, \IC(_{\lambda,r}\Pi^0)$ with the intersection cohomology sheaf of $_{\lambda,r}\Pi^0$. 
\end{Pp}
\begin{Prf} 
Note that $_{\lambda,r}\Pi^0\subset{_{\lambda,r}\Pi}$ is an open subscheme. The map $_{0,r}\Pi\ttimes \ov{\Gr}^{\lambda}_G\toup{\pi} {_{\lambda,r}\Pi}$ is an isomorphism over $_{\lambda,r}\Pi^0$, in particular $\dim{_{\lambda,r}\Pi^0}=rnm+\<\lambda, 2\check{\rho}_G-m\check{\omega}_n\>$. 

 The scheme $_{\lambda,r}\Pi$ is stratified by locally closed subschemes $_{\mu,r}\Pi^0$, where $\mu\in\Lambda^+_G$ satisfies 
(\ref{eq_lambda_stratification_NrPi}) and 
\begin{equation}
\label{equation_stratification_GL_n_strange}
\<w_0^G(\lambda-\mu), \check{\omega}_i\>\le 0
\end{equation}
for $i=1,\ldots,n$. Further, $_{0,r}\Pi\ttimes \ov{\Gr}^{\lambda}_G$ is stratified by locally closed subschemes $_{0,r}\Pi\ttimes \Gr^{\mu}_G$ with $\mu\in\Lambda^+_G$, $\mu\le\lambda$. Let us show that $\pi$ is stratified small (in the sense of \cite{MV}) with respect to these stratifications. 

 Let $\mu\in\Lambda^+_G$ satisfy (\ref{eq_lambda_stratification_NrPi}) and (\ref{equation_stratification_GL_n_strange}), take $v\in {_{\mu,r}\Pi^0}$. Let $Y$ be the fibre of $\pi: {_{0,r}\Pi\ttimes \Gr^{\lambda}_G}\to {_{\lambda,r}\Pi}$ over $v$. We must show that $2\dim Y\le \<\lambda-\mu, 2\check{\rho}_G-m\check{\omega}_n\>$. 
 
 From (\ref{equation_stratification_GL_n_strange}) it follows that $\<\lambda-\mu, \check{\omega}_n\>\le 0$. 
So, to finish the proof it suffices to show that $2\dim Y\le \<\lambda-\mu, 2\check{\rho}_G-n\check{\omega}_n\>$. 

 The scheme $Y$ classifies $\cO$-lattices $L'\subset L(F)$ such that $L'\in\Gr_G^{\lambda}$ and $U_{v,r}\subset L'$. Stratify $Y$ by locally closed subschemes $Y_{\tau}$ indexed by $\tau\in\Lambda^+_G$, which are \select{very positive}. We call 
$$
\tau=(b_1\ge\ldots\ge b_n)
$$ 
very positive iff $b_n\ge 0$. By definition, the subscheme $Y_{\tau}$ classifies $L'\in Y$ such that $U_{v,r}$ is in the position $\tau$ with respect to $L'$. 
Now by (\cite{MV}, Lemma~4.4), if $Y_{\tau}$ is nonempty then
$$
\dim Y_{\tau}\le \<\lambda+\tau-\mu, \check{\rho}_G\>
$$
So, we have to show that $\<\tau, 2\check{\rho}_G\>\le \<\lambda-\mu, -n\check{\omega}_n\>$. The formula for virtual dimensions $\dim(L/L')+\dim(L'/U_{v,r})=\dim(L/U_{v,r})$ reads 
$\<\tau+\lambda-\mu, \check{\omega}_n\>=0$. Thus, we are reduced to show that 
\begin{equation}
\label{equality_miracle_GL_n}
\<\tau, n\check{\omega}_n-2\check{\rho}_G\>\ge 0
\end{equation}
This inequality follows from the fact that $\tau$ is very positive, because 
$$
n\check{\omega}_n-2\check{\rho}_G=(1,3,5,\ldots, 2n-1)
$$ 
is very positive. Moreover, the inequality (\ref{equality_miracle_GL_n}) is strict unless $\tau=0$. Since $\tau=0$ iff $\lambda=\mu$, we are done.
\end{Prf}

\begin{Cor} 
\label{Cor_Hecke_GL_n}
i) The functor $\Sph_G\to \D_{(G\times H)(\cO)}(\Pi(F))$ given by $\cT\mapsto \H^{\la}_G(\cT, I_0)$ takes values in $\P_{(G\times H)(\cO)}(\Pi(F))$. The corresponding functor 
$$
\Sph_G\to \P_{(G\times H)(\cO)}(\Pi(F))
$$ 
is fully faithful, its image is the full subcategory $\P^{ss}_{(G\times H)(\cO)}(\Pi(F))$ of semi-simple objects in $\P_{(G\times H)(\cO)}(\Pi(F))$. \\
ii) For any $\lambda\in\Lambda^+_G$ we have 
$\Ext^1_{(G\times H)(\cO)}(\IC(_{\lambda,r}\Pi^0), \IC(_{\lambda,r}\Pi^0))=0$.
\end{Cor}
\begin{Prf}
For  $\lambda\in\Lambda^+_G$ let $r=\<\lambda, \check{\omega}_1\>$ and $N=\<-w_0^G(\lambda), \check{\omega}_1\>$. Pick $s\ge N+r$. The stabilizor, say $K$, in $(G\times H)(\cO/t^s)$ of a point of $_{\lambda,r}\Pi^0$ is connected. So, the irreducible objects of $\P_{(G\times H)(\cO)}(\Pi(F))$ are exactly $\IC(_{\lambda,r}\Pi^0)$, $\lambda\in\Lambda^+_G$. Part i) follows. 

  We have a canonical equivalence $\P_{(G\times H)(\cO)}(_{\lambda,r}\Pi^0)\,\iso\, \P_K(\Spec k)$. 
By (\cite{G}, Lemma~4.8), the connectedness of $K$ implies that $\P_K(\Spec k)$ is equivalent to the category of vector spaces. If $0\to \IC(_{\lambda,r}\Pi^0)\to \cK\to \IC(_{\lambda,r}\Pi^0)\to 0$ is an exact sequence in $\P_{(G\times H)(\cO/t^s)}(_{N,r}\Pi)$ then $\cK$ is the intermediate extension from $_{\lambda,r}\Pi^0$. Part ii) follows. 
\end{Prf}

\begin{Rem}
\label{Rem_Aut_O_equivariance}
 i) As in Section~4.5, one may strengthen Corollary~\ref{Cor_Hecke_GL_n}i) saying that the functor $\cT\mapsto \H^{\la}_G(\cT, I_0)$ takes values in the category $\P_{(G\times H)(\cO)\rtimes\Aut^0\cO}(\Pi(F))$. The functors in Proposition~\ref{Pp_action_Hecke_local_GL_m_GL_n} may also be seen as taking values in the $\Aut^0(\cO)$-equivariant version of the corresponding category.\\
ii) The category $\P_{(G\times H)(\cO)}(\Pi(F))$ is not semi-simple in general. To have an example, take $n=m=2$ and $\lambda=(1,0)$. Let $Y\subset {_{0,1}\Pi}$ be the support of $\IC(_{\lambda,1}\Pi^0)$ then $\dim Y=3$ and $\dim{_{0,1}\Pi}=4$. The restriction to $Y$ yields a nontrivial map $I_0\to \IC(_{\lambda,1}\Pi^0)[1]$ in $\D_{(G\times H)(\cO)}(_{0,1}\Pi)$.
\end{Rem}

\medskip

\begin{Prf}\select{of Proposition~\ref{Pp_action_Hecke_local_GL_m_GL_n}}
 
\Step 1 Assume first $n=m$. Interchanging $U_0$ and $L_0$, one derives from Proposition~\ref{Pp_action_Hecke_GL_n_smallest}
that the functors $\Rep(\GL_n)\to \P_{(G\times H)(\cO)}(\Pi(F))$ given by 
$$
V\mapsto \H^{\la}_H(V, I_0)\;\;\;\;\mbox{and}\;\;\;\;
V\mapsto \H^{\la}_G(V, I_0)
$$
are canonically isomorphic. For $n=m$ we are done.

\Step 2 For $m\ge n$ consider the Jacquet functors 
$$
J_P^*: \D_{(G\times H)(\cO)}(\Pi(F))\to \D_{(G\times M)(\cO)}(U_1\otimes L_0(F))
$$ 
We have $J_P^*(I_0)\,\iso\, I_0$ canonically. The action of $\GL(U_2)$ on $U_1\otimes L_0$ is trivial, so
$\cS\in \Sph_{\GL(U_2)}$ acts on $I_0\in \D_{(G\times M)(\cO)}(U_1^*\otimes L_0(F))$ as 
$$
\H^{\la}_{\GL(U_2)}(\cS, I_0)\,\iso\, I_0\otimes\RG(\Gr_{\GL(U_2)}, \cS)
$$ 
As a representation of $H$, $\det(U_0\otimes L_0)$ is the character $n\check{\omega}_m\in\check{\Lambda}^+_H$. As a representation of $M$, $\det(U_1\otimes L_0)$ is the character $n\check{\omega}_n\in\check{\Lambda}^+_H$. Thus, let $\kappa_1: \check{\GL}(U_1)\times\Gm\to \check{H}$ be the composition
$$
 \check{\GL}(U_1)\times\Gm\toup{i\times (2\check{\rho}_H-2\check{\rho}_{\GL(U_1)}-n\check{\omega}_n+n\check{\omega}_m)} \check{M}\hook{}\check{H},
$$
where $i:\check{GL}(U_1)\hook{}\check{M}$ is the natural inclusion. Let $\gRes^{\kappa_1}:\Sph_H\to\D\Sph_{\GL(U_1)}$ be the corresponding restriction functor. From Corollary~\ref{Cor_actually_used_Jacquet} and Lemma~\ref{Lm_splitting_unique} 
we get for $\cT\in\Sph_H$ a canonical isomorphism
$$
J_P^*\H^{\la}_H(\cT, I_0)\,\iso\, \H^{\la}_{\GL(U_1)}(\gRes^{\kappa_1}(\cT), I_0)
$$

 Let $\kappa_2: \check{G}\times\Gm\to \check{H}$ be the map obtained from $\kappa_1$ via the canonical identification $\check{\GL}(U_1)\,\iso\,\check{G}$. By Step 1, we have a canonical isomorphism
$$
\H^{\la}_{\GL(U_1)}(\gRes^{\kappa_1}(\cT), I_0)\,\iso\,
\H^{\la}_G(\gRes^{\kappa_2}(\cT), I_0)
$$
in $\D_{(G\times M)(\cO)}(U_1\otimes L_0(F))$. 

 Further, we may think of $J_P^*$ as the Jacquet functor corresponding to the parabolic subgroup $P\times G$ of $H\times G$. As a representation of $G$, 
$\det(U_0\otimes L_0)\otimes{\det(U_1\otimes L_0)}^{-1}$
is the character $(m-n)\check{\omega}_n\in\check{\Lambda}^+_G$. So, let $\kappa_3:\check{G}\times\Gm\to\check{G}\times\Gm$ be the map, whose second component $\check{G}\times\Gm\to\Gm$ is the projection, and the first component is 
$$
\check{G}\times\Gm\toup{(\id, (m-n)\check{\omega}_n)}\check{G}
$$
Write $\gRes^{\kappa_3}:\D\Sph_G\to\D\Sph_G$ for the corresponding geometric restriction functor. By Corollary~\ref{Cor_actually_used_Jacquet} and Lemma~\ref{Lm_splitting_unique}, 
for $\cS\in\D\Sph_G$ we get a canonical isomorphism
$$
J_P^*\H^{\la}_G(\cS, I_0)\,\iso\, \H^{\la}_G(\gRes^{\kappa_3}(\cS), I_0)
$$
in $\D_{(G\times M)(\cO)}(U_1\otimes L_0(F))$.  From Corollary~\ref{Cor_Hecke_GL_n} we conclude that 
\begin{equation}
\label{equivalence_for_GL_n_Jacquet}
J_P^*: \DP^{ss}_{(G\times H)(\cO)}(\Pi(F))\to \DP^{ss}_{(\GL(U_1)\times G)(\cO)}(U_1\otimes L_0(F))
\end{equation}
is an equivalence. The equality
$$
2\check{\rho}_H-2\check{\rho}_{\GL(U_1)}-2\check{\rho}_{\GL(U_2)}+n\check{\omega}_m-m\check{\omega}_n=0
$$
shows that  the composition $\check{G}\times\Gm\toup{\kappa_3}\check{G}\times\Gm\toup{\kappa}\check{H}$ equals $\kappa_2$. 

 Summarizing, for $\cT\in\Sph_H$ we get a canonical isomorphism
$$
J_P^*\H^{\la}_H(\cT, I_0)\,\iso\, J_P^*\H^{\la}_G(\gRes^{\kappa}(\cT), I_0)
$$
in $\D_{(G\times M)(\cO)}(U_1\otimes L_0(F))$, it could only depend on the decomposition $U_0\,\iso\, U_1\oplus U_2$, which is fixed. 
Now (\ref{equivalence_for_GL_n_Jacquet}) garantees that this isomorphism can be lifted to the desired isomorphism of functors (\ref{functors_GL_n_GL_m_Section5}). 
\end{Prf}

\medskip

 We will need the following version of Proposition~\ref{Pp_action_Hecke_local_GL_m_GL_n}. Set $\Pi_1=U_0^*\otimes L_0$. Recall the functor $\ast: \Sph(\check{H})\,\iso\,\Sph(\check{H})$ from Section~2.2.1. 

\begin{Cor} 
\label{Cor_used_GL_n_GL_m}
The two functors $\Sph_H\to \D_{(G\times H)(\cO)}(\Pi_1(F))$
given by
$$
\cT\mapsto \H^{\la}_H(\ast\cT, I_0)=\H^{\ra}_H(\cT, I_0)\;\;\;\;\mbox{and}\;\;\;\;
\cT\mapsto \H^{\la}_G(\gRes^{\kappa}(\cT), I_0)
$$
are canonically isomorphic, they are compatible with the tensor structures on $\Sph_H,\Sph_G$.
\QED
\end{Cor}

\bigskip\bigskip

\centerline{\scshape 6. Geometric model of the Weil representation of $\SO_{2m}\times\Sp_{2n}$}

\bigskip\noindent
6.1 Let $U_0=k^m, L_0=k^n$. Set $V_0=U_0\oplus U_0^*$, we equip it with the symmetric form $\Sym^2 V_0\to k$ as in Section~3.2. Set $H=\SO(V_0)$. 

 Let $P_H\subset H$ be the parabolic subgroup preserving $U_0$, $U_H\subset P_H$ be its unipotent radical. Write $Q_H=\GL(U_0)\,\iso\, \GL_m$ for the standard Levi factor of $P_H$. We equip it with the maximal torus $T_H$ of diagonal matrices and the Borel subgroup of upper-triangular matrices (its preimage in $P_H$ is a Borel subgroup, which yields our choice of positive roots). 

  Set $M_0=L_0\oplus L_0^*$, we equip it with the symplectic form $\wedge^2 M_0\to k$ 
arising from pairing between $L_0$ and $L_0^*$, so $L_0$ and $L_0^*$ are lagrangian subspaces in $M_0$. Set $G=\Sp(M_0)$. Let $P_G\subset G$ be the parabolic subgroup preserving $L_0$, $U_G\subset P_G$ its unipotent radical. Write $Q_G=\GL(L_0)\,\iso\,\GL_n$ for the natural Levi factor of $P_G$. We equip $Q_G$ with the maximal torus $T_G$ of diagonal matrices and the Borel subgroup of upper-triangular matrices (the preimage of the latter in $P_G$ is a Borel subgroup, which yields our choice of positive roots). 

 Keep the notation of Section~4, in particular, $\cO=k[[t]]$ and $\Omega$ is the completed module of relative differential of $\cO$ over $k$. Set $L=L_0(\cO), U=U_0(\cO), V=V_0(\cO)$ and $M=L\oplus L^*\otimes\Omega$. The isomorphism $\cO\,\iso\,\Omega$ sending $1$ to $dt$ yields an isomorphism of group schemes $G\,\iso\, \Sp(M)$ over $\Spec\cO$. So, we offen think of $G$ as the group acting on $M$.   
  
 Set $\Upsilon=L^*\otimes V\otimes\Omega$ and $\Pi=U^*\otimes M$. 

\begin{Rem} In general, $L_0^*\otimes V_0$ is not a spherical $Q_G\times H$-variety. By (\cite{GN}, Theorem~3.2.1), in this case the set of $Q_G(\cO)\times H(\cO)$-orbits on $\Upsilon(F)$ is not countable. Indeed, already for the open $Q_G\times H$-orbit $(Q_G\times H)/R$ in $L_0^*\otimes V_0$, the set of $R(F)$-orbits on $\Gr_{Q_G\times H}$ is not countable. 

 Similarly, in general $U_0^*\otimes M_0$ is not a spherical $Q_H\times G$-variety, and  the set of $Q_H(\cO)\times G(\cO)$-orbits on $\Pi(F)$ is not countable. 
 \end{Rem}
 
\smallskip\noindent
6.2.1  As in Section~3.2, we define the functor
\begin{equation}
\label{zeta_geom} 
\zeta: \D_{(Q_G\times Q_H)(\cO)}(\Upsilon(F))\,\iso\, \D_{(Q_H\times Q_G)(\cO)}(\Pi(F))
\end{equation}
as the partial Fourier transform (\ref{Fourier_transform_for_def_equiv}) with respect to the decomposition $\Upsilon(F)\,\iso\, L^*\otimes U\otimes\Omega(F)\oplus  L^*\otimes U^*\otimes\Omega(F)$.  

\begin{Def} 
\label{Def_Weil_category}
The \select{Weil category for $G\times H$} is the category $\Weil_{G,H}$ of triples $(\cF_1,\cF_2,\beta)$, where $\cF_1\in \P_{(Q_G\times H)(\cO)}(\Upsilon(F))$, $\cF_2\in \P_{(Q_H\times G)(\cO)}(\Pi(F))$, and $\beta: \zeta(f(\cF_1))\,\iso\, f(\cF_2)$ is an isomorphism for the diagram
$$
\begin{array}{ccc}
\P_{(Q_G\times H)(\cO)}(\Upsilon(F)) && 
\P_{(Q_H\times G)(\cO)}(\Pi(F))\\
\downarrow\lefteqn{\scriptstyle f} && \downarrow\lefteqn{\scriptstyle f}\\
\P_{(Q_G\times Q_H)(\cO)}(\Upsilon(F)) & \toup{\zeta} & \P_{(Q_H\times Q_G)(\cO)}(\Pi(F)),
\end{array}
$$
where $f$ are forgetful functors. Write 
$$
f_G:\Weil_{G,H}\to \P_{(Q_H\times G)(\cO)}(\Pi(F))\;\;\;\; \mbox{and}\;\;\;\; f_H:\Weil_{G,H}\to \P_{(Q_G\times H)(\cO)}(\Upsilon(F))
$$ 
for the functors  
sending $(\cF_1,\cF_2,\beta)$ to $\cF_2$ and $\cF_1$ respectively. Write $\D\!\Weil_{G,H}$ for the category obtained by replacing in the above definition $\P$ by $\DP$ everywhere.  
\end{Def}

 By (\cite{G}, Lemma~4.8), both functors $f$ in the above diagram are fully faithful, and their image is stable under subquotients. It follows that $\Weil_{G,H}$ is abelian, and both $f_G$ and $f_H$ are fully faithful.  
Write $\Weil^{ss}_{G,H}\subset\Weil_{G,H}$ for the full subcategory of semi-simple objects. We write
$$
\D\!\Weil^{ss}_{G,H}\subset\D\!\Weil_{G,H}
$$ 
for the full subcategory of objects of the form $\oplus_{i\in\ZZ}\, K_i[i]$ with $K_i\in \Weil^{ss}_{G,H}$ for all $i$.

Since $\zeta(I_0)\,\iso\, I_0$ canonically, $I_0$ is naturally an object of $\Weil^{ss}_{G,H}$. Combining the decomposition theorem (\cite{BBD}, Corollary~5.4.6) with the fact that $G$ and $H$-actions in the Weil representation commute with each other, one gets the following. 
 
\begin{Pp} 
\label{Pp_rich_Hecke_functors}
There exist natural functors $\D\Sph_G\to \D\!\Weil^{ss}_{G,H}$ and $\D\Sph_H\to \D\!\Weil^{ss}_{G,H}$ such that the diagrams commute
$$
\begin{array}{ccccccc}
\D\Sph_G & \to & \D\!\Weil^{ss}_{G,H} && \D\Sph_H & \to & \D\!\Weil^{ss}_{G,H}\\
 & \searrow\lefteqn{\scriptstyle a_G} & \downarrow\lefteqn{\scriptstyle f_G} &&& \searrow\lefteqn{\scriptstyle a_H} & \downarrow\lefteqn{\scriptstyle f_H}\\
 && \DP^{ss}_{(Q_H\times G)(\cO)}(\Pi(F)) &&&& \DP^{ss}_{(Q_G\times H)(\cO)}(\Upsilon(F))
\end{array}
$$ 
Here the functor $a_G$ (resp., $a_H$) sends $\cT$ to $\H^{\la}_G(\cT, I_0)$ (resp., to
$\H^{\la}_H(\cT, I_0)$). 
\end{Pp}
\begin{Prf}
The arguments for both functors being similar, we give a proof only for the second one. Given $\cT\in\Sph_H$, by decomposition theorem $\H^{\la}_H(\cT, I_0)\in \D_{(Q_G\times H)(\cO)}(\Upsilon(F))$ identifies with the direct sum of its (shifted) perverse cohomology sheaves. It suffices to show that each perverse cohomology sheaf $K$ of $\zeta(f \H^{\la}_H(\cT, I_0))$ actually lies in the full subcategory $\P_{(Q_H\times G)(\cO)}(\Pi(F))$ of $\P_{(Q_H\times Q_G)(\cO)}(\Pi(F))$. 

 Denote by $P_G^- \subset G$ the parabolic subgroup preserving $L_0^*$, write $U_G^-$ for its unipotent radical. By Lemma~\ref{Lm_extending_actions_on_Schwarz}, it suffices to show that $K$ admits a $U_G(\cO)$ and $U_G^-(\cO)$-equivariant structures. For $v\in \Upsilon(F)$ write $s_{\cL}(v)$ for the composition
$$
\Sym^2 L(F)\,\toup{\Sym^2 v}\, \Sym^2 (V\otimes\Omega)(F)\to\Omega^2(F)
$$
Let $\Char(\Upsilon)\subset\Upsilon(F)$ be the ind-subscheme of $v\in \Upsilon(F)$ such that $s_{\cL}(v): \Sym^2 L\to \Omega^2$ is regular. The $U_G(\cO)$-equivariance of $K$ is equivalent to the fact that $\zeta^{-1}(K)$ is the extension by zero from $\Char(\Upsilon)$. But the complex $\H^{\la}_H(\cT, I_0)$ itself satisfies this property, so its direct summand  also does. 
 
  To get a $U_G^-(\cO)$-action on $K$, consider the commutative diagram
$$
\begin{array}{ccc}
\P_{(Q_G\times Q_H)(\cO)}(\Upsilon(F)) & \toup{\zeta} & \P_{(Q_H\times Q_G)(\cO)}(\Pi(F))\\
 & \searrow\lefteqn{\scriptstyle \Four_{\psi}} & \downarrow \lefteqn{\scriptstyle \zeta_1}\\
 && P_{(Q_G\times Q_H)(\cO)}(L\otimes V(F)),
\end{array}
$$
where $\zeta_1$ is a partial Fourier transform (with respect to $\psi$), and $\Four_{\psi}$ is the complete Fourier transform (cf. Remark~\ref{Rem_comp_partial_Fourier}). By Lemma~\ref{Lm_Fourier_Hecke_commute}, 
$$
\zeta_1\zeta (f \H^{\la}_H(\cT, I_0))\,\iso\,
\Four_{\psi}(f \H^{\la}_H(\cT, I_0))\,\iso\, f_1\H^{\la}_H(\cT, I_0),
$$
where we have denote by 
$f_1: \P_{(Q_G\times H)(\cO)}(L\otimes V(F))\to \P_{(Q_G\times Q_H)(\cO)}(L\otimes V(F))$ the forgetful functor. 

 Let $\Char(L\otimes V(F))\subset L\otimes V(F)$ be the ind-subscheme of $v\in L\otimes V(F)$ such that
the composition 
$\Sym^2 L^*\,\toup{\Sym^2 v}\, \Sym^2 V(F)\to F$
factors through $\cO\subset F$. Note that 
$$
\H^{\la}_H(\cT, I_0)\in \P_{(Q_G\times H)(\cO)}(L\otimes V(F))
$$ 
is the extension by zero from $\Char(L\otimes V(F))$. The $U^-_G(\cO)$-equivariance of $K$ is equivalent to the fact that $\zeta_1(K)$ is the extension by zero from $\Char(L\otimes V(F))$. We are done.
\end{Prf}

\medskip

\begin{Rem} 
\label{Rem_Aut_O_equiv_for_SpSO}
i) Arguing as in Proposition~\ref{Pp_rich_Hecke_functors}, one may   construct the commuting actions of $\D\Sph_G$ and $\D\Sph_H$ on $\D\Weil_{G,H}$, this is left to a reader.

\smallskip\noindent
ii) As in Section~4.5, one may strengthen Proposition~\ref{Pp_rich_Hecke_functors} saying that the functors $a_G$, $a_H$ actually take values in the $\Aut^0(\cO)$-equivariant versions of the corresponding categories.

\smallskip\noindent
iii) Arguing as in Proposition~\ref{Pp_rich_Hecke_functors}, one shows that $\Weil_{G,H}$ identifies with the category of 
$$
F\in P_{(Q_G\times Q_H)(\cO)}(\Upsilon(F)),
$$ 
which are supported on $Char(\Upsilon)$ and such that $\Four_{\psi}(F)\in \P_{(Q_G\times Q_H)(\cO)}(L\otimes V(F))$ is the extension by zero from $Char(L\otimes V(F))$. Note that $\Weil_{G,H}\hook{} P_{(Q_G\times Q_H)(\cO)}(\Upsilon(F))$ is a full subcategory stable under subquotients. 
\end{Rem}

\medskip\noindent
6.2.2. By abuse of notation, we simply write $\H^{\la}_G(\cT, I_0)\in \D\!\Weil^{ss}_{G,H}$ (resp., $\H^{\la}_H(\cT, I_0)\in \D\!\Weil^{ss}_{G,H}$) for $\cT\in\D\Sph_G$ (resp., $\cT\in\D\Sph_H$). 

 Recall the definition of the homomorphisms $\kappa$ from Section~2.4.2.
For $m\le n$ we have $\kappa: \check{H}\times\Gm\to\check{G}$. We write $\gRes^{\kappa}: \Sph_G\to\D\Sph_H$ for the corresponding geometric restriction functor. 

For $m>n$ we have $\kappa: \check{G}\times\Gm\to \check{H}$. Write $\gRes^{\kappa}: \Sph_H\to\D\Sph_G$ for the corresponding geometric restriction functor. 

Here is our main local result. 

\begin{Th} 
\label{Th_main_local_Sp_SO}
1) Assume $m\le n$. The functors $\Sph_G\to \D\!\Weil^{ss}_{G,H}$ given by
\begin{equation}
\label{iso_forThmainlocal_m_less_n}
\cS\mapsto \H^{\la}_G(\cS, I_0)\;\;\;\;\mbox{and}\;\;\;\;
\cS\mapsto \H^{\la}_H(\ast\gRes^{\kappa}(\cS), I_0)
\end{equation}
are isomorphic. This isomorphism is compatible with the tensor structures on $\Sph_G,\Sph_H$. 

2) Assume $m>n$. The two functors $\Sph_H\to \D\!\Weil^{ss}_{G,H}$ given by
\begin{equation}
\label{iso_forThmainlocal_m_greater_n}
\cT\mapsto \H^{\la}_H(\cT, I_0)\;\;\;\;\mbox{and}\;\;\;\;
\cT\mapsto \H^{\la}_G(\gRes^{\kappa}(\ast\cT), I_0)
\end{equation}
are isomorphic. This isomorphism is compatible with the tensor structures on $\Sph_G,\Sph_H$. 
\end{Th}

 The proof will be given in Section~6.4. 
 
\medskip\noindent
6.3.1 In this subsection we assume $m\le n$ and analyse the action of $\Sph_H$ on $\D_{(Q_G\times H)(\cO)}(\Upsilon(F))$ in more details. 

Write $V^{\check{\lambda}}$ for the irreducible $H$-module with h.w. $\check{\lambda}\in\check{\Lambda}^+_H$. 
For $1\le i< m$ let $\check{\alpha}_i\in\check{\Lambda}^+_H$ denote the h.w. of the $H$-module $\wedge^i V_0$. 
Recall that 
$$
\wedge^m V_0\,\iso\, V^{\check{\alpha}_m}\oplus V^{\check{\alpha}'_m}
$$ 
is a direct sum of two irreducible representations, this is our definition of $\check{\alpha}_m, \check{\alpha}'_m$. Say that a maximal isotropic subspace $\cL\subset V_0$ is \select{$\check{\alpha}_m$-oriented} (resp., \select{$\check{\alpha}'_m$-oriented}) if $\wedge^m\cL\subset V^{\check{\alpha}_m}$ (resp., $\wedge^m\cL\subset V^{\check{\alpha}'_m}$). The group $H$ has two orbits on the scheme of maximal isotropic subspaces in $V_0$ given by their orientation. 

 For $v\in \Upsilon(F)$ let $s_{\cL}(v):\Sym^2 L\to \Omega^2(F)$ be the composition 
$$
\Sym^2 L\toup{\Sym^2 v} \Sym^2(V\otimes\Omega)(F)\to\Omega^2(F)
$$ 
For $\lambda\in\Lambda^+_H$ let $N=\<\lambda,\check{\alpha}_1\>$, define a closed subscheme $_{\lambda}\Upsilon\subset {_N\Upsilon}=t^{-N}\Upsilon$ as follows. A point $v\in {_N\Upsilon}$ lies in $_{\lambda}\Upsilon$ iff the following conditions hold:

\begin{itemize}
\item[C1)] $s_{\cL}(v):\Sym^2 L\to\Omega^2$ is regular;
\item[C2)] for $1\le i<m$ the map $\wedge^i L\;\toup{\wedge^i v}\; (\Omega^i\otimes \wedge^i V)(\<-w_0(\lambda), \check{\alpha}_i\>)$ is regular;
\item[C3)] the map $\wedge^m L\;\toup{v_m\oplus v'_m}\; (\Omega^m\otimes V^{\check{\alpha}_m})(\<-w_0(\lambda),\check{\alpha}_m\>)\oplus (\Omega^m\otimes V^{\check{\alpha}'_m})(\<-w_0(\lambda),\check{\alpha}'_m\>)$ induced by $\wedge^m v$ is regular. 
\end{itemize}
The scheme $_{\lambda}\Upsilon$ is stable under translations by $_{-N}\Upsilon$, so there is a closed subscheme $_{\lambda,N}\Upsilon\subset {_{N,N}\Upsilon}$ whose preimage under the projection $_N\Upsilon\to {_{N,N}\Upsilon}$ is $_{\lambda}\Upsilon$. 

 As in Section~4.4, we have a map
$\pi: {_{0,N}\Upsilon}\ttimes\ov{\Gr}^{\lambda}_H\to {_{N,N}\Upsilon}$,  
it factors through the closed immersion $_{\lambda,N}\Upsilon\hook{} {_{N,N}\Upsilon}$, and 
$$
\H^{\lambda}_H(I_0)\,\iso\, \pi_!(\Qlb\tboxtimes \cA^{\lambda}_H)[\dim {_{0,N}\Upsilon}] \in \D_{(Q_G\times H)(\cO)}(_{N, N}\Upsilon)
$$ 

 Let $\Char(\Upsilon)\subset \Upsilon(F)$ be the ind-subscheme of $v\in \Upsilon(F)$ satisfying C1). Note that $\Char(\Upsilon)$ is preserved by the $H(F)$-action. For $v\in \Char(\Upsilon)$ let $L_v=v(L)+V\otimes\Omega$ and 
$$
L_v^{\perp}=\{v\in V\otimes\Omega\mid \<v, u\>\in\Omega^2\;\;\;\mbox{for any}\;\; u\in L_v\}
$$
Let $V_v\subset V(F)$ be defined by $V_v\otimes\Omega=v(L)+L_v^{\perp}$, then $V_v$ is an orthogonal lattice in $V(F)$, that is, a point of $\Gr_H$. We stratify $\Char(\Upsilon)$ by locally closed subschemes $_{\lambda}\Char(\Upsilon)$ indexed by $\lambda\in\Lambda^+_H$. Namely, for $v\in\Char(\Upsilon)$ we let $v\in {_{\lambda}\Char(\Upsilon)}$ iff $V_v\in\Gr^{\lambda}_H$. 

 We have $_{\lambda}\Char(\Upsilon)\subset {_{\lambda}\Upsilon}$. Moreover, for $N=\<\lambda,\check{\alpha}_1\>$ there is a unique open subscheme $_{\lambda,N}\Upsilon^0\subset {_{\lambda,N}\Upsilon}$ whose preimage under the projection $_{\lambda}\Upsilon\to {_{\lambda,N}\Upsilon}$ identifies with $_{\lambda}\Char(\Upsilon)$. 
 
  Write $\IC(_{\lambda,N}\Upsilon^0)\in\P_{(Q_G\times H)(\cO)}(_{N,N}\Upsilon)$ for the intersection cohomology sheaf of $_{\lambda,N}\Upsilon^0$.
 
\begin{Lm} 
\label{Lm_very_nice_m_less_n}
Assume $m\le n$. 
1) The map 
$$
\pi: {_{0,N}\Upsilon}\ttimes\ov{\Gr}^{\lambda}_H\to {_{\lambda,N}\Upsilon}
$$ 
is an isomorphism over the open subscheme $_{\lambda,N}\Upsilon^0$. So, $\dim {_{\lambda,N}\Upsilon^0}=2Nnm+\<\lambda, 2\check{\rho}_H\>$.\\
2) For $\lambda\in \Lambda^+_H$ we have $\H^{\lambda}_H(I_0)\,\iso\,\IC(_{\lambda,N}\Upsilon^0)$ canonically.
\end{Lm}
\begin{Prf} 1) The fibre of $\pi$ over $v\in {_{\lambda,N}\Upsilon^0}$ is the scheme classifying orthogonal lattices $V'\subset V(F)$ such that $V'\in\ov{\Gr}^{\lambda}_H$ and $v(L)\subset V'\otimes\Omega$. Given such lattice $V'$, the inclusion $v(L)+V\otimes\Omega\subset V'\otimes\Omega+V\otimes\Omega$ must be an equality, because for $V'\in\Gr_H^{\mu}$ with $\mu\le\lambda$ we get
$$
\dim(V'+V)/V=\epsilon(\mu)\le \epsilon(\lambda)=
\dim(v(L)+V\otimes\Omega/V\otimes\Omega)
$$
We have set here $\epsilon(\mu)=
\max\{\<\mu,\check{\alpha}_m\>, \<\mu, \check{\alpha}'_m\>\}$. Thus, $V_v\in\Gr_H^{\lambda}$ is the unique preimage of $v$ under $\pi$. The first assertion follows. 

 For the convenience of the reader recall that $T_H\,\iso\,\Gm^m$ is the torus of diagonal matrices in $\GL(U_0)$, and 
\begin{equation}
\label{dom_coweights_for_SO} 
 \Lambda^+_H=\{\mu=(a_1\ge\ldots\ge a_m)\in\ZZ^m\mid a_{m-1}\ge \mid\! a_m\!\mid \}
\end{equation} 
In these notations one may take $\check{\alpha}_m=(1,\ldots,1)$ and $\check{\alpha}'_m=(1,\ldots,1,-1)$. So, if $\mu=(a_1,\ldots, a_m)\in \Lambda^+_H$ then $\epsilon(\mu)=a_1+\ldots+a_{m-1}+\mid\! a_m\!\mid$.

\smallskip\noindent
2) From 1) we learn that $\IC(_{\lambda,N}\Upsilon^0)$ appears in $\H^{\lambda}_H(I_0)$ with multiplicity one. So, it suffices to show that 
$$
\Hom(\H^{\lambda}_H(I_0), \H^{\lambda}_H(I_0))=\Qlb,
$$
where $\Hom$ is taken in the homotopy category of $\D_{(Q_G\times H)(\cO)}(\Upsilon(F))$. By adjointness, 
$$
\Hom(\H^{\lambda}_H(I_0), \H^{\lambda}_H(I_0))\,\iso\, 
\Hom(\H^{-w_0(\lambda)}_H \H^{\lambda}_H(I_0), I_0).
$$
So, it suffices to show that for any $\lambda\in\Lambda^+_H$ with $\lambda\ne 0$ one has 
$\Hom(\H^{\lambda}_H(I_0), I_0)=0$ in the homotopy category of $\D_{(Q_G\times H)(\cO)}(\Upsilon(F))$. As above, set $N=\<\lambda,\check{\alpha}_1\>$. We will show that 
$$
\Hom(\H^{\lambda}_H(I_0), I_0)=0
$$ 
in the homotopy category of $\D_{(Q_G\times H)(\cO/t^d)}(_{N,N}\Upsilon(F))$ for $d$ large enough.  

 Let $i: {_{0,N}\Upsilon}\to {_{N,N}\Upsilon}$ denote the natural closed immersion. Recall that $I_0=i_!p^!\Qlb[-2Nnm]$ on $_{N,N}\Upsilon$, where $p: {_{0,N}\Upsilon}\to\Spec k$ is the projection. By adjointness, we are reduced to show that 
$$
\Hom(p_!i^*\H^{\lambda}_H(I_0)[2Nnm], \Qlb)=0
$$ 
in the homotopy category of $\D_{(Q_G\times H)(\cO/t^d)}(\Spec k)$. It suffices to show that $p_!i^*\H^{\lambda}_H(I_0)[2Nnm]$ is placed in degrees $<0$.   

 Denote by $\cY^{\lambda}$ the preimage of $_{0,N}\Upsilon$ under 
$$
\pi: {_{0,N}\Upsilon\ttimes \ov{\Gr}^{\lambda}_H}\to {_{N,N}\Upsilon}
$$ 
Then $\cY^{\lambda}$ is the scheme classifying $V'\in\ov{\Gr}^{\lambda}_H$ and $v\in {_{0,N}\Upsilon}$ such that $v(L)\subset (V'/t^NV)\otimes\Omega$. 
 
 Stratify $\cY^{\lambda}$ by locally closed subschemes $\cY^{\lambda,\mu}$ indexed by $\mu\in\Lambda^+_H$ with $\mu\le\lambda$. The subscheme $\cY^{\lambda,\mu}\subset \cY^{\lambda}$ is given by the condition $V'\in\Gr^{\mu}_H$. 
 
 Recall that $\H^{\lambda}_H(I_0)=\pi_!(I_0\tboxtimes\cA^{\lambda}_H)$, where $I_0\tboxtimes\cA^{\lambda}_H$ is perverse. It remains to show that for each stratum $\cY^{\lambda,\mu}$ the complex 
\begin{equation}
\label{complex_contribution_cY_lambda_mu}
\RG_c(\cY^{\lambda,\mu}, (I_0\tboxtimes\cA^{\lambda}_H)\mid_{\cY^{\lambda,\mu}})[2Nnm]
\end{equation}
is placed in degrees $<0$. The key observation is that the map $\cY^{\lambda,\mu}\to \Gr_H^{\mu}$ sending $(v,V')$ to $V'$ is a vector bundle, its rank equals $n(2mN-\epsilon(\mu))$. Here $\epsilon(\mu)$ is the expression defined in 1). Indeed, for any lattice $V'\in \Gr^{\mu}_H$, the fibre of this vector bundle over $V'$ is 
$$
\Hom_{\cO}(L, ((V'\cap V)/t^NV)\otimes\Omega)
$$
and $\dim_k (V/(V'\cap V))=\dim(V+V')/V=\epsilon(\mu)$. 
So, (\ref{complex_contribution_cY_lambda_mu}) identifies with
\begin{equation}
\label{complex_already_on_Gr_mu}
\RG_c(\Gr_G^{\mu}, \cA^{\lambda}_H\mid_{\Gr^{\mu}_H})[2n\epsilon(\mu)].
\end{equation}

 By definition of the intersection cohomology sheaf, $\cA^{\lambda}_H\mid_{\Gr^{\mu}_H}$ has usual cohomology sheaves in degrees $\le -\<\mu,2\check{\rho}_H\>$, and the inequality is strict unless $\mu=\lambda$. So, (\ref{complex_already_on_Gr_mu}) is placed in degrees $\le \<\mu,2\check{\rho}_H\>-2n\epsilon(\mu)$, and the inequality is strict unless $\mu=\lambda$. 
 
  One checks that for any $\tau\in\Lambda^+_H$ we have 
$\<\tau,2\check{\rho}_H\>-2m\epsilon(\tau)\le 0$, and the inequality is strict unless $\tau=0$. Our assertion follows, because $n\ge m$.   

 For the convenience of the reader note that in the notation (\ref{dom_coweights_for_SO}) for $\tau=(a_1\ge\ldots\ge a_m)\in \Lambda^+_H$ one gets
$$
\<\tau,2\check{\rho}_H\>-2m\epsilon(\tau)=-2a_1-4a_2-\ldots -(2m-2)a_{m-1}-2m\!\mid\! a_m\!\mid
$$ 
\end{Prf}
 
\medskip

 Write $\cP$ for the composition of functors
$$
\D\Weil_{G,H}\;\toup{f_H}\; \D\P_{(Q_G\times H)(\cO)}(\Upsilon(F))\toup{J^*_{P_H}}
 \D^b_{(Q_G\times Q_H)(\cO)}(U\otimes L^*\otimes\Omega(F))
$$  
Let $k_0\subset k$ be a finite subfield. Assume that 
all the objects of Section~6.1 are defined over $k_0$. 
Set $\cO_0=k_0[[t]]$ and $F_0=k_0((t))$. 

 Write $\D\!\Weil_{G,H,k_0}$ for the category of triples $(\cF_1, \cF_2,\beta)$ as in Definition~\ref{Def_Weil_category}, where now 
$$
\cF_1\in \D^b_{(Q_G\times H)(\cO_0)}(\Upsilon(F_0)), \;\;\;\cF_2\in \D^b_{(Q_H\times G)(\cO_0)}(\Pi(F_0))
$$ 
are pure complexes of weight zero, and 
and $\beta: \zeta(f(\cF_1))\,\iso\, f(\cF_2)$ is an isomorphism as above (it is understood that (\ref{zeta_geom}) is normalized to preserve purity). Let $\Weil_{G,H,k_0}\subset \D\!\Weil_{G,H,k_0}$ be the full subcategory given by the condition that $\cF_i$ is perverse.

 We have a natural functor $\iota: \D\!\Weil_{G,H,k_0}\to \D\!\Weil_{G,H}$ that restricts to $\iota: \Weil_{G,H,k_0}\to\Weil_{G,H}$. 

 Note that any object of $\Weil_{G,H}$ is $\Gm$-equivariant with respect to the homotheties on $L^*\otimes V\otimes\Omega(F)$, because it is equivariant with respect to the action of the center of $Q_G(\cO)$. Write 
\begin{equation}
\label{functor_cP_0}
\cP_0: \D\!\Weil_{G,H,k_0}\to  \D^b_{(Q_G\times Q_H)(\cO_0), mixed}(U\otimes L^*\otimes\Omega(F_0))
\end{equation}
for the natural lifting of the functor $\cP$, here we have denoted by 
$$
\D^b_{(Q_G\times Q_H)(\cO_0), mixed}(U\otimes L^*\otimes\Omega(F_0))\subset \D^b_{(Q_G\times Q_H)(\cO_0)}(U\otimes L^*\otimes\Omega(F_0))
$$
the full subcategory of mixed complexes (\cite{BBD}, 5.1.5). By Corollary~\ref{Cor_actually_used_Jacquet}, if $K\in \D\!\Weil_{G,H,k_0}$ then $\cP_0(K)$ is pure of weight zero. The Grothendieck group of $\D^b_{(Q_G\times Q_H)(\cO_0), mixed}(U\otimes L^*\otimes\Omega(F_0))$ is denoted by $K'_{k_0}$. Remind that $K'_{k_0}$ is $\ZZ$-graded by weights, and its component of weight $i$ is a free abelian group with base consisting of irreducible $(Q_G\times Q_H)(\cO_0)$-equivariant pure perverse sheaves of weight $i$ on $U\otimes L^*\otimes\Omega(F_0)$.
 
\begin{Pp} 
\label{Pp_functor_cP_injection}
For $i=1,2$ let $K_i\in \D\!\Weil_{G,H, k_0}$. If $\cP_0(K_1)=\cP_0(K_2)$ in $K'_{k_0}$ then 
$\iota(K_1)\,\iso\, \iota(K_2)$ in $\D\!\Weil_{G,H}$. 
\end{Pp}
\begin{Prf} 
Write $K_{k_0}$
for the Grothendieck group of the category $\D\!\Weil_{G,H,k_0}$. 
It is $\ZZ$-graded by weights, and its component of weight $i$ is a free abelian group with base consisting of triples $(\cF_1,\cF_2,\beta)$ as in the definition of $\Weil_{G,H,k_0}$ such that $\cF_i$ is an irreducible perverse sheaf pure of weight $i$.
By Corolary~\ref{Cor_Hecke_GL_n}, 
\begin{equation}
\label{category_over_k_0_Levi}
\P^{ss}_{(Q_G\times Q_H)(\cO)}(U\otimes L^*\otimes\Omega(F))
\end{equation}
identifies with $\Sph_{Q_H}$ (resp., with $\Sph_{Q_G}$)
for $m\le n$ (resp., for $m>n$). All the objects of the latter category are defined over $k_0$ and as such are pure of weight zero. 

 The functor (\ref{functor_cP_0}) yields a homomorphism $J^*_{P_H}: K_{k_0}\to K'_{k_0}$. Let us show that it is injective. Let $F$ be an object in its kernel. For any finite subfield $k_0\subset k_1\subset k$ we have the $\Qlb$-vector space $\Weil_{G,H}(k_1)$ introduced in Section~3.2. The map $\tr_{k_1}$ trace of Frobenius over $k_1$ fits into the diagram
$$
\begin{array}{ccc}
K_{k_0} & \toup{J^*_{P_H}} & K'_{k_0}\\
\downarrow\lefteqn{\scriptstyle \tr_{k_1}} && \downarrow\lefteqn{\scriptstyle \tr_{k_1}}\\
\Weil_{G,H}(k_1) & \toup{J_{k_1}} & \Funct_{k_1},  
\end{array}
$$
where $\Funct_{k_1}$ is the non ramified Hecke algebra $\cH(Q_H)$ (resp., $\cH(Q_G)$) of $Q_H$ for $m\le n$ (resp., of $Q_G$ for $m>n$). By Proposition~\ref{Pp_local_main_classical}, the low horizontal arrow $J_{k_1}$ is injective. So, $\tr_{k_1}(F)=0$ for any finite extension $k_0\subset k_1$. By a result of Laumon (\cite{Lam}, Theorem~1.1.2) this implies $F=0$ in $K_{k_0}$.  
Since $K_i$ is pure of weight zero, the semi-simplifications of $K_1$ and $K_2$ are isomorphic in $\D\Weil_{G,H,k_0}$ by Remark~\ref{Rem_semi-simplifications_iso} below, and $\iota(K_1)\,\iso\, \iota(K_2)$ in $\D\!\Weil_{G,H}$.
\end{Prf} 
 
\medskip\noindent
\begin{Rem} 
\label{Rem_semi-simplifications_iso} Let $Y$ be a scheme of finite type over $k_0$ and $K_1,K_2$ pure complexes of weight zero in $\D^b(Y)$. Write $K_{mixed}(Y)$ for the Grothendieck group of the subcategory $\D^b_{mixed}(Y)\subset \D^b(Y)$ of mixed complexes.
If $K_1=K_2$ in $K_{mixed}(Y)$ then the semi-simplifications of $K_1$ and $K_2$ are isomorphic in $\D^b(Y)$. Indeed, $K_{mixed}(Y)$ is $\ZZ$-graded by weights, and its component of weight $i$ is a free abelian group with a base consisting of irreducible perverse sheaves pure of weight $i$ on $Y$.
\end{Rem}  
 
\begin{Con} 
Assume $m\le n$. The irreducible objects of  
$\Weil_{G,H}$ are exactly $\IC(_{\lambda,N}\Upsilon^0)$, $\lambda\in\Lambda^+_H$. 
The functor $\cT\mapsto \H^{\la}_H(\cT, I_0)$ yields an equivalence of categories
$$
\Sph_H\,\iso\, \Weil^{ss}_{G,H}
$$
\end{Con}

\medskip\noindent 
6.3.2 In this subsection we assume $m> n$ and analyse the action of $\Sph_G$ on $\D_{(Q_H\times G)(\cO)}(\Pi(F))$ in more details. 

 Let $\check{\omega}_i\in\check{\Lambda}^+_G$ be the highest weight of the fundamental representation of $G$ that appear in $\wedge^i M_0$ (all the weights of $\wedge^i M_0$ are less of equal to $\check{\omega}_i$). For $v\in \Pi(F)$ write $s_{\cU}(v):\wedge^2 U(F)\to \Omega(F)$ for the composition
$$
\wedge^2 U(F)\toup{\wedge^2 v} \wedge^2 M(F)\to\Omega(F)
$$ 
Write $\Char(\Pi)\subset\Pi(F)$ for the ind-subscheme of $v\in\Pi(F)$ such that $s_{\cU}(v):\wedge^2 U\to \Omega$ is regular. 

 For $\lambda\in\Lambda^+_G$ let $N=\<\lambda,\check{\omega}_1\>$, define the closed subscheme $_{\lambda}\Pi\subset {_N\Pi}=t^{-N}\Pi$ as follows. 
A point $v\in {_N\Pi}$ lies in $_{\lambda}\Pi$ if the following conditions hold:

\begin{itemize}
\item[C1)] $v\in\Char(\Pi)$;
\item[C2)] for $i=1,\ldots,n$ the map $\wedge^i U\, \toup{\wedge^i v}\, (\wedge^i M)(-\<w_0(\lambda), \check{\omega}_i\>)$ is regular.
\end{itemize}

 The scheme $_{\lambda}\Pi$ is stable under the translations by $t^N\Pi$, so there is a closed subscheme $_{\lambda,N}\Pi\subset {_{N,N}\Pi}$ such that $_{\lambda}\Pi$ is the preimage of $_{\lambda,N}\Pi$ under the projection $_N\Pi\to {_{N,N}\Pi}$. As in Section~4.4, we have a map
\begin{equation}
\label{map_proper_for_Hecke_G}
\pi: {_{0,N}\Pi}\ttimes \ov{\Gr}_G^{\lambda}\to {_{N,N}\Pi}
\end{equation}
and, by definition,
$$
\H^{\lambda}_G(I_0)\,\iso\, \pi_!(\Qlb\tboxtimes\cA^{\lambda}_G)[\dim {_{0,N}\Pi}]\in \D_{(Q_H\times G)(\cO)}(_{N,N}\Pi)
$$
Since all the weights of the $G$-module $\wedge^i M_0$ are less or equal to $\check{\omega}_i$, the map (\ref{map_proper_for_Hecke_G}) factors through the closed subscheme $_{\lambda,N}\Pi\hook{} {_{N,N}\Pi}$.  
 
 For $v\in\Char(\Pi)$ let $U_v=v(U)+M$ and 
$$
U_v^{\perp}=\{m\in M(F)\mid \<m, m_1\>\in\Omega\;\;\mbox{for any}\;\; m_1\in U_v\}
$$
Let $M_v=v(U)+U_v^{\perp}$. Note that $U_v/U_v^{\perp}$ is naturally a symplectic vector space, and $M_v/U_v^{\perp}\subset  U_v/U_v^{\perp}$ is a lagrangian subspace. 
So, $M_v\subset M(F)$ is a symplectic lattice, that is, $M_v\in\Gr_G$. 

 Stratify $\Char(\Pi)$ by locally closed subschemes $_{\lambda}\Char(\Pi)$ indexed by $\lambda\in\Lambda^+_G$. Namely, for $v\in\Char(\Pi)$ we let $v\in {_{\lambda}\Char(\Pi)}$ iff $M_v\in\Gr^{\lambda}_G$. The condition $M_v\in\Gr^{\lambda}_G$ is also equivalent to requiring that there is an isomorphism of $\cO$-modules
$$
U_v/M\,\iso\, \cO/t^{a_1}\oplus\ldots\oplus\cO/t^{a_n}
$$
for $\lambda=(a_1\ge\ldots\ge a_n\ge 0)\in\Lambda^+_G$. So, the stratification in question is by the isomorphism classes of the $\cO$-module $U_v/M$. 

 Clearly, $_{\lambda}\Char(\Pi)\subset {_{\lambda}\Pi}$, and there is a unique open subscheme $_{\lambda,N}\Pi^0\subset {_{\lambda,N}\Pi}$ whose preimage under the projection $_{\lambda}\Pi\to {_{\lambda,N}\Pi}$ identifies with $_{\lambda}\Char(\Pi)$. 

 Write $\IC(_{\lambda,N}\Pi^0)$ for the intersection cohomology sheaf of $_{\lambda,N}\Pi^0$. 
 
\begin{Lm} Assume $m> n$. 1) For any $\lambda\in\Lambda^+_G$ the map 
$$
\pi: {_{0,N}\Pi}\ttimes \ov{\Gr}_G^{\lambda}\to {_{\lambda,N}\Pi}
$$
is an isomorphism over the open subscheme $_{\lambda,N}\Pi^0$. So, $\dim {_{\lambda,N}\Pi^0}=2Nmn+\<\lambda, 2\check{\rho}_G\>$.\\
2) For $\lambda\in\Lambda^+_G$ we have $\H^{\lambda}_G(I_0)\,\iso\, \IC(_{\lambda,N}\Pi^0)$ canonically.
\end{Lm}
\begin{Prf}
1) The fibre of $\pi$ over $v\in {_{\lambda,N}\Pi^0}$ is the scheme classifying symplectic lattices $M'\subset M(F)$ such that $M'\in\ov{\Gr}^{\lambda}_G$ and $v(U)\subset M'$. Given such lattice $M'$, the inclusion $U_v\subset M'+M$ must be an equality, because for $M'\in \Gr_G^{\mu}$ with $\mu\le\lambda$ we get
$$
\dim(M'+M/M)=\epsilon(\mu)\le \epsilon(\lambda)=\dim(U_v/M).
$$
We have denoted here $\epsilon(\mu):=\<\mu, \check{\omega}_n\>$ for $\mu\in\Lambda^+_G$. 
Thus, $M'=M_v$ is the unique preimage of $v$ under $\pi$. 
The first assertion follows.

\smallskip\noindent
2) is completely analogous to the proof of the second part of Lemma~\ref{Lm_very_nice_m_less_n}.
\end{Prf}

\begin{Con} 
Assume $m>n$. The irreducible objects of $\Weil_{G,H}$ are exactly $\IC(_{\lambda,N}\Pi^0)$, $\lambda\in\Lambda^+_G$. The functor $\cT\mapsto \H^{\la}_G(\cT, I_0)$ yields an equivalence of categories
$$
\Sph_G\,\iso\, \Weil^{ss}_{G,H}
$$
\end{Con}

\begin{Rem} i) For $n=1$ and $m>n$ the isomorphism
$\H^{\lambda}_G(I_0)\,\iso\, \IC(_{\lambda,N}\Pi^0)$ for 
$\lambda\in\Lambda^+_G$ can also be obtained from Proposition~\ref{Pp_action_Hecke_GL_n_smallest}. Indeed, in this case $\Gr_G$ identifies with a connected component of $(\Gr_{\GL_2})_{red}$. The desired irreducibility of $\H^{\lambda}_G(I_0)$ becomes a particular case of Proposition~\ref{Pp_action_Hecke_GL_n_smallest}.\\
ii) For $m=1$ and $m\le n$ it is evident that $\H^{\lambda}_H(I_0)\,\iso\, \IC(_{\lambda,N}\Upsilon^0)$ for $\lambda\in\Lambda^+_H$.
\end{Rem}

\medskip\noindent
6.4 \select{Proof of Theorem~\ref{Th_main_local_Sp_SO}}

\Step 1 The following property of the Fourier transform functors  follows from base change for proper morphisms. If $\cV\to S\gets \cV^*$ is a diagram of dual vector bundles over a scheme $S$, let $\cV'\to S'\gets \cV'^*$ be the diagram obtained from it by the base change with respect to a closed immersion $S'\hook{} S$. Then for the inclusions $i_1: \cV'\hook{} \cV$ and $i_2: \cV'^*\hook{} \cV^*$ we have $i_2^*\comp\Four_{\psi}\,\iso\, \Four_{\psi}\comp i_1^*$.
Thus, the following diagram of functors commutes
$$
\begin{array}{ccccc}
&& \D\!\Weil_{G,H} \\
& \swarrow\lefteqn{\scriptstyle f_H} && \searrow\lefteqn{\scriptstyle f_G}\\
\DP_{(Q_G\times H)(\cO)}(\Upsilon(F)) &&&& \DP_{(Q_H\times G)(\cO)}(\Pi(F))\\
\downarrow\lefteqn{\scriptstyle J_{P_H}^*} &&&& \downarrow\lefteqn{\scriptstyle J_{P_G}^*}\\
\D^b_{(Q_G\times Q_H)(\cO)}(U\otimes L^*\otimes\Omega(F)) && \toup{\Four_{\psi}} && 
\D^b_{(Q_G\times Q_H)(\cO)}(U^*\otimes L(F))
\end{array}
$$

 Let $\kappa_H:\check{Q}_H\times\Gm\to \check{H}$ be the map, whose first component $\check{Q}_H\to\check{H}$ is the natural inclusion, and second component $\Gm\to\check{H}$ is
$2(\check{\rho}_H-\check{\rho}_{Q_H})-n\check{\omega}_m$. Here $\check{\omega}_m$ is the h.w. of the $Q_H$-module $\det U_0$. The corresponding geometric restriction functor is denoted by
$$
\gRes^{\kappa_H}:\Sph_H\to\D\Sph_{Q_H}
$$ 

 Let $\kappa_G: \check{Q}_G\times\Gm\to\check{G}$ be the map, whose first component is the natural inclusion
$\check{Q}_G\hook{} \check{G}$, and the second component is $2(\check{\rho}_G-\check{\rho}_{Q_G})-m\check{\omega}_n$. Here $\check{\omega}_n$ is the h.w. of the $Q_G$-module $\det L_0$. The corresponding geometric restriction functor is denoted by
$$
\gRes^{\kappa_G}:\Sph_G\to\D\Sph_{Q_G}
$$ 
Note that $J^*_{P_H}(I_0)\,\iso\, I_0$, $J^*_{P_G}(I_0)\,\iso\, I_0$ and $\Four_{\psi}(I_0)\,\iso\, I_0$ canonically. 
 
  Assume that all the objects of Section~6.1 are defined over some finite subfield $k_0\subset k$. Remind our notation $\cO_0=k[[t]]$ and $F_0=k_0((t))$.   
Combining Lemmas~\ref{Lm_Jacquet_functors}, \ref{Lm_splitting_unique}, \ref{Lm_Fourier_Hecke_commute} and Corollary~\ref{Cor_actually_used_Jacquet}
for $\cT\in\Sph_H$ and $\cS\in \Sph_G$ we get canonical isomorphisms
$$
\Four_{\psi} J^*_{P_H} \H^{\la}_H(\cT, I_0)\,\iso\, \H^{\la}_{Q_H}(\gRes^{\kappa_H}(\cT), I_0)
$$
and
$$
J^*_{P_G}\H^{\la}_G(\cS, I_0)\,\iso\, \H^{\la}_{Q_G}(\gRes^{\kappa_G}(\cS), I_0)
$$
in $\D^b_{(Q_G\times Q_H)(\cO_0), mixed}(U^*\otimes L(F_0))$.   

Remind the functor $\cP_0$ given by (\ref{functor_cP_0}). To summarize, for 
$\cT\in\Sph_H$ and $\cS\in \Sph_G$ we get canonical isomorphisms
\begin{equation}
\label{iso_H_for_Th_local}
\cP_0 (\H^{\la}_H(\cT, I_0))\,\iso\, \H^{\la}_{Q_H}(\gRes^{\kappa_H}(\cT), I_0)
\end{equation}
and 
\begin{equation}
\label{iso_G_for_Th_local}
\cP_0 (\H^{\la}_G(\cS, I_0))\,\iso\, \H^{\la}_{Q_G}(\gRes^{\kappa_G}(\cS), I_0)
\end{equation}
in $\D^b_{(Q_G\times Q_H)(\cO_0), mixed}(U\otimes L^*\otimes\Omega(F_0))$. 

\medskip

\Step 2 
CASE $m\le n$. Let 
$$
\kappa_Q: \check{Q}_H\times\Gm\to \check{Q}_G\times\Gm
$$ 
be the map whose second component $\check{Q}_H\times\Gm\to \Gm$ is the projection, and 
the first component $\check{Q}_H\times\Gm\to \check{Q}_G$ is the composition
$$
\check{Q}_H\times\Gm\toup{\id\times 2\check{\rho}_{\GL_{n-m}}} \check{Q}_H\times\GL_{n-m}\hook{\on{Levi}} \check{Q}_G
$$
Write $\gRes^{\kappa_Q}:\D\Sph_{Q_G}\to\D\Sph_{Q_H}$ for the corresponding geometric restriction functor. 

 Now (\ref{iso_G_for_Th_local}) and Corollary~\ref{Cor_used_GL_n_GL_m} yield for $\cS\in\Sph_G$ canonical isomorphisms
$$
\cP_0 (\H^{\la}_G(\cS, I_0))\,\iso\, \H^{\la}_{Q_G}(\gRes^{\kappa_G}(\cS), I_0)\,\iso\, \H^{\la}_{Q_H}(\gRes^{\kappa_Q}(\ast\gRes^{\kappa_G}(\cS)), I_0)
$$
(they could depend on the decomposition $L_0\,\iso\, L_1\oplus L_2$ with $\dim L_1=m$, which is assumed fixed). On the other hand, (\ref{iso_H_for_Th_local}) yields canonical isomorphisms
$$
\cP_0(\H^{\la}_H(\ast\gRes^{\kappa}(\cS), I_0))\,\iso\, \H^{\la}_{Q_H}(\gRes^{\kappa_H}(\ast\gRes^{\kappa}(\cS)), I_0)
$$

 Let $\sigma: \check{Q}_G\,\iso\,\check{Q}_G$ be the automorphism sending $g$ to $^tg^{-1}$ for $g\in \check{Q}_G=\GL_n$. The restriction functor with respect to $\sigma\times\id: \check{Q}_G\times\Gm\,\iso\, \check{Q}_G\times\Gm$ identifies with $\ast:\D\Sph_{Q_G}\,\iso\,\D\Sph_{Q_G}$. 

 We will define an automorphism $\sigma_H$ of $\check{H}$ inducing $\ast: \Rep(\check{H})\,\iso\,\Rep(\check{H})$ and $\kappa$ making  the following diagram commmutative
\begin{equation}
\label{diag_crucial_m_less_n}
\begin{array}{ccccc}
\check{H}\times\Gm & \toup{\sigma_H\times\id} &
\check{H}\times\Gm & \toup{\kappa} & \check{G}\\
\uparrow\lefteqn{\scriptstyle \kappa_H}&& && \uparrow\lefteqn{\scriptstyle \kappa_G}\\
\check{Q}_H\times\Gm & \toup{\kappa_Q} & \check{Q}_G\times\Gm & \toup{\sigma\times\id} & \check{Q}_G\times\Gm,
\end{array}
\end{equation}
This yields for $\cS\in\Sph_G$ a canonical isomorphism
\begin{equation}
\label{iso_for_Step2_Thmainlocal_m_less_n}
\cP_0 \H^{\la}_G(\cS, I_0)\,\iso\, \cP_0\H^{\la}_H(\ast \gRes^{\kappa}(\cS), I_0)
\end{equation}
in $\D^b_{(Q_G\times Q_H)(\cO_0), mixed}(U\otimes L^*\otimes\Omega(F_0))$.

 We set $W_0=\Qlb^n$. Let $W_0=W_1\oplus W_2$ be the decomposition, where $W_1$ (resp., $W_2$) is generated by the first $m$ (resp., last $n-m$) base vectors. Equip $W_0\oplus W_0^*\oplus \Qlb$ with the symmetric form given by the matrix
$$
\left(
\begin{array}{ccc}
0 & E_n & 0\\
E_n & 0 & 0\\
0 & 0 & 1
\end{array}
\right),
$$
where $E_n\in\GL_n(\Qlb)$ is the unity.  Realize $\check{G}$ as $\SO(W_0\oplus W_0^*\oplus \Qlb)$.

 Equip the subspace $W_1\oplus W_1^*\subset W_0\oplus W_0^*\oplus\Qlb$ with the induced symmetric form, realize $\check{H}$ as $\SO(W_1\oplus W_1^*)$, this yields the inclusion $\check{H}\hook{}\check{G}$. Let $\sigma_H:\check{H}\to\check{H}$ be the automorphism sending $g$ to $^tg^{-1}$. It is understood that $\check{Q}_G=\Aut(W_0)$ and $\check{Q}_H=\Aut(W_1)$ canonically. Let $\kappa$ be the composition
$$
\check{H}\times\Gm\toup{\id\times \alpha_{\kappa}} \check{H}\times\GL(W_2)\hook{} \check{G},
$$
where $\alpha_{\kappa}:\Gm\to\GL(W_2)$ is $(n-m+1)(\check{\omega}_n-\check{\omega}_m)-2\check{\rho}_{\GL_{n-m}}$. The equality
$$
\alpha_{\kappa}-2(\check{\rho}_H-\check{\rho}_{Q_H})+n\check{\omega}_m=2(\check{\rho}_G-\check{\rho}_{Q_G})-2\rho_{\GL(W_2)}-m\check{\omega}_n
$$
shows that (\ref{diag_crucial_m_less_n}) commutes. If $m=n$ then $\kappa$ is trivial on $\Gm$. 

 Now by Proposition~\ref{Pp_functor_cP_injection}, 
(\ref{iso_for_Step2_Thmainlocal_m_less_n}) can be lifted to the desired isomorphism (\ref{iso_forThmainlocal_m_less_n}) in $\D\!\Weil^{ss}_{G,H}$. 

  To show that these isomorphisms are compatible with the tensor structures on $\Sph_G,\Sph_H$ we argue as follows. Let $\Weil^{ss, 0}_{G,H}$ be the full subcategory of $\Weil^{ss}_{G,H}$  whose objects are isomorphic to direct sums of objects $\IC(_{\lambda, N}\Upsilon^0)$ with $\lambda\in\Lambda^+_H$. Let $\D\Weil^{ss,0}_{G,H}$ be defined similarly as a full subcategory of $\D\Weil^{ss}_{G,H}$. 
  
  The functors in 1) of Theorem~\ref{Th_main_local_Sp_SO} take values in $\D\Weil^{ss,0}_{G,H}$. Indeed, the functor $\Sph_H\to \D\Weil^{ss}_{G,H}$, $\cT\mapsto \H^{\la}_H(\cT, I_0)$ takes values in $\D\Weil^{ss,0}_{G,H}$ by Lemma~\ref{Lm_very_nice_m_less_n}. Now from the already obtained isomorphism (\ref{iso_forThmainlocal_m_less_n}) we see that for $\cS\in \Sph_G$, $\H^{\la}_G(\cS, I_0)$ lies in $\D\Weil^{ss,0}_{G,H}$. We have got the commuting actions of $\Rep(\check{G})$ and of $\Rep(\check{H}\times\Gm)$ on $\D\Weil^{ss,0}_{G,H}$. By Remark~\ref{Rem_tensor_functors_appear}, we get a tensor functor $\epsilon: \Rep(\check{G})\to \Rep(\check{H}\times\Gm)$ such that the action of $\Rep(\check{G})$ on $\D\Weil^{ss,0}_{G,H}$ factors through $\epsilon$. 
  
   Equip $\D^b_{(Q_G\times Q_H)(\cO)}(U\otimes L^*\otimes\Omega(F))$ with the action of $\Rep(\check{G})$ via the restriction through 
$$
\kappa(\sigma_H\times\id)\kappa_H: \check{Q}_H\times\Gm\to \check{G}
$$ 
and the action of $\Rep(\check{Q}_H\times\Gm)$ by Hecke functors. We also equip the same category with the action of $\Rep(\check{H}\times\Gm)$ via the restriction through $\kappa_H: \check{H}\times\Gm\to \check{H}\times\Gm$. 
Then 
$$
J^*_{P_H}\comp f_H: \D\Weil^{ss,0}_{G,H}\to \D^b_{(Q_G\times Q_H)(\cO)}(U\otimes L^*\otimes\Omega(F))
$$ 
is canonically a morphism of $\Rep(\check{G})$-modules and also canonically a morphism of $\Rep(\check{H}\times\Gm)$-modules. So, the diagram of tensor functors commutes
$$
\begin{array}{cccccc}
\Rep(\check{G}) & \toup{\epsilon} &&&&\Rep(\check{H}\times\Gm)\\
& \searrow\lefteqn{\scriptstyle \Res^{\kappa(\sigma_H\times\id)\kappa_H}} &&&& \downarrow\lefteqn{\scriptstyle \Res^{\kappa_H}}\\
&&&&& \Rep(\check{Q}_H\times\Gm)
\end{array}
$$
This completes the proof.

\medskip\noindent
CASE $m>n$. Let
$$
\kappa_Q: \check{Q}_G\times\Gm\to \check{Q}_H\times\Gm
$$ 
be the map whose second component $\check{Q}_G\times\Gm\to \Gm$ is the projection, and the
first component $\check{Q}_G\times\Gm\to \check{Q}_H$ is 
the composition
$$
\check{Q}_G\times\Gm\toup{\id\times 2\check{\rho}_{\GL_{m-n}}} \check{Q}_G\times \GL_{m-n}\toup{\on{Levi}} \check{Q}_H
$$
 Write $\gRes^{\kappa_Q}:\D\Sph_{Q_H}\to\D\Sph_{Q_G}$ for the corresponding geometric restriction functor. 
 
 Now (\ref{iso_H_for_Th_local}) and Corollary~\ref{Cor_used_GL_n_GL_m} yield for $\cT\in\Sph_H$ canonical isomorphisms
$$
\cP_0 \H^{\la}_H(\cT, I_0)\,\iso\, \H^{\la}_{Q_H}(\gRes^{\kappa_H}(\cT), I_0)\,\iso\,
\H^{\la}_{Q_G}(\gRes^{\kappa_Q}(\ast \gRes^{\kappa_H}(\cT)), I_0)
$$
in $\D^b_{(Q_G\times Q_H)(\cO_0), mixed}(U\otimes L^*\otimes\Omega(F_0))$, here we assume that a decomposition $U_0\,\iso\,U_1\oplus U_2$ of vector spaces with $\dim U_1=n$ is fixed.
 On the other hand, (\ref{iso_G_for_Th_local}) yields a canonical isomorphism
$$
\cP_0(\H^{\la}_G(\gRes^{\kappa}(\ast\cT), I_0))\,\iso\, 
\H^{\la}_{Q_G}(\gRes^{\kappa_G}\gRes^{\kappa}(\ast\cT), I_0)
$$

 Let $\sigma$ be the automorphism of $\check{Q}_H=\GL_m$ sending $g$ to $^tg^{-1}$. The automorphism $\sigma\times\id$ of $\check{Q}_H\times\Gm$ induces the equivalence $\ast: \D\Sph_{Q_H}\,\iso\,\D\Sph_{Q_H}$. We will define an automorphism $\sigma_H$ of $\check{H}$ inducing $\ast: \Rep(\check{H})\,\iso\,\Rep(\check{H})$ and $\kappa$ making the following diagram commutative
\begin{equation}
\label{diag_crucial_n_less_m}
\begin{array}{ccccc}
\check{G}\times\Gm & \toup{\kappa} & \check{H} & \toup{\sigma_H} & \check{H}\\
\uparrow\lefteqn{\scriptstyle \kappa_G}&& && \uparrow\lefteqn{\scriptstyle \kappa_H}\\
\check{Q}_G\times\Gm & \toup{\kappa_Q} & \check{Q}_H\times\Gm & \toup{\sigma\times\id} & \check{Q}_H\times\Gm
\end{array}
\end{equation}
This will provide for $\cT\in\Sph_H$ a canonical isomorphism
\begin{equation}
\label{iso_for_Step2_Thmainlocal_n_less_m}
\cP_0 \H^{\la}_H(\cT, I_0)\,\iso\, \cP_0\H^{\la}_G(\gRes^{\kappa}(\ast\cT), I_0)
\end{equation}

 Let $W_0=\Qlb^m$, let $W_1$ (resp., $W_2$) be the subspace of $W_0$ spanned by the first $n$ (resp., last $m-n$) base vectors. Equip $W_0\oplus W_0^*$ with the symmetric form given by the matrix
$$ 
\left(
\begin{array}{cc}
0 & E_m\\
E_m & 0\\
\end{array}
\right),
$$
where $E_m\in\GL_m(\Qlb)$ is the unity. Realize $\check{H}$ as $\SO(W_0\oplus W_0^*)$. Let $\sigma_H$ be the automorphism of $\check{H}$ sending $g$ to $^tg^{-1}$. 
Let $\bar W\subset W_2\oplus W_2^*$ be the subspace spanned by $e_{n+1}+e_{n+1}^*$, equip $W_1\oplus W_1^*\oplus \bar W$ with the induced form and realize $\check{G}$ as $\SO(W_1\oplus W_1^*\oplus \bar W)$. Thus, the inclusion $i_{\kappa}:\check{G}\hook{}\check{H}$ is fixed. There is a unique $\alpha_{\kappa}: \Gm\to\check{H}$ such that for $\kappa=(i_{\kappa}, \alpha_{\kappa}):\check{G}\times\Gm\to \check{H}$ the diagram (\ref{diag_crucial_n_less_m}) commutes. Actually, $\alpha_{\kappa}=2\rho_{\GL(W_2)}+(m-n-1)(\check{\omega}_n-\check{\omega}_m)$. Note that if $m=n+1$ then $\alpha_{\kappa}$ is trivial. 
 
 By Proposition~\ref{Pp_functor_cP_injection},  (\ref{iso_for_Step2_Thmainlocal_n_less_m}) can be lifted to the desired isomorphism (\ref{iso_forThmainlocal_m_greater_n}) in $\D\Weil_{G,H}^{ss}$.  The compatibility with tensor structures on $\Sph_H, \Sph_G$ is obtained as in the case $m\le n$. 
\QED

\bigskip

\begin{Rem} 
\label{Rem_tensor_functors_appear}
Let $\check{G}_i$ be a reductive groups over $\Qlb$. Let $\cC$ be an abelian $\Qlb$-linear category with an exact $\Qlb$-linear commuting actions of $\Rep(\check{G}_i)$. Let $c\in\cC$ be such that the functor $\Rep(\check{G}_1)\to\cC$, $V\mapsto V\star c$ is an equivalence, here $\star$ denotes the action. Then there is a tensor functor $a: \Rep(\check{G}_2)\to \Rep(\check{G}_1)$ such that the action of $\Rep(\check{G}_2)$ on $\cC$ factors through $a$ and the given $\Rep(\check{G}_1)$-action on $\cC$.
\end{Rem}

\begin{Rem} 
\label{Rem_special_cases_n_m}
In the special case $m=1$ we have $H=Q_H$. So, in this case $\Weil_{G,H}$ is equivalent to the category $\P_{(H\times G)(\cO)}(\Pi(F))$, and one need not glue the categories as in Definition~\ref{Def_Weil_category}. The proof of Theorem~\ref{Th_main_local_Sp_SO} can be simplified in this case. 
For an integer $N$ let $_N\IC\in \P_{(H\times G)(\cO)}(\Pi(F))$ denote the constant perverse sheaf on $t^{-N}\Pi$. The irreducible objects of $\Weil_{G,H}$ in this case are exactly $_N\IC$, $N\in\ZZ$. For a dominant coweight $\lambda=N$ of $H$ in this case we get $\H^{\lambda}_H(I_0)\,\iso\, {_N\IC}$. 
\end{Rem}  

\bigskip\medskip

\centerline{\scshape 7. Global theta-lifting for the dual pair $\GL_n,\GL_m$}

\bigskip\noindent
In this section we prove Theorems~\ref{Th_main_global_GL_m_GL_n} and \ref{Th_global_Hecke_property_I_GL_n_GL_m}. 

\medskip\noindent
\begin{Prf}\select{of Theorem~\ref{Th_global_Hecke_property_I_GL_n_GL_m}}  

\smallskip\noindent
Recall the notation $U_0=k^m$, $L_0=k^n$, and the groups $G=\GL(L_0), H=\GL(U_0)$. Set $M_0=L_0\otimes U_0$ and $M=M_0(\cO)$ for $\cO=k[[t]]$. Viewing $M_0$ as a representation of $G\times H$, one defines the functor $\glob_{\infty}: \D_{(G\times H)(\cO)\rtimes \Aut^0(\cO)}(M(F))\to \D(_{\infty}\cW_{n,m})$ as in Section~4.6. 
One gets $\glob_{\infty}(I_0)\,\iso\, {_{\infty}\cI}$. The Hecke functors (\ref{Hecke_functor_H_forW_mn}) and (\ref{Hecke_functor_G_forW_mn}) are a particular case of those defined in Section~4.6. 
Since $\glob_{\infty}$ commutes with Hecke functors, our assertion follows from $\Aut^0(\cO)$-equivariant version of Proposition~\ref{Pp_action_Hecke_local_GL_m_GL_n} (cf. Remark~\ref{Rem_Aut_O_equivariance}).
\end{Prf}

\medskip
\begin{Prf}\select{of Theorem~\ref{Th_main_global_GL_m_GL_n}}

\smallskip\noindent
The argument below mimics that of (\cite{BG}, Section~4.1.8). To simplify notation, we will establish for $\cS\in\Sph_H$ and $K\in \D(\Bun_n)$ an isomorphism
\begin{equation}
\label{iso_for_Th5_tempo_two}
_x\H^{\la}_H(\cS, F_{n,m}(K))\,\iso\, F_{n,m}(_x\H^{\ra}_G(\Res^{\kappa}(\cS), K))
\end{equation}
for a given $x\in X$. The proof of the original statement is analogous. 

 Let $_{x,\infty}Z_H$ denote the stack classifying $(U,U', \beta: U'\,\iso\,U\mid_{X-x})\in{_x\cH_H}$, $L\in\Bun_n$ and $s:\cO_X\to L\otimes U'(\infty x)$. We have a diagram, where both squares are cartesian
$$
\begin{array}{ccccc}
_{x,\infty}\cW_{n,m} & \getsup{h^{\la}_{Z,H}} & {_{x,\infty}Z_H} & \toup{h^{\ra}_{Z,H}} & {_{x,\infty}\cW_{n,m}}\\
\downarrow\lefteqn{\scriptstyle h_m} && \downarrow && \downarrow\lefteqn{\scriptstyle h_m}\\
\Bun_m & \getsup{h^{\la}_H} & {_x\cH_H} & \toup{h^{\ra}_H} & \Bun_m
\end{array}
$$
Here $h^{\ra}_H$ (resp., $h^{\la}_H$) sends $(U,U')$ to $U'$ (resp., to $U$). The map $h^{\ra}_{Z,H}$ (resp., $h^{\la}_{Z,H}$) sends $(U,U',L, s:\cO\to L\otimes U'(\infty x))$ to $(L,U', s:\cO\to L\otimes U'(\infty x))$ (resp., to 
$(L,U, \beta\comp s:\cO\to L\otimes U(\infty x))$). 
 
  Using base change and the projection formula, one gets an isomorphism 
$$
_x\H^{\la}_H(\cS, F_{n,m}(K))\,\iso\, (h_m)_!(h_n^*K\otimes {_x\H^{\la}_H}(\cS, \cI))[-\dim\Bun_n]
$$  
Here $h_n: {_{x,\infty}\cW_{n,m}}\to \Bun_n$ is the corresponding projection. By Theorem~\ref{Th_global_Hecke_property_I_GL_n_GL_m}, this complex identifies with 
\begin{equation}
\label{complex_tempo_Th5}
(h_m)_!(h_n^*K\otimes {_x\H^{\la}_G}(\gRes^{\kappa}(\cS), \cI))[-\dim\Bun_n]
\end{equation}
 
 Let $_{x,\infty}Z_G$ be the stack classifying $(L,L',\beta: L'\,\iso\, L\mid_{X-x})\in{_x\cH_G}$, $U\in\Bun_m$ and $s:\cO_X\to L'\otimes U(\infty x)$. As above, one has the diagram
$$
\begin{array}{ccccc}
_{x,\infty}\cW_{n,m} & \getsup{h^{\la}_{Z,G}} & {_{x,\infty}Z_G} & \toup{h^{\ra}_{Z,G}} & {_{x,\infty}\cW_{n,m}}\\
\downarrow\lefteqn{\scriptstyle h_n} && \downarrow && \downarrow\lefteqn{\scriptstyle h_n}\\
\Bun_n & \getsup{h^{\la}_G} & {_x\cH_G} & \toup{h^{\ra}_G} & \Bun_n
\end{array}
$$
Here  $h^{\ra}_G$ (resp., $h^{\la}_G$) sends $(L,L')$ to $L'$ (resp., to $L$). The map $h^{\ra}_{Z,G}$ (resp., $h^{\la}_{Z,G}$) sends $(U,L',L, s:\cO\to L'\otimes U(\infty x))$ to $(L',U, s:\cO\to L'\otimes U(\infty x))$ (resp., to the collection $(L,U, \beta\comp s:\cO\to L\otimes U(\infty x))$). 
 
 The maps $h_m\comp h^{\la}_{Z,G}$ and $h_m\comp h^{\ra}_{Z,G}$ coincide. So, by base change and projection formula, (\ref{complex_tempo_Th5}) identifies with 
$$
F_{n,m}(_x\H^{\ra}_G(\gRes^{\kappa}(\cS), K))
$$ 
This yields the desired isomorphism (\ref{iso_for_Th5_tempo_two}). 
\end{Prf}

\bigskip\bigskip

\centerline{\scshape 8. Global theta-lifting for the dual pair $\SO_{2m},\Sp_{2n}$}

\bigskip\noindent
8.1 In this subsection we derive Theorem~\ref{Th_main_global_symplectic_orthogonal} from Theorem~\ref{Th_Hecke_property_Aut_GH}.  
 We give the argument for $m\le n$ (the case $m>n$ is completely similar). 
 
  By base change theorem, for $\cS\in \Sph_G$, $K\in\D(\Bun_H)$ we get
$$
\H^{\la}_G(\cS, F_G(K))\,\iso\, (\id\times\gp)_!(\gq^*K\otimes \H^{\la}_G(\cS, \Aut_{G,H}))[-\dim\Bun_H], 
$$
where $\id\times\gp: X\times\Bun_{G,H}\to X\times\Bun_G$ and $\gq: \Bun_{G,H}\to\Bun_H$ are the projections. By Theorem~\ref{Th_Hecke_property_Aut_GH}, the latter complex identifies with 
\begin{equation}
\label{complex_for_Sect_8_1}
(\id\times\gp)_!(\gq^*K\otimes \H^{\ra}_H(\gRes^{\kappa}(\cS), \Aut_{G,H}))[-\dim\Bun_H] 
\end{equation}
Now the diagram
$$
\begin{array}{ccccc}
X\times\Bun_H & \getsup{\supp\times h^{\la}_H} & \cH_H & \toup{h^{\ra}_H} & \Bun_H\\
\uparrow\lefteqn{\scriptstyle\id\times\gq} && \uparrow && \uparrow\lefteqn{\scriptstyle\gq}\\
X\times\Bun_{G,H} & \getsup{\supp\times h^{\la}_H} & \cH_H\times\Bun_G & \toup{\supp\times h^{\ra}_H} & X\times \Bun_{G,H}\\
&&&& \downarrow\lefteqn{\scriptstyle \id\times\gp}\\
&&&& X\times\Bun_G
\end{array}
$$
and the projection formula show that (\ref{complex_for_Sect_8_1}) identifies with
$$
(\id\times\gp)_!((\id\times\gq)^*\H^{\la}_H(\gRes^{\kappa}(\cS), K)\otimes \Aut_{G,H})[-\dim\Bun_H]
$$
This is what we had to prove. \QED

\medskip\noindent
8.2 In this section we we derive Theorem~\ref{Th_Hecke_property_Aut_GH} from Theorem~\ref{Th_main_local_Sp_SO}. To simplify notations, fix $x\in X$, we will establish isomorphisms (\ref{iso_Th_Hecke_Aut_GH_m_less_n}) and (\ref{iso_Th_Hecke_Aut_GH_m_grater_n}) over $x\times\Bun_{G,H}$. The fact that these isomorphisms depend on $x$ as expected is left to the reader (one uses Remark~\ref{Rem_Aut_O_equiv_for_SpSO} to check that the isomorphisms we obtaine are independent of a trivialization $\cO_x\,\iso\, k[[t]]$).

 Keep notations of Section~2. As in Section~6.1, let $L=\cO_x^n$, set $M=L\oplus L\otimes\Omega_x$ with the corresponding symplectic form $\wedge^2 M\to\Omega_x$. Let $U=\cO_x^m$, set $V=U\oplus U^*$ with the corresponding symmetric form $\Sym^2 V\to\cO_x$. Sometimes we view $M$ (resp., $V$) as the trivial $G$-torsor (resp., $H$-torsor) over $\Spec\cO_x$. 
 
  In view of Theorem~\ref{Th_main_local_Sp_SO},  Theorem~\ref{Th_Hecke_property_Aut_GH} is reduced to the following result, which we actually prove. 
  
\begin{Pp} 
\label{Pp_functor_LW}
There is a natural functor $\LW: \D\!\Weil_{G,H}^{ss,0}\to \D(\Bun_{G,H})$ commuting with the actions of both $\D\Sph_G$ and $\D\Sph_H$. There is an isomorphism $\LW(I_0)\,\iso\,\Aut_{G,H}$.
\end{Pp}
 
\medskip\noindent 
8.2.1 The proof is based on the following construction from \cite{L2}. Let $\cL_d(M\otimes V(F_x))$ denote the scheme of discrete lagrangian lattices in $M\otimes V(F_x)$. Let $\cA_d$ be the line bundle on $\cL_d(M\otimes V(F_x))$ with fibre $\det(M\otimes V:R)$ at $R\in \cL_d(M\otimes V(F_x))$ (cf. \select{loc.cit.} for the definition of this relative determinant). Note that $\cA_d$ is $(G\times H)(\cO_x)$-equivariant, so it can be viewed as a line bundle on the stack quotient 
$$
\cL_d(M\otimes V(F_x))/(G\times H)(\cO_x)
$$ 
Let $\wt\cL_d(M\otimes V(F_x))$ denote the $\mu_2$-gerbe of square roots of $\cA_d$. 
Write 
$
\wt\cL_d(M\otimes V(F_x))/(G\times H)(\cO_x)
$ 
for the corresponding $\mu_2$-gerbe over $\cL_d(M\otimes V(F_x))/(G\times H)(\cO_x)$.   
 
 Let $\Bun_{G,H}^x$ be the stack classifying $\cM\in\Bun_G, \cV\in\Bun_H$ and isomorphisms 
$$
\gamma_G: \cM\mid_{\Spec\cO_x}\,\iso\, M\mid_{\Spec\cO_x},\;\;\; \;\;\; \gamma_H: \cV\mid_{\Spec\cO_x}\,\iso\, V\mid_{\Spec\cO_x}
$$ 
of the corresponding $G$-torsors and $H$-torsors over $\Spec\cO_x$. One has a morphism of schemes
\begin{equation}
\label{map_0_xi_x}
_0\xi_x: \Bun_{G,H}^x\to \cL_d(M\otimes V(F_x))
\end{equation}
sending the above point of $\Bun_{G,H}^x$ to the image of 
$\H^0(X-x, \cM\otimes\cV)$ in $M\otimes V(F_x)$. 
 
 The group $(G\times H)(F_x)$ acts on $\Bun_{G,H}^x$ as follows. An element $g\in G(F_x)$ sends the above point of $\Bun_{G,H}^x$ to $(\cM', \gamma'_G, \cV, \gamma_H)$, where $\gamma'_G=g\gamma_G$ and $\cM'$ is the $\cO_X$-module whose sections over an open subset $U\subset X$ are $s\in\H^0(U-x, \cM)$ such that $g\gamma_G(s)\in M(\cO_x)$. The action of $H(F_x)$ is similar.
 
The morphism $_0\xi_x$ is equivariant with respect to natural actions of $(G\times H)(F_x)$. Taking the stack quotient by $(G\times H)(\cO_x)$, it yields a morphism of stacks
$$
\xi_x: \Bun_{G,H}\to \cL_d(M\otimes V(F_x))/(G\times H)(\cO_x)
$$

 We have canonically $\xi_x^*\cA_d\,\iso\, \tau^*\cA_{G_{2nm}}$, where $\tau$ is defined in Section~2.4.1. We lift $\xi_x$ to a morphism 
$$
\tilde \xi_x: \Bun_{G,H}\to \wt\cL_d(M\otimes V(F_x))/(G\times H)(\cO_x)
$$
sending $(\cM\in\Bun_G, \cV\in\Bun_H)$ to $(\xi_x(\cM,\cV), \cB)$, where 
\begin{equation}
\label{def_cB_for_xi}
\cB=\frac{\det\RG(X,\cV)^n\otimes\det\RG(X, \cM)^m}{\det\RG(X,\cO)^{2nm}}
\end{equation}
is equipped with an isomorphism $\cB^2\,\iso\, \det\RG(X, \cM\otimes\cV)$ given by Lemma~\ref{Lm_description_line_bundle_on_Bun_GH}. 

 For $r\ge 0$ write $_{rx}\Bun_{G,H}\subset\Bun_{G,H}$ for the open substack given by $\H^0(X, \cM\otimes\cV(-rx))=0$ for $\cM\in\Bun_G,\cV\in\Bun_H$. If $r'\ge r$ then $_{rx}\Bun_{G,H}\subset {_{r'x}\Bun_{G,H}}$ is an open substack, and we have 
$$
\lim_{r\in \NN^{op}} \D(_{rx}\Bun_{G,H})\,\iso\, \D(\Bun_{G,H})
$$ 

\smallskip\noindent
8.2.2 As in (\cite{L2}, Section~7.2) define the $\DG$-category $\D_{(G\times H)(\cO_x)}(\wt\cL_d(M\otimes V(F_x)))$ and the restriction functor 
\begin{equation}
\label{functor_tilde_xi_x_pull_back}
\tilde\xi_x^*: \D_{(G\times H)(\cO_x)}(\wt\cL_d(M\otimes V(F_x))) \to \D(\Bun_{G,H})
\end{equation}
as follows. For $N,r\in\ZZ$ with $N+r\ge 0$ and a free $\cO_x$-module $\cL$ of finite rank write $_{N,r}\cL=t^{-N}\cL/t^r\cL$. Let $\cL(_{N,N}M\otimes V)$ denote the scheme of lagrangian subspaces in the symplectic $k$-space $_{N,N}M\otimes V$. For $N\ge r\ge 0$ let 
$$
_r\cL(_{N,N}M\otimes V)\subset \cL(_{N,N}M\otimes V)
$$ 
be the open subscheme of $R\in \cL(_{N,N}M\otimes V)$ such that $R\cap {_{-r,N}(M\otimes V)}=0$. 
For $r_1\ge 2N$ let $\cA_N$ be the $\ZZ/2\ZZ$-graded line bundle on the stack quotient
$$
_r\cL(_{N,N}M\otimes V)/(G\times H)(\cO/t^{r_1})
$$
whose fibre at a lagrangian subspace $R$ is $\det(_{0,N}M\otimes V)\otimes \det R$. Write
$$
(_r\cL(_{N,N}M\otimes V)/(G\times H)(\cO_x/t^{r_1}))^{\tilde{}}
$$
for the gerbe of square roots of this line bundle. The $\DG$-categories of $\Qlb$-sheaves on these gerbes for all $r_1\ge 2N$ are canonically equivalent to each other (compatibly with the perverse t-structures) and are denoted
$$
\D_{(G\times H)(\cO_x)}(_r\wt\cL(_{N,N}M\otimes V))
$$
For $N_1\ge N\ge r\ge 0$ we have a projection 
$$
p: {_r\cL(_{N_1,N_1}M\otimes V)}\to {_r\cL(_{N,N}M\otimes V)}
$$ 
sending $R$ to $R\cap {_{N,N_1}M\otimes V}$. There is a canonical $\ZZ/2\ZZ$-graded isomorphism $p^*\cA_N\,\iso\,\cA_{N_1}$. For $r_1\ge 2N_1$ it yields a morphism of stacks
$$
\tilde p: (_r\cL(_{N_1,N_1}M\otimes V)/(G\times H)(\cO_x/t^{r_1}))^{\tilde{}}\to (_r\cL(_{N,N}M\otimes V)/(G\times H)(\cO_x/t^{r_1}))^{\tilde{}}
$$
The latter gives rise to the transition functor
\begin{equation}
\label{transition_functor_for_8.2.2}
\D_{(G\times H)(\cO_x)}(_r\wt\cL(_{N,N}M\otimes V))\to 
\D_{(G\times H)(\cO_x)}(_r\wt\cL(_{N_1,N_1}M\otimes V))
\end{equation}
sending $K$ to $\tilde p^*K[\dimrel(\tilde p)]$, it is exact for the perverse t-structures and a fully faithful embedding. The colimit in $\DGCat_{cont}$ 
$$
\mathop{\colim}\limits_{N\in\NN} \D_{(G\times H)(\cO_x)}(_r\wt\cL(_{N,N}M\otimes V))
$$ 
is denoted $\D_{(G\times H)(\cO_x)}(_r\wt\cL_d(M\otimes V(F_x)))$.
For $N\ge r'\ge r$ and $r_1\ge 2N$ we have an open immersion
$$
\tilde j: (_r\cL(_{N,N}M\otimes V)/(G\times H)(\cO_x/t^{r_1}))^{\tilde{}}\hook{} (_{r'}\cL(_{N,N}M\otimes V)/(G\times H)(\cO_x/t^{r_1}))^{\tilde{}}
$$
hence the restriction functors
$$
\tilde j^*: \D_{(G\times H)(\cO_x)}(_{r'}\wt\cL(_{N,N}M\otimes V))\to
\D_{(G\times H)(\cO_x)}(_r\wt\cL(_{N,N}M\otimes V))
$$
compatible with the transition functors (\ref{transition_functor_for_8.2.2}). Passing to the colimit over $N\in\NN$ we get the functors
$$
\tilde j_{r',r}: \D_{(G\times H)(\cO_x)}(_{r'}\wt\cL_d(M\otimes V(F_x)))\to
\D_{(G\times H)(\cO_x)}(_r\wt\cL_d(M\otimes V(F_x)))
$$
By definition, $\D_{(G\times H)(\cO_x)}(\wt\cL_d(M\otimes V(F_x)))$ is the limit in $\DGCat_{cont}$
$$
\lim_{r\in\NN^{op}}\D_{(G\times H)(\cO_x)}(_r\wt\cL_d(M\otimes V(F_x)))
$$ 
(cf. also \select{loc.cit.}). The category $\D(\wt\cL_d(M\otimes V(F_x)))$ is defined along the same lines.

The map $p$ fits into the diagram
$$
\begin{array}{ccc}
_{rx}\Bun_{G,H} & \toup{\xi_N} & {_r\cL(_{N,N}M\otimes V)/(G\times H)(\cO_x/t^{r_1})}\\
 & \searrow\lefteqn{\scriptstyle \xi_{N_1}} & \uparrow\lefteqn{\scriptstyle p}\\
&&{_r\cL(_{N_1,N_1}M\otimes V)/(G\times H)(\cO_x/t^{r_1})}
\end{array}
$$
where $\xi_N$ sends  $(\cM,\cV)$ to the lagrangian subspace $\H^0(X, \cM\otimes\cV(Nx))\subset {_{N,N}M\otimes V}$. As above, it is understood that one first picks a trivialization 
$$
\cM\otimes\cV\mid_{\Spec\cO_x/t^{r_1}}\,\iso\, M\otimes V\mid_{\Spec\cO_x/t^{r_1}}
$$
of the corresponding $G\times H$-torsor over $\Spec\cO_x/t^{2N}$ and further takes the stack quotients by $(G\times H)(\cO_x/t^{r_1})$.

We have canonically $\xi_N^*\cA_N\,\iso\, \tau^*\cA_{G_{2nm}}$.
So, we get a similar diagram between the gerbes
$$
\begin{array}{ccc}
_{rx}\Bun_{G,H} & \toup{\tilde\xi_N} & (_r\cL(_{N,N}M\otimes V)/(G\times H)(\cO_x/t^{r_1}))^{\tilde{}}\\
 & \searrow\lefteqn{\scriptstyle \tilde\xi_{N_1}} & \uparrow\lefteqn{\scriptstyle p}\\
&&(_r\cL(_{N_1,N_1}M\otimes V)/(G\times H)(\cO_x/t^{r_1}))^{\tilde{}}
\end{array}
$$
The functors $K\mapsto \tilde\xi_N^* K[\dimrel(\xi_N)]$ are compatible with the transition functors (\ref{transition_functor_for_8.2.2}) so yield a functor
$$
_r\xi_x^*: \D_{(G\times H)(\cO_x)}(_r\wt\cL_d(M\otimes V(F_x)))\to \D(_{rx}\Bun_{G,H})
$$
Passing to the limit by $r$, one gets the desired functor (\ref{functor_tilde_xi_x_pull_back}).  

\medskip\noindent
8.2.3 Recall that $\P(\Bun_{G,H})\,\iso\,\mathop{\lim}\limits_{r\in\NN^{op}} \P(_{rx}\Bun_{G,H})\subset \mathop{\lim}\limits_{r\in \NN^{op}} \D(_{rx}\Bun_{G,H})$ is a full subcategory. 
Let $S_{M\otimes V(F_x)}$ denote the theta-sheaf on $\wt\cL_d(M\otimes V(F_x))$ introduced in (\cite{L2}, Section~6.5). It is naturally $(G\times H)(\cO_x)$-equivariant, and we have $\tilde\xi_x^* S_{M\otimes V(F_x)}\,\iso\, \Aut_{G,H}$ by (\select{loc.cit.}, Theorem~3).

\medskip\noindent
8.2.4  As in (\cite{L2}, Section~5.4) let $\wt\Sp(M\otimes V)(F_x)$ denote the metaplectic group corresponding to the c-lattice $M\otimes V$ in $M\otimes V(F_x)$. This is a group stack classifying collections 
$$
(g\in \Sp(M\otimes V)(F_x), \cB, \; \cB^2\,\iso\, \det(M\otimes V: g(M\otimes V))),
$$
where $\cB$ is a 1-dimensional $k$-vector space. The product map sends 
$$
(g_1, \cB_1, \sigma_1: \cB_1^2\,\iso\, \det(M\otimes V: gM\otimes V)), (g_2, \cB_2, \sigma_2: \cB_2^2\,\iso\, \det(M\otimes V: gM\otimes V))  
$$
to $(g_1g_2, \cB, \sigma: \cB^2\,\iso\, \det(M\otimes V: g_1g_2M\otimes V))$, where $\cB=\cB_1\otimes\cB_2$ 
and $\sigma$ is the composition
\begin{multline*}
(\cB_1\otimes\cB_2)^2\toup{\sigma_1\otimes\sigma_2} \det(M\otimes V: g_1M\otimes V)\otimes \det(M\otimes V: g_2M\otimes V)
\toup{\id\otimes g_1} \\
\det(M\otimes V:g_1M\otimes V)\otimes\det(g_1M\otimes V: g_1g_2M\otimes V)\,
\iso\, \det(M\otimes V: g_1g_2M\otimes V)
\end{multline*}
  
\begin{Lm} 
\label{Lm_det_relative_for_SL}
Let $M_i$ be a free $\cO_x$-module of rank $n_i$. Then for $g\in \Gr_{\SL(M_0)}$ there is a canonical isomorphism of $\ZZ/2\ZZ$-graded lines 
$$
\det(M_0\otimes M_1: (gM_0)\otimes M_1)\,\iso\, \det(M_0: gM_0)^{n_1}
$$
\end{Lm}
\begin{Prf}
Let $\cA_0$ be the line bundle on $\Gr_{\SL(M_0)}$ with fibre $\det(M_0: gM_0)$ at $g\in \Gr_{\SL(M_0)}$. It is known that $\Pic\Gr_{\SL(M_0)}\,\iso\,\ZZ$ is generated by $\cA_0$. Let $\cL$ be the line bundle on $\Gr_{\SL(M_0)}$ with fibre 
$$
\det(M_0\otimes M_1: gM_0\otimes M_1)
$$ 
at $g$. A choice of a base in $M_1$ yields a $\ZZ/2\ZZ$-graded isomorphism $\cL\,\iso\, \cA_0^{n_1}$. Thus, the line bundle $\cL\otimes \cA_0^{-n_1}$ on $\Gr_{\SL(M_0)}$
is constant. Its fibre at $1\in\Gr_{\SL(M_0)}$ is canonically trivialized, so the line bundle itself is canonically trivialized.
\end{Prf}   
  
\bigskip

 By Lemma~\ref{Lm_det_relative_for_SL}, for $g\in G(F_x), h\in H(F_x)$ we have a canonical $\ZZ/2\ZZ$-graded isomorphism
\begin{multline*}
\det(M\otimes V: gM\otimes hV)\,\iso\, \det(M\otimes V: gM\otimes V)\otimes \det(gM\otimes V: gM\otimes hV)\,\iso\\ 
\det(M:gM)^{2m}\otimes \det(V: hV)^{2n}
\end{multline*} 
It yields a canonical section $(G\times H)(F_x)\to \wt\Sp(M\otimes V)(F_x)$ sending $(g\in G(F_x), h\in H(F_x)$ to $(g\otimes h, \cB, \,\cB^2\,\iso\,\det(M\otimes V: gM\otimes hV))$, where
$$
\cB=\det(M:gM)^{m}\otimes \det(V: hV)^{n}
$$
The canonical sections $(G\times H)(F_x)\to \wt\Sp(M\otimes V)(F_x)$ and $\Sp(M\otimes V)(\cO_x)\to \wt\Sp(M\otimes V)(F_x)$ are compatible over $(G\times H)(\cO_x)$. 

 The group $\wt\Sp(M\otimes V)(F_x)$ acts naturally on $\wt\cL_d(M\otimes V(F_x))$. Let $(G\times H)(F_x)$ act on $\wt\cL_d(M\otimes V(F_x))$ via the canonical section 
$(G\times H)(F_x)\to \wt\Sp(M\otimes V)(F_x)$. Let 
$$
_0\tilde\xi_x: \Bun_{G,H}^x\to \tilde\cL_d(M\otimes V(F_x))
$$
be the morphism sending $(\cM,\cV)$ to $(_0\xi_x(\cM, \cV),  \cB)$ with $\cB$ given by (\ref{def_cB_for_xi}). 
Then $_0\tilde\xi_x$ is naturally $(G\times H)(F_x)$-equivariant. For this reason, (\ref{functor_tilde_xi_x_pull_back}) commutes with the actions of $\D\Sph_G$ and $\D\Sph_H$ on both sides (cf. Proposition~\ref{Pp_11}).

\medskip\noindent
8.2.5 \select{Proof of Proposition~\ref{Pp_functor_LW}}

\medskip\noindent
Recall the notation $\Upsilon=L^*\otimes V\otimes\Omega_x$ and $\Pi=U^*\otimes M$ from Section~6.1. Consider the decomposition $M\otimes V=(L\otimes V)\oplus (L^*\otimes V\otimes\Omega_x)$. In (\cite{L2}, Definition~5) we associated to this decomposition a functor 
$$
\cF_{L\otimes V(F_x)}: \D(\Upsilon(F_x))\to \D(\wt\cL_d(M\otimes V(F_x)))
$$ 
exact for the perverse t-structures (its definition is also found in Appendix~A). Let also
$$
\cF_{M\otimes U(F_x)}: \D(\Pi(F_x))\to \D(\wt\cL_d(M\otimes V(F_x)))
$$
denote the corresponding functor for the decompositioon $M=(M\otimes U)\oplus (M\otimes U^*)$. By (\cite{L2}, Theorem~2 and Proposition~5), the diagram is canonically commutative
$$
\begin{array}{ccc}
\D(\Upsilon(F_x)) & \toup{\cF_{L\otimes V(F_x)}} & \D(\wt\cL_d(M\otimes V(F_x)))\\
\downarrow\lefteqn{\scriptstyle \zeta}
& \nearrow\lefteqn{\scriptstyle \cF_{M\otimes U(F_x)}}\\
\D(\Pi(F_x)),
\end{array}
$$
where $\zeta$ is the partial Fourier transform (\ref{zeta_geom}). 
 
 Let $(G\times H)(F_x)$ act on $\wt\cL_d(M\otimes V(F_x))$ via the canonical section 
$$
(G\times H)(F_x)\to \wt\Sp(M\otimes V)(F_x)
$$ 
Then $\cF_{L\otimes V(F_x)}$ commutes with the action of $H(F_x)$, and $\cF_{M\otimes U(F_x)}$ commutes with the action of $G(F_x)$ (cf. \cite{L2}, Section~6.6).  So, if $(\cF_1,\cF_2,\beta)$ is an object of $\Weil_{G,H}$ as in Definition~\ref{Def_Weil_category} then 
$\cF_{L\otimes V(F_x)}(\cF_1)$ is $H(\cO_x)$-equivariant, and $\cF_{M\otimes U(F_x)}(\cF_2)$ is $G(\cO_x)$-equivariant. We get a functor
$$
\LW_d: \Weil_{G,H}\to \P_{(G\times H)(\cO_x)}(\wt\cL_d(M\otimes V(F_x)))
$$
sending $(\cF_1,\cF_2,\beta)$ to $\cF_{L\otimes V(F_x)}(\cF_1)$. By definition, we get $\LW_d(I_0)\,\iso\, S_{M\otimes V(F_x)}$. Extend $\LW_d$ to a functor 
$$
\LW_d: \D\!\Weil_{G,H}\to \D_{(G\times H)(\cO_x)}(\wt\cL_d(M\otimes V(F_x)))
$$
by $\LW_d(K[r])=\LW_d(K)[r]$. By Proposition~\ref{Pp_appendix} from Appendix A, $\LW_d$ commutes with the actions of both $\D\Sph_G$ and $\D\Sph_H$. Finally, we set $\LW=\tilde\xi_x^*\comp\LW_d$. 
Our assertion follows. 

\smallskip

 Thus,  Proposition~\ref{Pp_functor_LW} and Theorem~\ref{Th_Hecke_property_Aut_GH} are proved. \QED

\bigskip\noindent
8.3 In this subsection we establish some additional properties of $\Aut_{G,H}$. Write $S_i$ for the stratum of $\Bun_{G,H}$ given by $\dim\H^0(X, M\otimes V)=i$ for $M\in\Bun_G, V\in\Bun_H$. 

\begin{Pp} i) The complex $\Aut_{G,H}$ is placed in perverse degrees $\ge 0$.\\
ii) For any $n,m$ the complex $\Aut_{G,H}$ has nontrivial perverse cohomologies in degrees larger than any given integer. 
\end{Pp}
\begin{Prf}
i) Recall that $\Bun_{P(G)}$ is the stack classifying $L\in\Bun_n$ and an exact sequence $0\to \Sym^2 L\to ?\to\Omega\to 0$ on $X$. Let $^0\Bun_{P(G)\times H}\subset \Bun_{P(G)\times H}$ be the open substack given by
\begin{equation}
\label{condition}
\H^0(X, \Sym^2 L)=0\;\;\;\; \mbox{and}\;\;\;\;\H^1(X, L^*\otimes\Omega\otimes V)=0
\end{equation}
for $V\in\Bun_H$, $(L\subset M)\in\Bun_{P(G)}$. 

Let $\cY_{P(G)}$ be the stack classifying $L\in\Bun_n$ and $s: \Sym^2 L\to\Omega^2$. So, $\cY_{P(G)}$ and $\Bun_{P(G)}$ are generalized dual vector bundles over $\Bun_n$.

 Denote by $^0(\cY_{P(G)}\times\Bun_H) \subset \cY_{P(G)}\times\Bun_H$ the open substack given by 
(\ref{condition}). We have the Fourier transform 
$$
\Four_{\psi}:\D(\cY_{P(G)}\times\Bun_H)\to\D(\Bun_{P(G)\times H})
$$
Write $\nu: \Bun_{P(G)\times H}\to\Bun_{G\times H}$ for the projection. Its restriction $^0\nu:{^0\Bun_{P(G)\times H}}\to\Bun_{G\times H}$ is smooth and surjective. 

 Let $\cV_{n,H}$ be the stack classifying $L\in\Bun_n$, $V\in\Bun_H$ and $v: L\to V\otimes\Omega$. Let $^0\cV_{n,H}\subset\cV_{n,H}$ be the open substack given by (\ref{condition}). We have an affine map
$\pi_{\cL}: \cV_{n,H}\to \cY_{P(G)}\times\Bun_H$  sending $v$ to the composition 
$$
\Sym^2 L\toup{v\otimes v} \Sym^2(V\otimes \Omega)\to\Omega^2
$$
Set 
$$
\cI=(\Qlb\mid_{\cV_{n,H}})[\dim\Bun_H+\dim\Bun_n+a_H],
$$
where $a_H$ is the function of a connected component of $\cV_{n,H}$ sending $(L,V,v)$ to $\chi(L^*\otimes\Omega\otimes V)$. By (\cite{L}, Proposition~1), there is an isomorphism over $\Bun_{P(G)\times H}$
\begin{equation}
\label{iso_for_v6}
\Four_{\psi}(\pi_{\cL})_!\cI\,\iso\, \nu^*\Aut_{G,H}[\dimrel(\nu)]
\end{equation}

 The complex $\cI$ is perverse over $^0\cV_{n,H}$ and coincides with $\Qlb[\dim(^0\cV_{n,H})]$. Since  $\pi_{\cL}$ is affine, $(\pi_{\cL})_!$ is left exact for the perverse t-structure. So, (\ref{iso_for_v6}) implies that 
$$
^0\nu^*\Aut_{G,H}[\dimrel(^0\nu)]
$$ 
is placed in perverse degrees $\ge 0$. 

\smallskip\noindent
ii) Consider two cases.

\smallskip\noindent
{\scshape CASE $2n>m-1$.} Pick a $k$-point $M$ of $\Bun_G$ and a rank $m$ vector bundle $E$ on $X$ of degree zero. Let $\cA$ be a line bundle on $X$ of degree $a>0$ large enough in the sense below. Set $U=E\otimes\cA$ and $V=U\oplus U^*$ with the symmetric form $\Sym^2 V\to\cO_X$ as in Section~3.2. 
The pair $(M,V)$ can be viewed as a closed substack
$Y:=B(\Aut(M)\times\Aut(V))$ of $\Bun_{G,H}$. 
We will let $a$ go to infinity, below we write $\const$ for quantities independent of $a$. 

 The dimension of the group of automorphisms $\Aut(M)$ is constant, whence $\dim\Aut(V)=m(m-1)a+\const$. Since $a$ is large enough, $Y\subset S_i$ for $i=2nma+\const$. Assuming that $\Aut_{G,H}$ is placed in perverse degrees $\le C$ for some $C>0$, we get
$$
\const-m(m-1)a=\dim Y\le \dim S_i\le \dim\Bun_{G,H}-i+C=\const-2nma
$$
Since $2n>m-1$, this is a contradiction.

\medskip
\noindent
{\scshape CASE $2m>n+1$.} Pick a $k$-point $V$ of $\Bun_H$ and a rank $n$ vector bundle $E$ on $X$ of degree zero. Let $\cA$ be a line bundle on $X$ of degree $a>0$ large enough. Set $L=E\otimes \cA$ and $M=L\oplus L^*\otimes\Omega$ with symplectic form $\wedge^2 M\to\Omega$ as in Section~6.1. As above, we get a closed substack $Y=B(\Aut(M)\times\Aut(V))$ of $\Bun_{G,H}$. Write $\const$ for quantities independent of $a$. 

 In this case $\dim\Aut(M)=n(n+1)a+\const$. For $a$ large enough we have $Y\subset S_i$ for $i=2nma+\const$. If $\Aut_{G,H}$ is placed in perverse degrees $\le C$ then 
$$
\const-n(n+1)a=\dim Y\le \dim S_i\le \dim\Bun_{G,H}-i+C=\const-2nma
$$ 
Since  $2m>n+1$, this is a contradiction. 

 Any pair of integers $n,m>0$ satisfies one of the above inequalities. We are done.
\end{Prf}

\medskip

\begin{Rem} The complex $\Aut_{G,H}$ on $\Bun_{G,H}$ is not pure in general. Let us show that for $m=1$ and any $n$ the complex $\Aut_{G,H}$ is not pure. 

 Recall that $\Bun_G$ is irreducible. We have $\Bun_H=\Pic X$. Consider the connected component 
$(\Pic^0 X)\times \Bun_G$, the intersection $\cU$ of this component with $S_0$ is nonempty. Indeed, 
take $V=\cO_X^2$, the corresponding point of $\Pic^0 X$ is $\cO_X$. There 
is $M\in \Bun_G$ with $\H^0(M)=0$.

 Over $\cU$ the complex $\Aut_{G,H}$ identifies with the constant perverse sheaf $\Qlb[\dim \cU]$. Its 
intermediate extension to $(\Pic^0 X)\times \Bun_G$ is the constant perverse sheaf. If $\Aut_{G,H}$ was 
pure, this constant perverse sheaf would be its direct summand. However, $(\Pic^0 
 X)\times \Bun_G$ contains points of $S_i$ for some $i>0$. The $*$-fibre of $\Aut_{G,H}$ at such point is 
the trivial one-dimensional space placed in usual degree $i-\dim \cU$. It can not contain $\Qlb[\dim \cU]$ as a direct summand.
\end{Rem}
 
\bigskip

\centerline{\bf Appendix A. Complement to \cite{L2}}

\bigskip\noindent
A.1  In this appendix independent of the rest of the paper we prove Proposition~\ref{Pp_appendix} below. We freely use the notations from \cite{L2}. 

 Let $\cO=k[[t]]\subset F=k((t))$, let $\Omega$ be the completed module of relative differentials of $\cO$ over $k$. Let $U$ be a free $\cO$-module of rank $d$, set $M=U\oplus U^*\otimes\Omega$. Equip $M$ with a natural symplectic form $\wedge^2 M\to \Omega$, so that $U$ and $U^*\otimes\Omega$ are lagrangian.

 For a c-lattice $R\subset R^{\perp}\subset M(F)$ remind that $\cA_R$ is the line bundle on $\cL(R^{\perp}/R)$ with fibre $\det L\otimes\det(M:R)$ at $L$. Write $\wt\cL(R^{\perp}/R)$ for the gerbe of square roots of $\cA_R$.

 Let $g\in G(F)$ be a $k$-point and $\tilde g=(g,\cB)\in \tilde G(F)$ with $\cB^2\,\iso\, \det(M:gM)$ its lifting to $\tilde G(F)$. One has a canonical isomorphism
\begin{equation}
\label{iso_tilde_g_for_Appendix}
\tilde g: \tilde\cL(R^{\perp}/R)\,\iso\, \tilde\cL(gR^{\perp}/gR) 
\end{equation}
defined as follows. It sends $(L, \cB_1)\in \tilde\cL(R^{\perp}/R)$ with $\cB_1^2\,\iso\, \det L\otimes \det(M:R)$ to $(gL, \cB\otimes\cB_1)$ equipped with 
\begin{multline*}
(\cB\otimes\cB_1)^2\,\iso\, \det L\otimes \det(M:R)\otimes\det(M: gM)\toup{g\otimes\id}\\
\det gL\otimes \det(gM:gR)\otimes\det(M: gM)\,\iso\, \det gL\otimes\det (M: gR)
\end{multline*}

 Let $H_R=(R^{\perp}/R)\times\A^1$ be the corresponding Heisenberg group. Then (\ref{iso_tilde_g_for_Appendix}) gives rise to the isomorphism
$$
\tilde g: \tilde\cL(R^{\perp}/R)\times \tilde\cL(R^{\perp}/R)\times H_R\,\iso\, \tilde\cL(gR^{\perp}/gR)\times \tilde\cL(gR^{\perp}/gR)\times H_{gR}
$$
for which one has canonically 
\begin{equation}
\label{iso_appendix_equivariance_G}
\tilde g^*F\,\iso\, F
\end{equation}
Here, by abuse of notation, for each finite-dimensional symplectic $k$-space $M_0$ we denote by $F$ the perverse sheaf on $\wt\cL(M_0)\times\wt\cL(M_0)\times H_{M_0}$ from (\select{loc.cit.}, Theorem~1). 
 
 For any $L^0_R, Q^0_R\in \wt\cL(R^{\perp}/R)$ the diagram is canonically commutative
$$
\begin{array}{ccc}
\D\cH_{Q_R} &\toup{\cF_{L^0_R, Q^0_R}} & \D\cH_{L_R}\\
\downarrow\lefteqn{\scriptstyle g} && \downarrow\lefteqn{\scriptstyle g}\\
\D\cH_{g(Q_R)} &\toup{\cF_{\tilde g (L^0_R), \tilde g (Q^0_R)}} & \D\cH_{g(L_R)}
\end{array}
$$ 
for the functors $\cF_{L^0_R, Q^0_R}$ from (\select{loc.cit.}, Section~6.2). Here the categories in the top (resp., low) row are categories of certain sheaves on $H_R$ (resp., on $H_{gR}$).

\medskip\noindent
A.2 We let $\GL(U)(F)$ act on $\wt\cL_d(M(F))$ via the section 
\begin{equation}
\label{section_for_GL(U)}
\GL(U)(F)\to \tilde G(F)
\end{equation}
defined in (\select{loc.cit.} just after Definition~5).

 Note that $U^*\otimes\Omega(F)$ is a placid ind-scheme, so we have the $\DG$-category $\D((U^*\otimes\Omega(F))$. For a $k$-point $L^0\in \wt\cL_d(M(F))$ in (\select{loc.cit.}, Proposition~5) we introduced the functor
$$
\cF_{U(F), L^0}: \D(U^*\otimes\Omega(F))\to \D\cH_L
$$  
Let us reformulate its definition. Given two c-lattices $T\subset S$ in $U(F)$, let $'T=\{x\in U^*\otimes\Omega(F)\mid \<x, u\>\in\Omega\;\mbox{for all}\; u\in T\}$, similarly for $'S$. Then $'S\subset {'T}$ are c-lattices in $U^*\otimes\Omega(F)$. Let $R=T\oplus {'S}$ then $R\subset R^{\perp}=S\oplus {'T}$, and $R^{\perp}/R$ is a symplectic space. 
  
  Set $U_R=S/T\in \cL(R^{\perp}/R)$. Let
$$
U_R^0=(U_R, \det(U: T))\in \wt\cL(R^{\perp}/R)
$$
equipped with a canonical $\ZZ/2\ZZ$-graded isomorphism
$\det(U: T)^2\,\iso\, \det U_R\otimes \det(M:R)$. One has a canonical isomorphism $\P('T/{'S})\,\iso\, \cH_{U_R}$ exact for the perverse t-structures. Then $\cF_{U(F), L^0}$ is the limit of the functors
$$
\D('T/{'S})\,\iso\,\D\cH_{U_R}\,\toup{\cF_{L^0_R, U^0_R}}\, \D\cH_{L_R}
$$
as $T$ becomes smaller and $S$ becomes bigger inside $U(F)$.
 
 Proposition~5 from \select{loc.cit.} can be strengthened as follows. 
\begin{Lm} 
\label{Lm_appendix}
Let $g\in\GL(U)(F)$. For a $k$-point $L^0\in \wt\cL_d(M(F))$ the diagram in $\DGCat_{cont}$ is canonically commutative
$$
\begin{array}{ccc}
\D(U^*\otimes\Omega(F)) & \toup{\cF_{U(F), L^0}} & \D\cH_L\\
\downarrow\lefteqn{\scriptstyle g} && \downarrow\lefteqn{\scriptstyle g}\\
\D(U^*\otimes\Omega(F)) & \toup{\cF_{U(F), \tilde gL^0}} & \D\cH_{gL},
\end{array}
$$
where $\tilde g$ is the image of $g$ under (\ref{section_for_GL(U)}). Here the left vertical arrow sends $K$ to $g_*K$.
\end{Lm}
\begin{Prf} Given c-lattices $T\subset S\subset U(F)$ set $R=T\oplus {'S}$. For the corresponding decomposition $gR=gT\oplus g{'S}$ one gets $g(U_R)=U_{gR}\in \cL(gR^{\perp}/gR)$. Moreover, $\tilde g(U^0_R)=U^0_{gR}$, where $\tilde g$ is the map (\ref{iso_tilde_g_for_Appendix}). As in Section~A.1, the following diagram is canonically commutative
$$
\begin{array}{ccccc}
\D('T/{'S}) & \iso &
\D\cH_{U_R} & \toup{\cF_{L^0_R, U^0_R}} & \D\cH_{L_R} \\
\downarrow\lefteqn{\scriptstyle g} && \downarrow\lefteqn{\scriptstyle g} && \downarrow\lefteqn{\scriptstyle g}\\
\D(g('T)/g('S)) & \iso &\D\cH_{g(U_R)} & \toup{\cF_{\tilde g(L^0_R), \tilde g (U^0_R)}} & \D\cH_{g(L_R)}
\end{array}
$$
Our assertion follows by going to the limit as $T$ becomes smaller and $S$ becomes bigger in $U(F)$.
\end{Prf}

  The diagram
$$
\begin{array}{ccc}
\cH_L & \toup{\cF_{L^0}} & \D(\wt\cL_d(M(F)))\\
\downarrow\lefteqn{\scriptstyle g} && \downarrow\lefteqn{\scriptstyle \tilde g} \\
\cH_{gL} & \toup{\cF_{\tilde g L^0}} & \D(\wt\cL_d(M(F))) 
\end{array}
$$
is canonically commutative, where the right vertical arrow sends $K$ to $(\tilde g^{-1})^*K$ (cf. \select{loc.cit.}, Section~6.3-6.4). Combining this fact with Lemma~\ref{Lm_appendix}, we get that the functor
$$
\cF_{U(F)}: \D(U^*\otimes\Omega(F))\to \D(\wt\cL_d(M(F)))
$$
commutes with the actions of $\GL(U)(F)$. (One may show that it actually commutes with the strong actions of $Shv(\GL(U)(F))$ in the sense of \cite{G5}, 1.4.3). 

 Let now $H\subset \GL(U)$ be a closed connected reductive subgroup. The line bundle $\cA_d$ on $\cL_d(M(F))$ is $H(\cO)$-equivariant, so we have a $\mu_2$-gerbe $\wt\cL_d(M(F))/H(\cO)$ of square roots of $\cA_d$ over the stack quotient $\cL_d(M(F))/H(\cO)$.
\begin{Pp} 
\label{Pp_appendix}
The functor 
$$
\cF_{U(F)}: \D_{H(\cO)}(U^*\otimes\Omega(F))\to \D_{H(\cO)}(\wt\cL_d(M(F)))
$$
commutes with the Hecke functors for $H$ acting on both sides.
\end{Pp}

 Let us first give a proof at the level of functions. Remind that if $R\subset M(F)$ is a c-lattice then 
$$
\wt\cL_d(M(F))_R\subset \wt\cL_d(M(F))
$$ 
is the open substack of $L^0\in \wt\cL_d(M(F))$ such that $L\cap R=0$. Let $T\subset S\subset U(F)$ be c-lattices, set $R=T\oplus {'S}$. By definition, if $K\in \D('T/{'S})$ then the restriction of $\cF_{U(F)}(K)$ to $\wt\cL_d(M(F))_R$ is given by
\begin{equation}
\label{value_dep_on_R}
L^0\mapsto \int_{y\in {'T}/{'S}} F_{L^0_R, U^0_R}((-y,0))K(y)dy
\end{equation}
for $L^0\in \wt\cL_d(M(F))_R$. Here $L^0_R$ is the image of $L^0$ under the natural map $\delta_R: \wt\cL_d(M(F))_R\to \wt\cL(R^{\perp}/R)$, and $(-y,0)\in (R^{\perp}/R)\times\A^1\,\iso\,H_R$. We denoted by $dy$ a suitable Haar measure. 
 
 Now let $\cT\in \Sph_H$ and $K\in \D_{H(\cO)}(U^*\otimes\Omega(F))$. Recall that, at the level of functions, $
\H^{\la}_H(\cT, K)$ is the function
$$
y\mapsto \int_{g\in H(F)/H(\cO)} K(g^{-1}y)\cT(g)dg,
$$ 
where $dg$ is a Haar measure on $\Gr_H$. 

 Let us check the equality of functions $\cF_{U(F)}\H^{\la}_H(\cT, K)=\H^{\la}_H(\cT, \cF_{U(F)}(K))$ at a given point $L^0\in \wt\cL_d(M(F))$. Pick $T$ sufficiently small and $S$ sufficiently large with respect to $\cT, K$ and such that $L^0\in \wt\cL_d(M(F))_R$. The value of 
$\cF_{U(F)}\H^{\la}_H(\cT, K)$ at $L^0$ is 
$$
\int_{y\in {'T}/{'S}, \; g\in H(F)} F_{L^0_R, U^0_R}((-y,0))K(g^{-1}y)\cT(g)dydg
$$
After the change of variables $g^{-1}y=z$, 
in view of (\ref{iso_appendix_equivariance_G})
the above expression becomes
\begin{multline*}
\int_{z\in g^{-1}('T)/g^{-1}('S), \; g\in H(F)} F_{L^0_R, U^0_R}((-gz,0))K(z)\cT(g)dzdg=\\
\int_{z\in g^{-1}('T)/g^{-1}('S), \; g\in H(F)} F_{\tilde g^{-1}(L^0_R), \tilde g^{-1}(U^0_R)}((-z,0))K(z)\cT(g)dzdg
\end{multline*}
Here $\tilde g$ is the image of $g$ under (\ref{section_for_GL(U)}).
Recall that $\tilde g^{-1}(U^0_R)=U^0_{g^{-1}R}$ for any $g\in\GL(U)(F)$. We have also $\tilde g^{-1}(L^0_R)=(\tilde g^{-1} L^0)_{g^{-1}R}$.

\medskip

  Since $g^{-1}L\in \cL_d(M(F))_{g^{-1}R}$, the value of $\H^{\la}_H(\cT, \cF_{U(F)}(K))$ at $L^0$ is
\begin{multline*}
\int_{g\in H(F)} (\cF_{U(F)}(K))(\tilde g^{-1} L^0)\cT(g)dg=\\
\int_{g\in H(F), z\in g^{-1}('T)/g^{-1}('S)} F_{(\tilde g^{-1} L^0)_{g^{-1}R}, U^0_{g^{-1}R}}((-z,0))K(z)\cT(g)dzdg
\end{multline*}
So, both expressions coincide. This completes the proof of Proposition~\ref{Pp_appendix} at the level of functions. 
 
\begin{Rem} 
\label{Rem_appendix}
Let $\bar T\subset T\subset S\subset \bar S$ be c-lattices in $U(F)$. Set $\bar R=\bar T\oplus {'\bar S}$ and $R=T\oplus {'S}$. The fact that (\ref{value_dep_on_R}) does not change if $R$ is replaced by $\bar R$ is a consequence of the following claim (obtained from \select{loc.cit.} Lemma~5 and Proposition~3). Let $V=\bar T\oplus {'S}$, so that $V^{\perp}=S\oplus {'\bar T}$. Let $H_{\bar R}=(\bar R^{\perp}/\bar R)\times \A^1$ and $H_V=(V^{\perp}/V)\times\A^1$
be the corresponding Heisenberg groups. Let $i_V: H^V=(V^{\perp}/{\bar R})\times\A^1\hook{} H_{\bar R}$, this is a closed subgroup. Let $\alpha_V: H^V\to H_V$ be the map $(u,a)\mapsto (u\!\mod V, a)$. Then
$$
\alpha_{V!}i_V^* F_{L^0_{\bar R}, U^0_{\bar R}}\,\iso\, F_{L^0_V, U^0_V}
$$
up to a shift. Further, let $H^R:=(R^{\perp}/V)\times\A^1\subset H_V$, let $\alpha_R: H^R\to H_R$ be the map sending $(u,a)$ to $(u\!\mod R, a)$. Write $i_R: H^R\hook{} H_V$ for the natural closed immersion. Then 
$$
\alpha_{R !}i_R^*F_{L^0_V, U^0_V}\,\iso\, F_{L^0_R, U^0_R}
$$
up to a shift (both claims are  also true in families as $L^0$ varies in $\wt\cL_d(M(F))_R$). Actually, the complex $i_R^*F_{L^0_V, U^0_V}$ is constant along the fibres of $\alpha_R$.
\end{Rem} 
\medskip\noindent
A.3 The precise definition of the Hecke functors 
$$
\H^{\la}_H: \Sph_H\times \D_{H(\cO)}(\wt\cL_d(M(F)))\to \D_{H(\cO)}(\wt\cL_d(M(F)))
$$
is left to a reader, let us only explain the idea.

 For $s_1+s_2\ge 0$ set $_{s_1,s_2}H(F)=\{g\in H(F)\mid t^{s_1}U\subset gU\subset t^{-s_2}U\}$. Let $_{s_1,s_2}\Gr_H={_{s_1,s_2}H(F)/H(\cO)}$. Let $K\in \D_{H(\cO)}(\wt\cL_d(M(F)))$. Assume that $\cT\in\Sph_H$ is the extension by zero from $_{s_1,s_1}\Gr_H$ for some $s_1>0$.
 
  For $r>0$ let $R=t^rM$ and $Z=t^{r+s_1}M$. One has a map 
$$
q_{\cL}: \tilde\cL_d(M(F))_R\times {_{s_1,s_1}H(F)} \to \wt\cL_d(M(F))_Z  
$$
sending $(L^0, g)$ to $\tilde g^{-1}L^0$. Set $\bar Y=t^{r+2s_1}M$. After taking the stack quotients by $H(\cO)$ on the LHS and by $H(\cO/t^s)$ on the RHS, the composition
$$
\tilde\cL_d(M(F))_R\times {_{s_1,s_1}H(F)} \toup{q_{\cL}} \wt\cL_d(M(F))_Z\toup{\delta_Z} \wt\cL(Z^{\perp}/Z)
$$ 
yields a map  
$$
\wt\cL_d(M(F))_R\times {_{s_1,s_1}\Gr_H}\to\wt\cL(Z^{\perp}/Z)/H(\cO/t^s)
$$
for $s$ large enough (actually $s\ge 2r+2s_1$). It factors as
$$
\wt\cL_d(M(F))_R\times {_{s_1,s_1}\Gr_H}\to \wt\cL(\bar Y^{\perp}/\bar Y)_R\times {_{s_1,s_1}\Gr_H}\toup{\tau_{\cL}} \wt\cL(Z^{\perp}/Z)/H(\cO/t^s)
$$
Here $\wt\cL(Z^{\perp}/Z)/H(\cO/t^s)$ denotes the corresponding $\mu_2$-gerbe over the stack $\cL(Z^{\perp}/Z)/H(\cO/t^s)$. 

 If $r$ is large enough then the restriction of $\H^{\la}_H(\cT, K)$ to $\wt\cL_d(M(F))_R$ is the inverse image under the natural map $\wt\cL_d(M(F))_R\to \wt\cL(\bar Y^{\perp}/\bar Y)_R$ of the complex
$$
\pr_{1!}(\cT\otimes \tau_{\cL}^*K)
$$
(up to a shift depending on a connected component of $_{s_1,s_1}\Gr_H$).
 
\medskip\noindent
A.4  The geometrization of the proof from Section~A.2 is quite formal, so we only give a sketch. Recall that for a free $\cO$-module $\cV$ and $N+r\ge 0$ we write $_{N,r}\cV=t^{-N}\cV/t^r\cV$. 

 Let $K\in \D_{H(\cO)}(_{N_1,r_1}U^*\otimes\Omega)$. Assume that $\cT\in\Sph_H$ is the extension by zero from $_{s_1,s_1}\Gr_H$.  
Let $T=t^rU$ and $S=t^{-r}U$, so that $'T=t^{-r}U^*\otimes\Omega$ and $'S=t^rU^*\otimes\Omega$. Set $R=T\oplus {'S}$. Assume that $r$ is large enough with respect to $K$, $\cT$. Set $Z=t^{s_1}R$.
Consider the diagram
$$
\begin{array}{cclc}
\tilde\cL_d(M(F))_R\times {_{s_1,s_1}H(F)} & \toup{q_{\cL}} & \wt\cL_d(M(F))_Z \\
\downarrow\lefteqn{\scriptstyle \delta_R\times\id}\\
\tilde\cL(R^{\perp}/R)\times {_{s_1,s_1}H(F)} & \getsup{\pr_{12}} & \tilde\cL(R^{\perp}/R)\times {_{s_1,s_1}H(F)}\times {_{r,r}U^*\otimes\Omega}& \toup{q} {_{r+s_1, r-s_1}U^*\otimes\Omega}
\end{array}
$$
Here $\pr_{12}$ is the projection, and $q$ sends $(L^0_R, g, y)$ to $g^{-1}y$. Let 
$$
\beta_R: \tilde\cL(R^{\perp}/R)\times {_{r,r}U^*\otimes\Omega}\to \tilde\cL(R^{\perp}/R)\times\tilde\cL(R^{\perp}/R)\times H_R
$$ 
be the map sending $(L^0_R, y)$ to $(L^0_R, U^0_R, (-y,0))$.

 From Remark~\ref{Rem_appendix} together with (\ref{iso_appendix_equivariance_G}) one gets the following. 
There is an isomorphism 
\begin{equation}
\label{iso_appendix_last}
 q_{\cL}^*(\cF_{U(F)}K)\,\iso\, (\delta_R\times\id)^*(\pr_{12})_!(\beta_R^*F\otimes q^*K)
\end{equation}
(up to a shift) geometrizing the equality
$$
(\cF_{U(F)}K)(\tilde g^{-1}L^0)=\int_{'T/{'S}} F_{L^0_R, U^0_R}((-y,0))K(g^{-1}y)dy
$$
for $L^0\in \tilde\cL_d(M(F))_R$, $g\in {_{s_1,s_1}H(F)}$.
 
 The above diagram gives rise, by taking the corresponding stack quotients, to the following one 
$$
\begin{array}{cclc}
\tilde\cL_d(M(F))_R\times {_{s_1,s_1}H(F)}\\
\downarrow\\
\wt\cL(\bar Y^{\perp}/\bar Y)_R\times {_{s_1,s_1}\Gr_H} & \toup{\tau_{\cL}} & \wt\cL(Z^{\perp}/Z)/H(\cO/t^s)\\
\downarrow\lefteqn{\scriptstyle \delta_{\bar Y,R}\times\id}\\
\tilde\cL(R^{\perp}/R)\times {_{s_1,s_1}\Gr_H} & \getsup{\pr_{12}} & \tilde\cL(R^{\perp}/R)\times {_{s_1,s_1}\Gr_H}\times {_{r,r}U^*\otimes\Omega}& \toup{\act_q} (_{r+s_1, r-s_1}U^*\otimes\Omega)/H(\cO/t^s)
\end{array},
$$
here  $\delta_{\bar Y, R}: \wt\cL(\bar Y^{\perp}/\bar Y)_R\to \tilde\cL(R^{\perp}/R)$ is the natural map. The isomorphism (\ref{iso_appendix_last}) descends to an isomorphism 
\begin{equation}
\label{iso_appendix_descended}
\tau_{\cL}^*(\cF_{U(F)}K)\,\iso\, (\delta_{\bar Y, R}\times\id)^*(\pr_{12})_!(\beta_R^*F\otimes \act_q^*K)
\end{equation}
over 
$\wt\cL(\bar Y^{\perp}/\bar Y)_R\times {_{s_1,s_1}\Gr_H}$.
  Tensor both sides of (\ref{iso_appendix_descended}) with $\cT$ and take the direct image under the projection 
$$\wt\cL(\bar Y^{\perp}/\bar Y)_R\times {_{s_1,s_1}\Gr_H}\to \wt\cL(\bar Y^{\perp}/\bar Y)_R
$$
One gets this way the desired isomorphism 
$$
\cF_{U(F)}\H^{\la}_H(\cT, K)\,\iso\,\H^{\la}_H(\cT, \cF_{U(F)}(K))
$$ 
over $\wt\cL_d(M(F))_R$. 
The other details of the proof of Proposition~\ref{Pp_appendix} are left to the reader.

\bigskip

\centerline{\bf Appendix B. Ind-pro systems and Hecke functors}

\medskip

\centerline{\bf B.1. Ind-schemes with an action of $G(F)$}

\bigskip
\noindent
B.1.1 On the referee's request we add this appendix where we generalize to some extent the ind-pro systems and Hecke  functors used in this paper. We will use the terminology from (\cite{K}, Section~4). 

 Let $\cI$ be a partially ordered set filtering in both directions, that is, for $\alpha_1,\alpha_2\in\cI$ there are $\alpha,\beta\in\cI$ such that $\beta\le \alpha_i\le\alpha$ for $i=1,2$. The limit by $\cI$ will be the limit as $\alpha$ becomes bigger in $\cI$. Let $\cI^0$ be $\cI$ with the reversed order. Set $I=\{(\alpha\in \cI, \beta\in \cI)\mid \beta\le\alpha\}$.
 
Assume given an ind-pro system $(X^{\alpha}_{\beta})$ indexed by pairs $(\alpha,\beta)\in I$. Assume that $X^{\alpha}_{\beta}$ is a $k$-scheme of finite type, and 
\begin{itemize}
\item[(1)] For each $\alpha\le\alpha'\in \cI$ the structure map $i^{\alpha\alpha'}_{\beta}: X^{\alpha}_{\beta}\to X^{\alpha'}_{\beta}$ is a closed embedding.
\item[(2)] For each $\beta\le\beta'\le\alpha\in \cI$ the structure map $p^{\alpha}_{\beta\beta'}: X^{\alpha}_{\beta'}\to X^{\alpha}_{\beta}$ is an affine morphism smooth of some relative dimension $d(\beta,\beta')$ independent of $\alpha$.
\item[(3)] For each $\beta\le\beta'\le \alpha\le\alpha'\in\cI$ the commutative square
\begin{equation}
\label{diag_X_beta_alpha}
\begin{array}{ccc}
X^{\alpha}_{\beta'} &\hook{} & X^{\alpha'}_{\beta'}\\
\downarrow && \downarrow\\
X^{\alpha}_{\beta} &\hook{} & X^{\alpha'}_{\beta} 
\end{array}
\end{equation}
is cartesian.
\item[(4)] there is $(\alpha, \beta)\in I$ such that $X^{\alpha}_{\beta}$ is smooth.
\end{itemize}

 By definition, $X^{\infty}_{\infty}=\varinjlim_{\alpha\in \cI}\varprojlim_{\beta\in \cI^0} X^{\alpha}_{\beta}$ is a placid ind-scheme. By (\select{loc.cit.}, Proposition~4.4.2), one has canonically
$$
X^{\infty}_{\infty}\,\iso\, \varprojlim_{\beta\in \cI^0}\varinjlim_{\alpha\in \cI} X^{\alpha}_{\beta}
$$
as prestacks. Set $X^{\infty}_{\beta}=\varinjlim_{\alpha\in\cI} X^{\alpha}_{\beta}$ and $X^{\alpha}_{\infty}=\varprojlim_{\beta\in \cI^0}X^{\alpha}_{\beta}$. Passing to limits, one gets the structure maps $i^{\alpha\alpha'}_{\infty}: X^{\alpha}_{\infty}\to X^{\alpha'}_{\infty}$ and $p^{\infty}_{\beta\beta'}: X^{\infty}_{\beta'}\to X^{\infty}_{\beta}$.

Assume in addition that each fibre of the map $p^{\alpha}_{\beta\beta'}$ is isomorphic to an affine space of dimension $d(\beta,\beta')$. One has the $\DG$-category $\D(X^{\infty}_{\infty})$ as in Section~2.1.2. 
  
 Namely, the functor $\D(X^{\alpha}_{\beta})\to \D(X^{\alpha}_{\beta'})$ given by $K\mapsto (p^{\alpha}_{\beta\beta'})^*K[d(\beta,\beta')]$ is exact for the perverse t-structures and fully faithful, and similarly for the functors $(i^{\alpha\alpha'}_{\beta})_!$. We let $\D(X^{\infty}_{\infty})$ be $\colim_{\alpha,\beta} \D(X^{\alpha}_{\beta})$, and similarly for the category $\P(X^{\alpha}_{\beta})$ of perverse sheaves.
 
  The function $d$ is initially defined on pairs $(\beta,\beta')\in\cI$. We extended it to a function $d: \cI\times\cI\to\ZZ$ as follows: if $\beta'\le\alpha$ and $\beta\le\alpha$ for some $\alpha\in\cI$ then $d(\beta,\beta')=d(\alpha,\beta')-d(\alpha,\beta)$.
 
\medskip\noindent 
B.1.2 Let now $G$ be a connected reductive group over $k$, set $\cO=k[[t]]\subset F=k((t))$. Assume that $G(\cO)$ acts (on the left) on each $X^{\alpha}_{\beta}$ via its finite-dimensional quotient $G(\cO/t^r)$ for some $r$ depending on $\alpha,\beta$. Assume that the maps $i^{\alpha\alpha'}_{\beta}$ and $p^{\alpha}_{\beta\beta'}$ are $G(\cO)$-equivariant.
Then $G(\cO)$ acts on each $X^{\alpha}_{\infty}$ by functoriality. The induced maps $X^{\alpha}_{\infty}\to X^{\alpha'}_{\infty}$ are $G(\cO)$-equivariant, and $G(\cO)$ acts on $X^{\infty}_{\infty}$ (cf. \select{loc.cit} and \cite{G4}, Appendix).

 One defines the category $\D_{G(\cO)}(X^{\infty}_{\infty})$ as follows. For $(\alpha,\beta)\in I$ the group $G(\cO)$ acts on each $X^{\alpha}_{\beta}$ via its finite-dimensional quotient $G(\cO/t^r)$, and for $r_1\ge r$ the projection between the stack quotients
$$
\phi: G(\cO/t^{r_1})\backslash X^{\alpha}_{\beta}\to G(\cO/t^r)\backslash X^{\alpha}_{\beta} 
$$
yields an exact for the perverse t-structures equivalence 
$$
\phi^*[\dimrel(\phi)]: \D_{G(\cO/t^r)}(X^{\alpha}_{\beta})\to \D_{G(\cO/t^{r_1})}(X^{\alpha}_{\beta})
$$
Denote by $\D_{G(\cO)}(X^{\alpha}_{\beta})$ the category $\D(G(\cO/t^{r_1})\backslash (X^{\alpha}_{\beta}))$ for any $r_1\ge r$. 

 For $\beta\le\beta'\in\cI^0$, $\alpha\le\alpha'\in \cI$ the diagram (\ref{diag_X_beta_alpha}) yields a diagram in $\DGCat_{cont}$
$$
\begin{array}{ccc}
\D_{G(\cO)}(X^{\alpha}_{\beta'}) & \to &  \D_{G(\cO)}(X^{\alpha'}_{\beta'})\\
\uparrow && \uparrow\\
\D_{G(\cO)}(X^{\alpha}_{\beta}) & \to & \D_{G(\cO)}(X^{\alpha'}_{\beta}),
\end{array}
$$ 
where each arrow is fully faithful and exact for the perverse t-structures functor. Then $\D_{G(\cO)}(X^{\infty}_{\infty})$ is $\colim_{\alpha,\beta}\; \D_{G(\cO)}(X^{\alpha}_{\beta})$ in $\DGCat_{cont}$.

 Along the same lines, one defines the category $\P_{G(\cO)}(X^{\infty}_{\infty})$ of perverse sheaves on $X^{\infty}_{\infty}$.
 
 Define also the $\DG$-category $\D_{G(\cO)}(X^{\infty}_{\infty}\times\Gr_G)$ as follows. Let $K\subset G(F)$ be a closed subscheme $G(\cO)$-invariant on the left and right and such that $K/G(\cO)$ is of finite type. Given $(\alpha,\beta)\in I$ let $r$ be large enough so that $G(\cO)$ acts on $X^{\alpha}_{\beta}\times K/G(\cO)$ via $G(\cO/t^r)$. For any $s\ge r$ one has a canonical equivalence exact for the perverse t-structures
$$
\D_{G(\cO/t^r)}(X^{\alpha}_{\beta}\times K/G(\cO))\iso
\D_{G(\cO/t^s)}(X^{\alpha}_{\beta}\times K/G(\cO))
$$
Define $\D_{G(\cO)}(X^{\alpha}_{\beta}\times K/G(\cO))$ as $\D_{G(\cO/t^s)}(X^{\alpha}_{\beta}\times K/G(\cO))$ for any $s\ge r$. Then  $\D_{G(\cO)}(X^{\infty}_{\infty}\times\Gr_G)$ is defined as $\colim_{\alpha,\beta, K} \;\D_{G(\cO)}(X^{\alpha}_{\beta}\times K/G(\cO))$ under $\alpha\in\cI, \beta\in\cI^0$ and $K$ becomming bigger and bigger. We similarly have the category $\P_{G(\cO)}(X^{\infty}_{\infty}\times\Gr_G)$ of perverse sheaves.
 
\medskip\noindent
B.1.3 Now assume in addition that each $X^{\alpha}_{\beta}$ is equidimensional then for $\alpha, \alpha'\in\cI$ one defines the relative dimension $\dim(X^{\alpha}_{\infty}: X^{\alpha'}_{\infty})$ as $\dim X^{\alpha}_{\beta}-\dim X^{\alpha'}_{\beta}$ for any $\beta\in\cI$ with $\beta\le\alpha, \beta\le \alpha'$.

Assume that the action of $G(\cO)$ on $X^{\infty}_{\infty}$ is extended to an action of $G(F)$ in the sense of prestacks.  Write $\act: X^{\infty}_{\infty}\times G(F)\to X^{\infty}_{\infty}$ for the map $(x,g)\mapsto g^{-1}x$.
We need the following well-known result.

\begin{Lm} 
\label{Lm_well-known_flattness_injlim}
If $A$ is a commutative ring and $A_i$ is an inductive system of flat $A$-algebras indexed by a filtering poset then $\varinjlim A_i$ is a flat $A$-algebra.\QED
\end{Lm}

Let $K\subset G(F)$ be a closed subscheme $G(\cO)$-invariant on the left and right and such that $K/G(\cO)$ is of finite type. For any such $K$ we assume the following.
\begin{itemize}
\item[(A)] For any $\alpha\in\cI$ there is $\alpha'\in \cI$ such that $\act: X^{\alpha}_{\infty}\times K\to X^{\infty}_{\infty}$ factors through $X^{\alpha'}_{\infty}$. For any $\beta'\le\alpha'\in\cI$ there is $\beta\le\alpha\in\cI$ and a commutative diagram
$$
\begin{array}{ccc}
X^{\alpha}_{\infty}\times K & \toup{\act} &X^{\alpha'}_{\infty}\\
\downarrow\lefteqn{\scriptstyle p_{\beta\infty}^{\alpha}\times\id} && \downarrow\lefteqn{\scriptstyle p_{\beta'\infty}^{\alpha'}}\\
X^{\alpha}_{\beta}\times K & \toup{\act_{\beta\beta'}} & X^{\alpha'}_{\beta'}
\end{array}
$$
(By Lemma~\ref{Lm_well-known_flattness_injlim}, $p^{\alpha}_{\beta\infty}$ is faithfully flat and quasi-compact. So, there may be at most one morphism $\act_{\beta\beta'}$ making the latter diagram commutative.
Though it is not reflected in the notation, the map $\act_{\beta\beta'}$ depends also on $K, \alpha,\alpha'$.)
 Besides, for any $\beta\le\alpha$ there is $\beta'\le \alpha'$ and a commutative diagram 
$$
\begin{array}{ccc}
X^{\alpha}_{\infty}\times K & \toup{\act} &X^{\alpha'}_{\infty}\\
\downarrow\lefteqn{\scriptstyle p^{\alpha}_{\beta\infty}} && \downarrow\lefteqn{\scriptstyle p^{\alpha'}_{\beta'\infty}} \\
X^{\alpha}_{\beta}\times K & \toup{\act_{\beta\beta'}} & X^{\alpha'}_{\beta'}
\end{array}
$$
\end{itemize}

\begin{Rem} 
\label{Rem_commutes_diag}
Let $\bar\beta\le\beta\le \alpha\in\cI$ and $\bar\beta'\le\beta'\le \alpha'\in\cI$, assume that $\act_{\beta\beta'}$ and $\act_{\bar\beta\bar\beta'}$ exist then the diagram commutes
\begin{equation}
\label{diag_beta_beta'_commutes}
\begin{array}{ccc}
X^{\alpha}_{\bar\beta}\times K & \toup{\act_{\bar\beta\bar\beta'}} & X^{\alpha'}_{\bar\beta'}\\
\downarrow && \downarrow\\
X^{\alpha}_{\beta}\times K & \toup{\act_{\beta\beta'}} & X^{\alpha'}_{\beta'}
\end{array}
\end{equation}
\end{Rem}

\medskip\noindent
B.1.4 Let 
$$
a: X^{\infty}_{\infty}\times G(F)\to X^{\infty}_{\infty}\times G(F)
$$ 
be the map $(m,g)\mapsto (g^{-1}m, g)$. Let $(a,b)\in G(\cO)\times G(\cO)$ act on the source sending $(m,g)$ to $(am, agb)$. Let it act on the target sending $(m',g')$ to $(b^{-1}m', ag'b)$. Then $a$ is equivariant for these actions, so yields a morphism of stack quotients
$$
_q\act: G(\cO)\backslash(X^{\infty}_{\infty}\times\Gr_G)\to (X^{\infty}_{\infty}/G(\cO))\times (G(\cO)\backslash \Gr_G),
$$
where the action of $G(\cO)$ on $X^{\infty}_{\infty}\times\Gr_G$ is the diagonal one.

 Assuming (A), in the rest of this subsection we define the inverse image functor
\begin{equation}
\label{inv_image_functor_under_A}
_q\act^*(\cdot,\cdot): \D_{G(\cO)}(X^{\infty}_{\infty})\times \D_{G(\cO)}(\Gr_G)\to \D_{G(\cO)}(X^{\infty}_{\infty}\times\Gr_G)
\end{equation}
satisfying the following properties.
\begin{itemize}
\item[A1)] For $S\in \D_{G(\cO)}(X^{\infty}_{\infty})^{constr}$ and $\cT\in \D_{G(\cO)}(\Gr_G)^{constr}$ one has
$\DD(_q\act^*(S,\cT))\,\iso\, {_q\act^*}(\DD(S),\DD(\cT))$ naturally.
\item[A2)] If both $S$, $\cT$ are perverse then $_q\act^*(S,\cT)$ is perverse.
\end{itemize}

 Let $K\subset G(F)$ be a left and right $G(\cO)$-invariant closed subscheme. The map $X^{\alpha}_{\beta}\times K \toup{\act_{\beta\beta'}} X^{\alpha'}_{\beta'}$ is $G(\cO)$-equivariant, where $h\in G(\cO)$ acts on $(x,g)\in X^{\alpha}_{\beta}\times K$ as $(x, gh)$ and on $y\in X^{\alpha'}_{\beta'}$ as $h^{-1}y$. Let $r$ be large enough so that $G(\cO)$ acts on $X^{\alpha'}_{\beta'}$ via $G(\cO/t^r)$. Then $\act_{\beta\beta'}$ induces a morphism of stack quotients
$$
X^{\alpha}_{\beta}\times K/G(\cO)\to G(\cO/t^r)\backslash X^{\alpha'}_{\beta'} 
$$ 
Assume in addition $r$ large enough so that $G(\cO)$ acts diagonally on $X^{\alpha}_{\beta}\times (K/G(\cO))$ via its quotient $G(\cO/t^r)$. The latter map being equivariant under this action, we get a diagram of stack quotients
$$
G(\cO/t^r)\backslash(K/G(\cO))\getsup{\pr_2}
G(\cO/t^r)\backslash(X^{\alpha}_{\beta}\times K/G(\cO))\toup{_r\act_{\beta\beta'}} G(\cO/t^r)\backslash X^{\alpha'}_{\beta'} 
$$
Consider the functor
$$
\begin{array}{ccc}
\D_{G(\cO/t^r)}(X^{\alpha'}_{\beta'})\times \D_{G(\cO/t^r)}(K/G(\cO)) & \toup{t_{\beta\beta'}} &\D_{G(\cO/t^r)}(X^{\alpha}_{\beta}\times K/G(\cO))\\
\parallel && \parallel\\
\D_{G(\cO)}(X^{\alpha'}_{\beta'})\times \D_{G(\cO)}(K/G(\cO)) & &\D_{G(\cO)}(X^{\alpha}_{\beta}\times K/G(\cO))
\end{array}
$$
sending $(S,\cT)$ to 
$$
(_r\act_{\beta\beta'})^*S\otimes \pr_2^*\cT[d(\beta',\beta)+\dim G(\cO/t^r)+\delta]
$$
here $\delta: \pi_1(G)\to\ZZ$ is some group homomorphism, we then view $\delta$ as a function of $K/G(\cO)$ sending $(K/G(\cO))\cap \Gr_G^{\theta}$ to $\delta(\theta)$. Here $\Gr_G^{\theta}$ is the connected component of $\Gr_G$ corresponding to $\theta\in\pi_1(G)$. We will precise $\delta$ later. 
  
  Now for any data as in Remark~\ref{Rem_commutes_diag} consider the commutative diagram (\ref{diag_beta_beta'_commutes}), it yields a commutative diagram
$$
\begin{array}{ccc}
\D_{G(\cO)}(X^{\alpha'}_{\bar\beta'})\times \D_{G(\cO)}(K/G(\cO)) &\toup{t_{\bar\beta\bar\beta'}} &\D_{G(\cO)}(X^{\alpha}_{\bar\beta}\times K/G(\cO))\\  
\uparrow && \uparrow  \\
\D_{G(\cO)}(X^{\alpha'}_{\beta'})\times \D_{G(\cO)}(K/G(\cO)) &\toup{t_{\beta\beta'}} &\D_{G(\cO)}(X^{\alpha}_{\beta}\times K/G(\cO)) , 
\end{array}
$$  
where vertical arrows are the transition functors.
Let $t'_{\beta\beta'}$ denote $t_{\beta\beta'}$ followed by the inclusion into $\D_{G(\cO)}(X^{\alpha}_{\infty}\times K/G(\cO))$.
Passing to the limit by $\beta'$, the functors $t'_{\beta\beta'}$ yield a functor
$$
t^{\alpha\alpha'}: \D_{G(\cO)}(X^{\alpha'}_{\infty})\times \D_{G(\cO)}(K/G(\cO)) \to\D_{G(\cO)}(X^{\alpha}_{\infty}\times K/G(\cO))
$$
Define a functor 
$$
t_K: \D_{G(\cO)}(X^{\infty}_{\infty})\times \D_{G(\cO)}(K/G(\cO)) \to\D_{G(\cO)}(X^{\infty}_{\infty}\times K/G(\cO))
$$ 
as follows. Let $S\in \D_{G(\cO)}(X^{\gamma}_{\infty})$ for some $\gamma\in\cI$, $\cT\in \D_{G(\cO)}(K/G(\cO))$. 
There are $\bar\alpha,\bar\alpha'$ sufficiently large in $\cI$ with the following properties:
\begin{itemize}
\item $\gamma\le\bar\alpha'$ in $\cI$, and there is a map $\act: X^{\bar\alpha}_{\infty}\times K\to X^{\bar\alpha'}_{\infty}$ as in (A).
\item for any $\alpha\ge\bar\alpha\in \cI$ and $\alpha'\ge \bar\alpha'\in\cI$ such that $\act: X^{\bar\alpha}_{\infty}\times K\to X^{\bar\alpha'}_{\infty}$ the images of $(S,\cT)$ under $t^{\alpha\alpha'}$ are canonically isomorphic in 
$\D_{G(\cO)}(X^{\infty}_{\infty}\times K/G(\cO))$.
\end{itemize}
This defines the desired functor $t_K$. The functors $t_K$ are compatible with the trasition functors, passing to the limit by $K$, one gets the desired functor (\ref{inv_image_functor_under_A}). One checks that there is a unique $\delta$ as above such that A2) holds for this functor.

\medskip\noindent 
B.1.5 Define an action of $\Sph_G$ on $\D_{G(\cO)}(X^{\infty}_{\infty})$ as follows. Let $K\subset G(F)$ be a left and right $G(\cO)$-invariant closed subscheme.
Assume that $r$ is large enough so that $G(\cO)$ acts on $X^{\alpha}_{\beta}\times K/G(\cO)$ via $G(\cO/t^r)$. For the projection 
$$
\pr: G(\cO/t^r)\backslash (X^{\alpha}_{\beta}\times K/G(\cO))\to G(\cO/t^r)\backslash X^{\alpha}_{\beta}
$$ 
the functors $\pr_!: \D_{G(\cO)}(X^{\alpha}_{\beta}\times K/G(\cO))\to \D_{G(\cO)}(X^{\alpha}_{\beta})$ are compatible with the transition functors, so yield a functor $\pr_!: \D_{G(\cO)}(X^{\infty}_{\infty}\times K/G(\cO))\to \D_{G(\cO)}(X^{\infty}_{\infty})$. Finally, we define the Hecke functor 
$$
\H^{\la}_G(\cdot,\cdot): \Sph_G\times \D_{G(\cO)}(X^{\infty}_{\infty})\to  \D_{G(\cO)}(X^{\infty}_{\infty})
$$
by $\H^{\la}_G(\cT, S)=\pr_!(_q\act^*(S,\cT))$ for $\cT\in\Sph_G, S\in  \D_{G(\cO)}(X^{\infty}_{\infty})$.

 By A2), they commute with the Verdier duality for constructible objects, namely
$$
\DD(\H^{\la}_G(\cT, S))\,\iso\, \H^{\la}_G(\DD\cT, \DD S)
$$
They are compatible with the tensor structure on $\Sph_G$ (as in \cite{BG}, Section~3.2.4).  For $\cT\in\Sph_G$ and $S\in \D_{G(\cO)}(X^{\infty}_{\infty})$ set $\H^{\ra}_G(\cT,S)=\H^{\la}_G(\ast \cT, S)$. Here $\ast: \Sph_G\to\Sph_G$ is the covariant equivalence of categories induced by by the map $G(F)\to G(F), g\mapsto g^{-1}$.
The functors 
$$
K\mapsto \H^{\la}_G(\cT,S)\;\;\;\mbox{and}\;\;\;
K\mapsto \H^{\ra}_G(\DD(\cT),S)
$$
are mutually (both left and right) adjoint. 

\bigskip

\centerline{\bf B.2. Ind-pro systems of scheme type}

\medskip\noindent
B.2.1 Now assume given an ind-pro system $(Y^{\alpha}_{\beta})$ indexed by $\cI$ as in Section~B.1.1 satisfying the properties (1)-(4) with the only difference that the maps $i^{\alpha\alpha'}_{\beta}$ are now open immersions.

 The for each $\alpha\le\alpha'\in\cI$, in the limit one gets an open immersion $i^{\alpha\alpha'}_{\infty}: Y^{\alpha}_{\infty}\hook{} Y^{\alpha'}_{\infty}$, and each $Y^{\alpha}_{\infty}$ is a scheme (not necessarily of finite type). Since $\cI$ is filtering, the colimit $Y^{\infty}_{\infty}=\colim_{\alpha\in\cI} Y^{\alpha}_{\infty}$ taken in $\PreStk$ is a scheme, which is a union of its open subschemes $Y^{\alpha}_{\infty}$ for $\alpha\in\cI$. We will say that $(Y^{\alpha}_{\beta})$ is an ind-pro system {\it of scheme type}.
 
  Assume in addition that each fibre of the map $p^{\alpha}_{\beta\beta'}$ is isomorphic to an affine space of dimension $d(\beta,\beta')$. 

 Under these assumptions, define the $\DG$-category $\D(Y^{\infty}_{\infty})$ as follows. By definition, $\D(Y^{\alpha}_{\infty})=\mathop{\colim}\limits_{\beta\in \cI^{op}}\D(Y^{\alpha}_{\beta})$ as in Section~B.1.1. Let $\D(Y^{\infty}_{\infty})=\mathop{\lim}\limits_{\alpha\in \cI^{op}}\D(Y^{\alpha}_{\infty})$ taken in $\DGCat_{cont}$. 

\medskip\noindent
B.2.2 Assume that the group $G(\cO)$ acts on each $Y^{\alpha}_{\beta}$ via its finite-dimensional quotient $G(\cO/t^r)$ for some $r$. Then as in Section~B.1.2 one defines the categories $\D_{G(\cO)}(Y^{\alpha}_{\beta})$, and 
$$
\D_{G(\cO)}(Y^{\alpha}_{\infty})\,\iso\, \colim_{\beta} \D_{G(\cO)}(Y^{\alpha}_{\beta})
$$  
Set $\D_{G(\cO)}(Y^{\infty}_{\infty})=\mathop{\lim}\limits_{\alpha\in\cI^{op}} \D_{G(\cO)}(Y^{\alpha}_{\infty})$.
 
 The category $\D_{G(\cO)}(Y^{\infty}_{\infty}\times\Gr_G)$ is defined as  
$\colim_{K}\; \D_{G(\cO)}(Y^{\infty}_{\infty}\times K/G(\cO))$ taken over $K$,  where $K/G(\cO)\subset\Gr_G$ is a closed $G(\cO)$-invariant subscheme of finite type (as $K$ becomes bigger and bigger). As in Section~B.1.2, 
$$
\D_{G(\cO)}(Y^{\alpha}_{\infty}\times K/G(\cO))\,\iso\,
\mathop{\colim}\limits_{\beta\in \cI} \D_{G(\cO)}(Y^{\alpha}_{\beta}\times K/G(\cO))
$$ 
Further, $\D_{G(\cO)}(Y^{\infty}_{\infty}\times K/G(\cO))\,\iso\,\mathop{\lim}\limits_{\alpha\in\cI^{op}} \D_{G(\cO)}(Y^{\alpha}_{\infty}\times K/G(\cO))$.

\medskip\noindent
B.2.3 Assume each $Y^{\alpha}_{\beta}$ equidimensional.
Assume that the action of $G(\cO)$ is extended to an action of $G(F)$ on $Y^{\infty}_{\infty}$ in the sense prestacks. Assume the condition (A). Then one defines a functor 
\begin{equation}
\label{q_act_def_for_scheme_type}
_q\act^*(\cdot,\cdot): \D_{G(\cO)}(Y^{\infty}_{\infty})\times \D_{G(\cO)}(\Gr_G)\to \D_{G(\cO)}(Y^{\infty}_{\infty}\times\Gr_G)
\end{equation} 
as follows.

 First, one defines the functors $t^{\alpha\alpha'}$ exactly as in Section~B.1.4. We view them as a functor
$$
t^{\alpha\alpha'}: \D_{G(\cO)}(Y^{\alpha'}_{\infty})\times \D_{G(\cO)}(K/G(\cO)) \to\D_{G(\cO)}(Y^{\alpha}_{\infty}\times (K/G(\cO)))
$$  
Further, if $\alpha\le\bar\alpha\in\cI$, $\alpha'\le\bar\alpha'\in\cI$ assume we have a diagram
$$
\begin{array}{ccc}
Y^{\bar\alpha}_{\infty}\times K & \toup{\act} & Y^{\bar\alpha'}_{\infty}\\
\uparrow\lefteqn{\scriptstyle i^{\alpha\bar\alpha}_{\infty}\times\id} && \uparrow\lefteqn{\scriptstyle i^{\alpha'\bar\alpha'}_{\infty}}\\
Y^{\alpha}_{\infty}\times K & \toup{\act} & Y^{\alpha'}_{\infty}
\end{array}
$$
Then $t^{\alpha\alpha'}$ commutes with the functors $(i^{\alpha\bar\alpha}_{\infty})^*$ for this diagram, so passing to the limit, one gets a functor
$$
t_K: \D_{G(\cO)}(Y^{\infty}_{\infty})\times \D_{G(\cO)}(K/G(\cO))\to \D_{G(\cO)}(Y^{\infty}_{\infty}\times (K/G(\cO)))
$$
Passing further to the colimit by $K$, one gets the desired functor (\ref{q_act_def_for_scheme_type}).

 If $\alpha\le\alpha'$ then for the cartesian square
$$
\begin{array}{ccc}
Y^{\alpha}_{\infty}\times (K/G(\cO))& \toup{i^{\alpha\alpha'}_{\infty}\times\id} & Y^{\alpha'}_{\infty}\times (K/G(\cO))\\
\downarrow\lefteqn{\scriptstyle\pr} &&\downarrow\lefteqn{\scriptstyle\pr}\\
Y^{\alpha}_{\infty} & \toup{i^{\alpha\alpha'}_{\infty}} & Y^{\alpha'}_{\infty}
\end{array}
$$
the functors $\pr_!: \D_{G(\cO)}(Y^{\alpha'}_{\infty}\times(K/G(\cO)) )\to \D_{G(\cO)}(Y^{\alpha'}_{\infty})$ commute with $(i^{\alpha\alpha'}_{\infty})^*$, so passing to the limit, one gets a functor
$$
\pr_!: \D_{G(\cO)}(Y^{\infty}_{\infty}\times (K/G(\cO))\to \D_{G(\cO)}(Y^{\infty}_{\infty})
$$
The Hecke functors 
$$
\H^{\la}_G(\cdot,\cdot): \Sph_G\times \D_{G(\cO)}(Y^{\infty}_{\infty})\to  \D_{G(\cO)}(Y^{\infty}_{\infty})
$$
for $\cT\in \Sph_G$, which is an extension by zero from $K/G(\cO)$, is then defined as
$$
\H^{\la}_G(\cT,S)=\pr_!(t_K(S,\cT))
$$

\medskip\noindent
B.2.4 Now assume we are given two ind-pro systems of scheme type $(Y^{\alpha}_{\beta})$ and $(\cY^{\alpha}_{\beta})$ with $\alpha,\beta\in\cI$ satisfying all the  assumptions of Section~B.2.

 Assume for each $(\alpha,\beta)\in\cI$ we are given a morphism $\pi^{\alpha}_{\beta}: Y^{\alpha}_{\beta}\to \cY^{\alpha}_{\beta}$ commuting with the actions of $G(\cO)$. Passing to the limit, we get a morphism $\pi^{\infty}_{\infty}: Y^{\infty}_{\infty}\to \cY^{\infty}_{\infty}$. 
 
Assume each $Y^{\alpha}_{\beta}$, $\cY^{\alpha}_{\beta}$ to be smooth of pure dimension, so we may consider the functors $(\pi^{\alpha}_{\beta})^*[\dimrel(\pi^{\alpha}_{\beta})]: \D_{G(\cO)}(\cY^{\alpha}_{\beta})\to \D_{G(\cO)}(Y^{\alpha}_{\beta})$. They are compatible with the transition functors, so yield a functor 
\begin{equation}
\label{functor_inv_image_for_app_B_Section2} 
(\pi^{\infty}_{\infty})^*: \D_{G(\cO)}(\cY^{\infty}_{\infty})\to \D_{G(\cO)}(Y^{\infty}_{\infty})
\end{equation}  
The proof of the following is left to a reader. 
\begin{Pp} 
\label{Pp_11}
Assume that $\pi$ is $G(F)$-equivariant. Then
(\ref{functor_inv_image_for_app_B_Section2}) commutes with the action of $\Sph_G$. \QED
\end{Pp}

\medskip

\centerline{\bf B.3. A functor given by a kernel}

\medskip\noindent
3.1 Assume given an ind-pro system $(X^{\alpha}_{\beta})$ satisfying all the assumptions of Section~B.1. Let $(Y^{\alpha}_{\beta})$ be an ind-pro system of scheme type satisfying all the assumptions of Section~B.2. They are indexed by the same set $(\alpha,\beta)\in I$.

 Assume for each $(\alpha,\beta)\in I$ we are given $K^{\alpha}_{\beta}\in \D_{G(\cO)}(X^{\alpha}_{\beta}\times Y^{\alpha}_{\beta})$. It gives rise to a functor
\begin{equation}
\cF^{\alpha}_{\beta}: \D_{G(\cO)}(X^{\alpha}_{\beta})\to \D_{G(\cO)}(Y^{\alpha}_{\beta})
\end{equation}
sending $S$ to $(\pr_2)_!(\pr_1^*S\otimes K^{\alpha}_{\beta})$
for the diagram of projections 
$$
G(\cO/t^r)\backslash X^{\alpha}_{\beta}\;\getsup{\pr_1} \; G(\cO/t^r)\backslash (X^{\alpha}_{\beta}\times Y^{\alpha}_{\beta})\;\toup{\pr_2}\; G(\cO/t^r)\backslash Y^{\alpha}_{\beta}
$$
for $r$ large enough (with respect to $\alpha,\beta$). The quotient here are stack quotients.
 
 In the only example we have the following happens. For $\alpha\le\alpha'\in\cI$ and $\beta'\le\beta\in\cI$ the diagram is canonically commutative
$$
\begin{array}{ccccc}
\D_{G(\cO)}(X^{\alpha}_{\beta}) & \toup{\cF^{\alpha}_{\beta}} & \D_{G(\cO)}(Y^{\alpha}_{\beta}) & \to & \D_{G(\cO)}(Y^{\alpha}_{\beta'})\\
\downarrow &&& \nearrow\lefteqn{\scriptstyle (i^{\alpha\alpha'}_{\beta'})^*}\\
\D_{G(\cO)}(X^{\alpha'}_{\beta'}) & \toup{\cF^{\alpha'}_{\beta'}} & \D_{G(\cO)}(Y^{\alpha'}_{\beta'}) 
\end{array}
$$
The two arrows which have no names in this diagram are the corresponding transition functors.

  The above diagram shows that in the limit by $\beta$ the functors $\cF^{\alpha}_{\beta}$ yield a functor 
$$
\cF^{\alpha}_{\infty}: \D_{G(\cO)}(X^{\alpha}_{\infty})\to \D_{G(\cO)}(Y^{\alpha}_{\infty})
$$  
 
  Define a functor $\cF^{\alpha}: \D_{G(\cO)}(X^{\alpha}_{\infty})\to \D_{G(\cO)}(Y^{\infty}_{\infty})$ as follows. Let $S\in \D_{G(\cO)}(X^{\alpha}_{\infty})$. For any $\alpha\le\alpha'\in\cI$ we declare the restriction of $\cF^{\alpha}(S)$ to $Y^{\alpha'}_{\infty}$ to be the image of $S$ under the composition 
$$
\D_{G(\cO)}(X^{\alpha}_{\infty})\toup{(i^{\alpha\alpha'}_{\infty})_!} \D_{G(\cO)}(X^{\alpha'}_{\infty})\toup{\cF^{\alpha'}_{\infty}} \D_{G(\cO)}(Y^{\alpha'}_{\infty})
$$
The corresponding projective system is naturally an object of $\D_{G(\cO)}(Y^{\infty}_{\infty})$.

 Passing to the colimit by $\alpha\in\cI$ the functors $\cF^{\alpha}$ yield a functor 
$$
\cF: \D_{G(\cO)}(X^{\alpha}_{\alpha})\to \D_{G(\cO)}(Y^{\infty}_{\infty})
$$
In the only example we have the functor $\cF$ commutes with the action of $G(F)$.
 
  An analog of Proposition~\ref{Pp_appendix} would be the claim that, possibly under some additional assumptions, the functor $\cF$ commutes with the action of the Hecke functors $\Sph_G$ on both sides.

\end{document}